\documentclass{article}

\usepackage[utf8]{inputenc}
\usepackage[T1]{fontenc}
\usepackage{geometry}
\usepackage{graphicx}
\usepackage{amsmath}
\usepackage{amssymb}
\usepackage{amsthm}
\usepackage{cases}
\usepackage{hyperref}
\usepackage{abstract}
\usepackage{float}
\usepackage{titling}
\usepackage[dvipsnames]{xcolor}
\usepackage[ruled,vlined]{algorithm2e}
\usepackage{enumerate}
\usepackage[
	safeinputenc,
	backend=biber,
	style=phys,
	biblabel=brackets,
	sorting=none,
	maxnames = 10,
	minnames = 1
  ]{biblatex}
\DeclareFieldFormat
  [article,inbook,incollection,inproceedings,patent,thesis,unpublished]
  {title}{#1\isdot}
\hypersetup{
    colorlinks=true,
    linkcolor=blue,
    filecolor=magenta,      
    urlcolor=cyan,
}
\newtheorem{Theorem}{Theorem}[section]
\newtheorem{Proposition}{Proposition}[section]
\newtheorem{Lemma}{Lemma}[section]

\newtheorem{Assumption}{Assumption}

\theoremstyle{definition}

\newtheorem{Definition}{Definition}[section]
\newtheorem{Remark}{Remark}[section]

\title{Sweeping Process Approach to Stress Analysis in
Elastoplastic Lattice Springs Models with Applications to Hyperuniform Network Materials
}

\date{}                     
\author{Ivan Gudoshnikov\thanks{Institute of Mathematics of the Czech Academy of Sciences, \v{Z}itn\'{a} 609/25, 115 67, Praha 1, Czech Republic, \url{gudoshnikov@math.cas.cz}} \and Yang Jiao\thanks{Materials Science and Engineering, Arizona State University,  Tempe, AZ 85287, United States, \url{yjiao13@asu.edu}} \and Oleg Makarenkov\thanks{Department of Mathematical Sciences, the University of Texas at Dallas, 800 West Campbell Road,
Richardson, TX 75080, United States, \url{makarenkov@utdallas.edu}} \and
Duyu Chen\thanks{Materials Research Laboratory, University of California, Santa Barbara, California 93106, United States, \url{duyu@alumni.princeton.edu}}}
\thanksmarkseries{arabic}
\begin{document}
\maketitle
\begin{abstract}

Disordered network materials abound in both nature and synthetic situations while rigorous analysis of their nonlinear mechanical behaviors remains challenging. The purpose of this paper is to connect the mathematical framework of sweeping process originally proposed by Moreau to a generic class of Lattice Spring Models with plasticity phenomenon. {\color{black} We derive the equations of quasistatic evolution of an elastic-perfectly plastic lattice and relate them to concepts from rigidity theory and structural mechanics. Then} we explicitly construct a sweeping process and provide numerical schemes to find the evolution of stresses in the model. In particular, we develop a highly efficient ``leapfrog'' computational framework that allow ones to rigorously track the progression of plastic events in the system based on the sweeping process theory. The utility of our framework is demonstrated by analyzing the elastoplastic stresses in a novel class of disordered network materials exhibiting the property of hyperuniformity, in which the infinite wave-length density fluctuations associated with the distribution of network nodes are completely suppressed. We find enhanced mechanical properties such as increasing stiffness, yield strength and tensile strength as the degree of hyperuniformity of the material system increases. Our results have implications for optimal network material design and our event-based framework can be readily generalized for nonlinear stress analysis of other heterogeneous material systems.

\end{abstract}

\section{Introduction}
\label{sect:introduction}

Disordered network materials such as collagen in extracellular matrix \cite{network01, network02, network03}, engineered cellular materials and foams \cite{network04, network05}, and certain amorphous 2D materials \cite{network06, network07, network08, network09}, abound in both nature and synthetic situations. Recent progress in advanced manufacturing {\color{black} such as} laser-based 3D printing allows salable production of a wide spectrum of complex network and cellular material systems, with desirable and optimized structural features. Microstructure-sensitive mechanical analysis of such materials, especially the non-linear elastoplastic behaviors, is crucial to establishing quantitative structure-property relations for material design and optimization.

Among the commonly used modeling frameworks, the Lattice Spring Models (LSM) represent the original material using an (ordered or disordered) network of springs, each possessing a nonlinear constitutive relation, which can naturally capture the complex geometrical and topological features of the material \cite{LSM01, LSM02, LSM03, LSM04, LSM05, Kale_OstojaStarzewski}. The preponderance of previous numerical solutions of Lattice Spring Models, especially when incorporating nonlinear spring models, typically employ a time-driven scheme with sufficiently small time steps in order to better capture the nonlinear behaviors (e.g., the transition and onset of plasticity, initialization of cracks etc.), which on the other hand, can be very computationally {\color{black} expensive}. In addition, even with very small time steps, there is no guarantee that all important plasticity events can be accurately captured. The purpose of this paper is to connect the Lattice Spring Models with plasticity phenomenon to the mathematical framework of sweeping process, which further enables us to devise a rigorous and efficient event-based ``leapfrog'' scheme for the elastoplastic stress analysis of complex disordered network materials.


The sweeping process is an important topic of contemporary research in mathematics of nonsmooth and nonlinear phenomena. Its purpose is to model the evolution of the processes with continuous time and firm one-sided (inequality) constraints on a state variable. A sweeping process can be described as a type of initial value problem governed by a time-depended (``moving'') convex set constraint, which ``sweeps'' a point (the state variable). The moving set as a function of time and the initial position of the point are the input data of the problem and the trajectory of the ``swept'' point is the solution. We will provide a short mathematical and visual introduction to the sweeping process in Section \ref{sect:sweeping_intro}.

The theory of sweeping process was founded by French mathematician and mechanics theorist J.-J. Moreau in early 1970's \cite{Moreau1973}. He employed it to describe nonsmooth phenomena in mechanics, such as elastoplasticity (e.g. one-dimensional continuous rod \cite{Moreau1976}), frictionless and frictional contact of rigid bodies \cite{Moreau1999, Moreau2004}.
Moreau's ideas are recognized as fundamental in contemporary literature on elastoplastic continuous media (e.g. \cite{han}) and nonsmooth mechanics (e.g. \cite{AcaryBrogliato2008,Adly2017,Brogliato2016}). 

In the recent decades the topic of sweeping process received exponentially increasing attention from researchers. One of the most important achievements in the field was the development of the theory of {\it optimal control} for the sweeping processes \cite{CHHM2012, CHHM2015, CHHM2016} which was later applied to robotics and traffic flow \cite{Colombo2019}, soft crawlers \cite{Colombo2021} and crowd motion \cite{CaoMordukhovich2019}. Optimal control of an elasto-plastic pseudo-rigid body is considered in \cite{CHHM2016} as a single-point toy model, and there is independent research available on the optimal control of elastoplastic continuous media \cite{Meyer,Meyer2} and an abstract rate-independent evolution variational inequality \cite{BrokateChristofPre}. Also, research on topology optimization based on quasi-static continuous elastoplasticity models \cite{Vermaak2021, AlmiStefanelli2022} became available recently. The sweeping process we construct in the current paper based on a discrete Lattice Springs Models is ready for future application of the optimal control theory to network-structured elastoplastic materials.

Another fruitful direction of research in sweeping processes is  the long-term asymptotic and stability analysis.  In particular, results offered by \cite{Krejci1996} and \cite{Adly2006} helped to establish the convergence of stresses to a periodic regime \cite{GM-1}, finite-time stability \cite{GMR} and structural stability \cite{GM-2} of periodic regimes in cyclically loaded rheological models of spatial dimension $1$. The present paper develops a framework that makes the results of \cite{GM-1,GM-2,GMR} applicable to rheological models of spatial dimensions higher than $1$.

The core of this paper is the construction of the sweeping process to model the stresses in the lattice, which can be summarized as the following. Let $m$ be the number of springs, then $\mathbb{R}^m$ represents all possible combinations of stresses in the lattice, while the set of stresses admissible by the elastic-perfectly plastic constitutive law is an $m$-dimensional rectangle in $\mathbb{R}^m$. On the other hand, the stresses which satisfy quasi-static equilibrium (the {\it self-stresses}) form a  hyperplane in $\mathbb{R}^m$. We use a change of variables, which converts external displacement and stress loads to parallel translation of the rectangle and then take the intersection of the (translated) rectangle with the hyperplane, associated with the self-stresses. The intersection is a polyhedron, which moves by parallel translation and changes its shape when, respectively, external displacement and stress loads vary. We construct the sweeping process with the intersection as its moving set, find the solution of the sweeping process, and recover the stress trajectory from the solution. This construction is provided in Section \ref{sect:sweeping_solution} with accompanying illustrations and examples.  

One can see from the above general explanation, that the properties of the graph structure of the lattice, such as the set of self-stresses, are as important in our construction as the elasto-plastic constitutive laws of individual springs. These properties directly influence qualitative and computational aspects of the problem, such as the external load we can impose and the overall dimension of the problem. In our previous works on the rheological models of spatial dimension $1$ \cite{GM-1,GM-2,GMR} it was enough to employ matrix graph theory (e.g. \cite{bapat}) to show that the properties of the corresponding sweeping process depend on the {\it cycle space} of the underlying graph and the amount of connected components in the graph. In this paper we consider lattices of spatial dimensions $2$ and higher, and the characteristics of their graph structure is a subject of {\it structural mechanics} and {\it rigidity theory}. The history of research in these areas goes back to James Clerk Maxwell \cite{Maxwell1864}, and they remain important topics in science even today due to their fundamental nature and abundant applications ranging from crystallography \cite{Giddy1993,Lubensky2015,Badri2014}, microstructures of metamaterials \cite{Meera2019} to sensor networks \cite{SensNets} and {\it tensegrity structures} \cite{RothWhiteley1981}, physically implemented in art and architecture \cite{TensegrityReview}. Alongside with the initial derivation of the equations of the Lattice Spring Model in Section \ref{sect:equations} we provide the related concepts from rigidity theory, which help us to rigorously explain {\color{black} various} aspects of our mathematical construction of the sweeping process, provide the motivation for the conditions we require and show how the dimension of the sweeping process depends on the graph structure of the lattice.

The paper is organized as following: after the introduction and preliminaries (Sections \ref{sect:introduction} and \ref{sect:preliminaries}) we establish the governing equations of the Lattice Spring Model (Section \ref{sect:equations}) accompanied by the concepts of structural mechanics and rigidity theory which have implications for our construction. In Section \ref{sect:sweeping_intro} we give a short presentation of the mathematical theory of the sweeping processes and describe the basic time-stepping numerical scheme, associated with the sweeping process, traditionally called the {\it catch-up algorithm}. Section \ref{sect:sweeping_solution} is a detailed guide on how to construct a sweeping process associated with the equations of the Lattice Spring Model, which is then used to compute the evolution of stresses via the adaptation of the catch-up algorithm. 
	
In Section \ref{sect:leapfrog} we discuss an event-based ``leapfrog'' numerical scheme, which can make the computation even more efficient in a special simple case of a sweeping process. In terms of the Lattice Spring Model this special case means that the stress load is constant and the displacement load changes at a constant rate. In particular, the event-based scheme allows to jump over the purely elastic phase of evolution in one step. The possibility to use an event-basted scheme under tighter regularity assumptions is common in simulations of nonsmooth systems, see e.g. the discussion in \cite{Moreau1999}. The utility of the event-based method is demonstrated by analyzing the stresses in the triangular grid with a defect (a hole) at its center, which is discussed in Section \ref{sect:example_with_the_hole}. Section \ref{sect:delaunay} is devoted to the analysis of elastoplastic stresses in a novel class of disordered hyperuniform network materials via the event-based scheme, which correspond to the Delaunay or stealthy hyperuniform point distributions with different degrees of disorder.  Section \ref{sect:conclusions} contains concluding remarks and Appendix \ref{sect:appendix_reduced_dimension} is devoted to more efficient versions of the time-stepping and event-based algorithms of reduced dimensions.

\section{Preliminaries}
\label{sect:preliminaries}
\subsection{Projection on a convex set and a normal cone}
Before we {\color{black} proceed to} the equations of Lattice Springs Model and the sweeping process, we would like to remind the reader some mathematical definitions and notations which we will rely on further down the text.

\begin{figure}[h]\center
\includegraphics[scale=0.5]{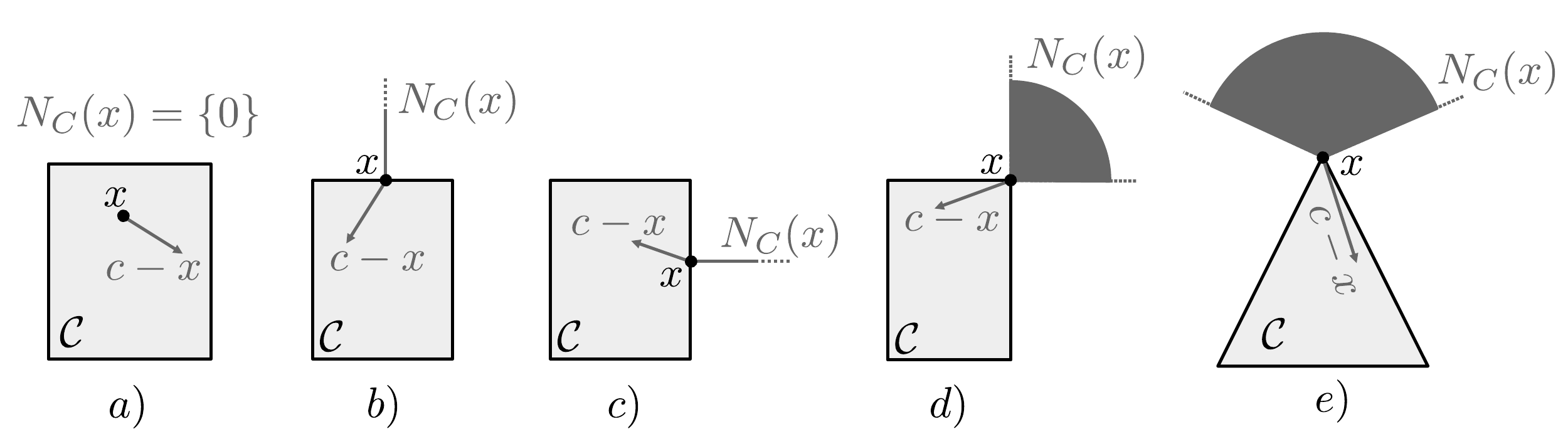}
\caption{\footnotesize a) when $x$ is in the interior of $\mathcal{C}$ the normal cone is always a singleton set of zero vector. b) and c) when $x$ lays on a boundary of $\mathcal{C}$ and the tangent to the boundary is well-defined, the normal cone is a one-parametric ray, orthogonal to the tangent and pointing outwards from the set. d) and e) when $x$ lays on a boundary which doesn't have a well-defined tangent, the normal cone fills an entire sector. } 
\label{fig:normal-cone-def}
\end{figure}

\begin{Definition} Let $\mathcal{C}$ be a nonempty closed convex set from $\mathbb{R}^n$. {\color{black} The {distance} from a point to the set is defined as
\begin{equation}
{\rm dist}(x, \mathcal{C}) = \underset{y\in \mathcal{C}}{{\rm min\,}}  \|x-y\|.
\label{eq:dist_def}
\end{equation}
In turn,} the {\it projection} of a point $x\in \mathbb{R}^n$ on a convex set $\mathcal{C}$ is the nearest point to $x$ among all the members of $\mathcal{C}$, i.e.
\begin{equation}
{\rm proj}(x,\mathcal{C}) = \underset{y\in \mathcal{C}}{{\rm arg\, min\,}}  \|x-y\|.
\label{eq:proj}
\end{equation}
\end{Definition}
\noindent The projection on a closed convex nonempty set always exists and is uniquely defined (see e.g. \cite[Section 3]{Hiriart-UrrutyANDLemarechal} or \cite[Th. 5.2]{Brezis2011}). \begin{Definition} Let $\mathcal{C}$ be a nonempty closed convex set from $\mathbb{R}^n$.  Given a point $x\in \mathcal{C}$, the {\it outward normal cone to $C$ at $x$} can be defined as a set of the vectors {\color{black} making an angle of at least $90^{\circ}$  (including zero vector)} with all the vectors of the type $c-x$, where $c\in \mathcal{C}$, i.e.
\begin{equation}
N_\mathcal{C} (x) =\{y\in \mathbb{R}^n:  y^T (c-x) \leqslant 0 \text{ for all }c\in \mathcal{C} \}.
\label{eq:normal_cone}
\end{equation}
\label{def:normal_cone}
\end{Definition}
\noindent Figure \ref{fig:normal-cone-def} illustrates Definition \ref{def:normal_cone} with several typical situations.
\noindent {\color{black} 
\begin{Remark}
We stress that the normal cone \eqref{fig:normal-cone-def} defined only for $x\in \mathcal{C}$, which is always assumed whenever notation $N_\mathcal{C}(x)$ is used.
\label{rem:hidden_inclusion}
\end{Remark}

\noindent The normal cone is a {\it convex cone}, i.e.
\[y,z\in N_\mathcal{C}(x) \Longrightarrow y+z\in N_\mathcal{C}(x),\]
\begin{equation}
y\in N_\mathcal{C}(x), \quad \lambda\geqslant 0 \Longrightarrow \lambda y \in N_\mathcal{C}(x).
\label{eq:normal_cone_is_a_cone}
\end{equation}
}
The projection on a convex set and the normal cone are related {\color{black} by }
\begin{equation}
y={\rm proj}(x,\mathcal{C}) \, \Longleftrightarrow \,  (x-y) \in N_{\mathcal{C}}(y).
\label{eq:nc_proj_relation}
\end{equation}

Given symmetric positive definite $n\times n$ matrix $S$ we can define a weighted inner product in $\mathbb{R}^n$:
\begin{equation}
x,y \mapsto x^T S y, \qquad \text{for all }x,y\in \mathbb{R}^n.
\label{eq:S-inner-prod}
\end{equation} 
{\color{black}
For a linear subspace $\mathcal{L}\subset \mathbb{R}^n$ we denote its orthogonal complement in sense of \eqref{eq:S-inner-prod} by $\mathcal{L}^{\perp_{S}}$, in the case of the standard inner product ($S=I_{n\times n}$) we write $\mathcal{L}^\perp$.
}
The definitions of distance \eqref{eq:dist_def} and projection \eqref{eq:proj} in sense of \eqref{eq:S-inner-prod} are modified as 
\begin{equation*}
{\color{black}
{\rm dist}^S(x,\mathcal{C}) = \underset{y\in \mathcal{C}}{{\rm min\,}}  \sqrt{(x-y)^T S(x-y)},
}
\end{equation*}
\begin{equation}
{\rm proj}^S(x,\mathcal{C}) = \underset{y\in \mathcal{C}}{{\rm arg\, min\,}}  (x-y)^T S(x-y).
\label{eq:S-proj}
\end{equation}
As the normal cone is defined via the inner product, in this case we will also use the notation
\begin{equation}
N^S_\mathcal{C} (x) =\{y\in \mathbb{R}^n:  y^T S(c-x) \leqslant 0 \text{ for all }c\in \mathcal{C} \}.
\label{eq:S-normal_cone}
\end{equation}   
In this text we will deal with a special case when set $\mathcal{C}$ is {\it polyhedral}, i.e. it can be written as
\begin{equation}
\mathcal{C}=\{x\in\mathbb{R}^n:Ax\leqslant b, A_{eq}x=b_{eq}\}
\label{eq:polyhedral_set}
\end{equation}
(here and for the rest of the paper vector inequalities are meant in the component-wise sense),
where $A, A_{eq}$ are fixed matrices of dimensions, respectively, $l\times n$ and $l_{eq}\times n$ for some $l, l_{eq}\in \mathbb{N}$, and $b$, $b_{eq}$ are vectors from $\mathbb{R}^l, \mathbb{R}^{l_{eq}}$ respectively. For example, in Fig. \ref{fig:normal-cone-def} c), d) set $\mathcal{C}$ is polyhedral. For a point $x\in \mathcal{C}$ we say that $i$-th constraint is {\it active} if and only if the inequality in \eqref{eq:polyhedral_set} is satisfied as an equality for $i$-th component, i.e.  $(Ax)_i=b_i$.  The projection \eqref{eq:S-proj} onto a polyhedral set takes the form
\begin{equation}
{\rm proj}^S (x, \mathcal{C})=\underset{\substack{y\in \mathbb{R}^n:\\  A y\leqslant b, \\ A_{eq} y = b_{eq}}}{\rm arg\, min} \frac{1}{2} y^TS y + (-Sx)^T y,  
\label{eq:S-proj-polyhedral}
\end{equation} 
which is  in the standard form of the quadratic programming problem, with well-developed numerical methods available in numerical  libraries.

\subsection{Moore-Penrose pseudoinverse matrix}
In this text we will use (real) Moore-Penrose pseudoinverse matrix as it is an important practical tool to solve {\color{black} linear algebraic} equations, and it is also available in many numerical packages. 
Here we remind the reader the definition of the Moore-Penrose pseudoinverse and some of its basic properties.
\begin{Proposition}{\bf \cite[p. 9]{bapat}}
Let $A$ be an $m \times n$-matrix. Then there exists a unique $n \times m$-matrix $A^+$, called {\it Moore-Penrose pseudoinverse of $A$}, such that all of the following hold:
\begin{align}
AA^+A&=A, \label{eq:MP1}\\
A^+AA^+&=A^+, \label{eq:MP2}\\
(AA^+)^T&=AA^+,\label{eq:MP3}\\
(A^+A)^T&=A^+A.\label{eq:MP4}
\end{align}
\label{prop:MP1}
\end{Proposition}
\begin{Proposition}{\bf \cite[Def. 1.1.2]{Campbell2008}}
A matrix $A^+$ is a Moore-Penrose pseudoinverse of $A$ if and only if $AA^+$ and $A^+A$ are orthogonal projection matrices onto, respectively, ${\rm Im}\,A$ and ${\rm Im}\,A^T$.
\label{prop:MP2}
\end{Proposition}
\noindent The following proposition can be verified by direct substitution into \eqref{eq:MP1}-\eqref{eq:MP4}:
\begin{Proposition}
If the columns of $A$ are linearly independent, then $A^+= (A^TA)^{-1}A^T$ and $A^+$ is a {\it left inverse} of $A$, i.e.
\begin{equation*}
A^+A=I_{n\times n}.
\end{equation*}
Similarly, if the rows of $A$ are linearly independent, then $A^+= A^T(AA^T)^{-1}$ and $A^+$ is a {\it right inverse} of $A$, i.e.
\begin{equation*}
AA^+=I_{m\times m}.
\end{equation*}
\label{prop:MP3}
\end{Proposition}

{\color{black}
\subsection{Directed graph and incidence matrix}
Consider a set $V$ of $n$ elements called {\it nodes (vertices)} and a set $E\subset V\times V$ of $m$ {\it edges}, which are ordered pairs $(v_1,v_2)$. Combined, $V$ and $E$ define a mathematical structure of a {\it directed graph}.  For a given edge $(v_1,v_2)$ node $v_1$ is called the {\it origin}, node $v_2$ is called the {\it terminus} \cite[Ch. 7]{ClarkHoltonGraphTheory}, and both $v_1$ and $v_2$ are called {\it endpoints} of the edge. Any directed graph can be described by an $n\times m$ matrix $Q$ called {\it incidence matrix}, provided that the origin and the terminus are distinct for each edge. The incidence matrix is constructed as the following \cite{bapat}:
for $i\in \overline{1,m},\, j\in\overline{1,n}$ set
\begin{equation*}
Q_{ji}= \begin{cases}0 &\text{ if node $j$ is not an endpoint of edge $i$,}\\
1 & \text{ if node $j$ is the origin of edge $i$,} \\
-1& \text{ if node $j$ is the terminus of edge $i$.}
\end{cases}
\end{equation*}
}

\section{Equations of the Lattice Spring Model}
\label{sect:equations}
\subsection{Geometry and linearized kinematics of lattices}
\label{ssect:geometry_of_the_lattice}
A Lattice Spring Model is given as a graph with $n$ nodes (vertices) and $m$ edges, where each edge is an elastic-perfectly plastic spring. The vertices are said to be from $\mathbb{R}^d$ representing a physical space, so we {\color{black} focus on} $d=1$, $d=2$ and $d=3$. 
	 
{\color{black}
In the model we assume that the graph structure of the lattice is given by an incidence matrix $Q$ of a directed graph, obtained by assigning an arbitrary orientation to each spring. The assigned orientations only serve accounting purposes and the model does not depend on their choice, as it will be evident from formulas below.
}

At any particular moment the coordinates of the vertices can be collected into a vector $\xi\in \mathbb{R}^{nd}$ so that $\xi_{d(j-1)+k}$ is $k$-th coordinate of node $j$ (where $j\in \overline{1,n}, k\in \overline{1,d}$). {\color{black} Let the origin and the terminus of spring $i$  (where $i\in \overline{1,m}$) be, respectively, $j'$ and $j''$, then the length of the spring is the norm of the vector
\[
\left(\xi_{d(j''-1)+k}\right)_{k\in \overline{1, d}} - \left(\xi_{d(j'-1)+k}\right)_{k\in \overline{1, d}}=-\left(\sum_{j=1}^n Q_{ji}\xi_{d(j-1)+k}\right)_{k\in \overline{1, d}}\]
}
\noindent The lengths of all $m$ springs can be collected in a vector from $\mathbb{R}^m$ {\color{black} and expressed as a value of the following function of $\xi$:}
\[
 \varphi:\mathbb{R}^{nd}\to \mathbb{R}^m
\]
\begin{equation}
\varphi(\xi)=\left(\varphi_i(\xi)\right)_{i\in\overline{1,m}}:=\left(\sqrt{\sum_{k=1}^d\left(\sum_{j=1}^n Q_{ji}\xi_{d(j-1)+k}\right)^2}\right)_{i\in\overline{1,m}}.
\label{eq:distanceFunction}
\end{equation}

We choose a reference configuration of nodes $\xi_0=(\xi^0_{d(j-1)+k})\in \mathbb{R}^{nd}$ and use the linearization of \eqref{eq:distanceFunction} at the reference configuration to write the first governing equation
\begin{equation}
\left(D_{\xi_0} \varphi\right) \zeta=x \tag{LSM1}
\label{eq:gc}
\end{equation}
in which $\zeta \in\mathbb{R}^{nd}$ is the vector of {\it displacements} of the nodes from the reference configuration $\xi_0$, $x\in\mathbb{R}^m$ is the vector of {\it total elongations} of the springs from the lengths $\varphi(\xi_0)$ and $D_{\xi_0} \varphi$ is the $m\times (nd)$ Jacobi matrix of  $\varphi$ at $\xi_0$. Specifically, the $(i, d(j-1)+k)$ entry of $D_{\xi_0} \varphi$ is given by
\begin{equation}
\left. \frac{\partial \varphi_i}{\partial \xi_{d(j-1)+k}} \right|_{\xi=\xi_0}= \frac{\left(\sum\limits_{\bar j=1}^n Q_{\bar ji}\xi^0_{d(\bar j-1)+k}\right)Q_{j i}}{\sqrt{\sum\limits_{\bar k=1}^d\left(\sum\limits_{\bar j=1}^n Q_{\bar ji}\xi^0_{d(\bar j-1)+\bar k}\right)^2}} = \mathcal{D}_{ik}Q_{ji},
\label{eq:D_xi_phi}
\end{equation}
where $\mathcal{D}$ is the $m\times d$-matrix with $(i,k)$ entry
\begin{equation}
\mathcal{D}_{ik}=\frac{1}{\varphi_i(\xi_0)}\sum\limits_{\bar j=1}^n Q_{\bar ji}\xi^0_{d(\bar j-1)+k}.
\label{eq:unit_vectors_along_springs}
\end{equation}
Observe that $i$-th row of $\mathcal{D}$ is the {\it unit vector in the direction from the terminus to the origin of spring $i$ in reference configuration $\xi_0$}, i.e. the direction of such unit vector is opposite to the chosen orientation in the geometric directed graph, corresponding to $Q$ with the nodes placement $\xi_0$. 

{\color{black} 
The geometric meaning of equation \eqref{eq:gc} is to guarantee that total elongations $x$ are geometrically possible (up to the linear approximation), so we call it the {\it geometric constraint}. Formula \eqref{eq:D_xi_phi} can be used to  compute matrix $D_{\xi_0} \varphi$ from incidence matrix $Q$ and reference configuration $\xi_0$. \eqref{eq:gc} also appears in the literature \cite[(2.6)]{Lubensky2015} and \cite[(3.17)]{Moreau1973} ({\color{black}formulated for individual springs in the latter}). 
The counterpart of \eqref{eq:gc} in classical continuum mechanics is the displacement-strain relation, see e.g. \cite[(3.7.15)]{MechanicsTextbook}, \cite[(2.55)]{han}.


\subsection{Overview of the rigidity properties which follow {\color{black}from} \eqref{eq:gc}.}
\label{ssect:rigidity_properties1}

Along with the derivation of the equations of the Lattice Spring Model, we would like to provide {\color{black} the} interested reader with some insights and common terminology coming from closely related areas of rigidity theory and structural mechanics. The terminology is taken in large part from the summary \cite[Sect. 2.1, 2.2]{Lubensky2015} and also from \cite[Ch. 8 and 9]{Alfakih2018_rigidity},\cite{RothWhiteley1981},\cite{AsimowRothRigidityI, AsimowRothRigidityII, Roth1981, Whiteley1992}. In the discussion of this section we do not concern ourself with the elasto-plastic properties of springs, and only focus on the kinematics of the lattice, related to the underlying graph structure, i.e. equation \eqref{eq:gc}. In sections to follow, when {\color{black} new} equations will be introduced in the model, we will {\color{black} relate such equations} to the corresponding rigidity terminology and properties.

Matrix $D_{\xi_0} \varphi$ is known in structural mechanics as the {\it compatibility matrix} \cite{Lubensky2015}. Its kernel (nullspace) ${\rm Ker} \, (D_{\xi_0}\varphi)$ is called {\it the set of zero modes} \cite[Sect. 2.1]{Lubensky2015}, and, as it can be observed from \eqref{eq:gc}, the set of zero modes is the space of all infinitesimal displacements of the nodes which do not change any lengths of the springs (up to the linear term). The dimension of the nullspace, i.e. the nullity of the compatibility matrix is called {\it the number of zero modes}.

When the lattice is considered in a $d$-dimensional Euclidean space, there always exists the linear space (denote it $\mathcal{L}$) of {\it infinitesimal rigid motions (isometries) of the lattice} within the set of zero modes:
\begin{equation}
\mathcal{L}\subset {\rm Ker}\, (D_{\xi_0} \varphi) \subset \mathbb{R}^{nd}.
\label{eq:zero_modes_space}
\end{equation}
Each vector $\zeta\in \mathcal{L}$ corresponds to a combination of a parallel translation and an infinitesimal rotation of the entire lattice, so that $\zeta_{d(j-1)+k}$ is the $k$-the component of the velocity vector of node $j$ when the lattice is subject to such combined motion ($j\in \overline{1,n},k\in \overline{1,d}$). Equivalently, one can describe $\mathcal{L}$ as the set of infinitesimal displacements of the nodes preserving, up to the linear term, the distances between {\it all} nodes (not only the adjacent ones).

Except for a special degenerate situation (specifically, when all of the nodes located along a single line, but $d=3$) the dimension of $\mathcal{L}$ ({\it the number of the rigid motions of the lattice}) coincides with the number of rigid motions in space $\mathbb{R}^d$, which is well known to be  (see e.g. \cite[p. 188]{Alfakih2018_rigidity}) $d$ parallel translations plus $d(d-1)/2$ rotations. Therefore
\[
{\rm dim}\, \mathcal{L}=\frac{d(d+1)}{2},
\]
and 
\[
{\rm dim}\, {\rm Ker}\, (D_{\xi_0} \varphi)\geqslant\frac{d(d+1)}{2}.
\]

{\color{black} Zero modes $\zeta \in {\rm Ker}\, (D_{\xi_0} \varphi) \setminus \mathcal{L}$ correspond to the infinitesimal motions of one part of the lattice relative to another. Such modes are called {\it mechanisms} in the engineering literature and {\it floppy modes} in physics literature \cite{Lubensky2015}. The structure of the lattice (specifically, $Q$ and $\xi_0$) dictates whether they are present or not.}

\begin{Definition}{\bf  \cite[Ch. 9]{Alfakih2018_rigidity},\cite{Whiteley1992}.} A lattice with no mechanisms, i.e. such that all its zero modes are rigid motions ($\mathcal{L}= {\rm Ker}\, (D_{\xi_0} \varphi)$) is called {\it infinitesimally rigid}. 
\label{def:inf_rigidity}
\end{Definition}
\noindent Fig. \ref{fig:rigidity} illustrates Definition \ref{def:inf_rigidity} with essential examples for the case $d=2$.

\begin{figure}[H]\center
\includegraphics{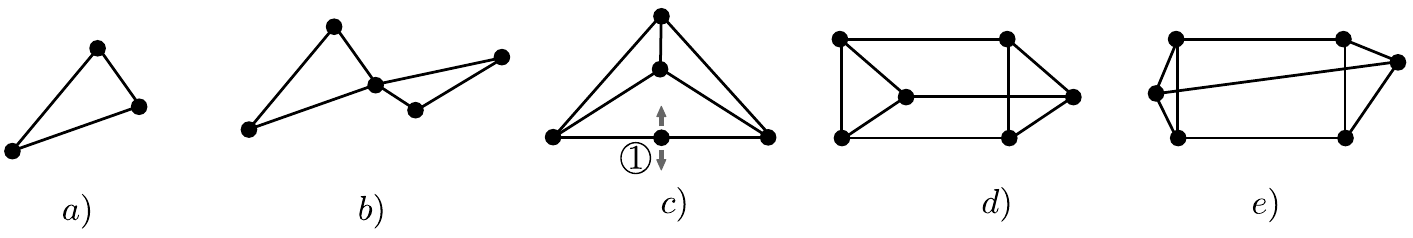}
\caption{
\footnotesize Examples to illustrate Definition \ref{def:inf_rigidity} in case $d=2$. In all of the examples ${\rm dim}\, \mathcal{L} = \frac{d(d+1)}{2}=3$.
a) A triangle makes the simplest infinitesimally rigid lattice,  ${\rm dim \, Ker}\,(D_{\xi_0} \varphi ) = 3$. 
b) Parts of the lattice can freely move relative to each other, which makes the lattice not infinitesimally rigid, ${\rm dim \, Ker}\,(D_{\xi_0} \varphi ) = 4 > 3$. 
c) This lattice is not infinitesimally rigid because a small displacement of node $1$ in the directions orthogonal to all of the adjacent edges is a floppy mode, ${\rm dim \, Ker}\,(D_{\xi_0} \varphi ) = 4 > 3$. 
d) When $d=2$ this lattice is is {\it not infinitesimally rigid} as the identical triangles can synchronously rotate in the plane without changing the lengths of the horizontal edges, ${\rm dim \, Ker}\,(D_{\xi_0} \varphi ) = 4 > 3$. 
e) This lattice has the same graph structure as d, however it is infinitesimally rigid when $d=2$, ${\rm dim \, Ker}\,(D_{\xi_0} \varphi ) = 3$.
} 
\label{fig:rigidity}
\end{figure} 
}

\subsection{Additional constraint and kinematic determinacy}
\label{ssect:additional_constraint}

Along with the geometric constraint we introduce an additional constraint of the form 
\begin{equation}
R(\zeta+\xi_0)+r(t)=0
\tag{LSM2}\label{eq:bc}
\end{equation}
in which $\zeta\in \mathbb{R}^{nd}$ is, again, the displacement vector of the nodes (so that $\zeta+\xi_0$ is a coordinate vector for the nodes),  $R$ is a given $q\times nd$-matrix and $r$ is a given function of time with $q$-vector values for some $q\in \mathbb{N}$. We call equation \eqref{eq:bc}, function $r$ and number $q$, respectively, the {\it external displacement constraint},  the {\it displacement load} and the {\it {\color{black}number} of external displacement constraints}. In turn, {\color{black} we say that $\zeta\in \mathbb{R}^{nd}$ is a {\it feasible displacement} when it satisfies \eqref{eq:bc}}. {\color{black} Naturally, we require the external displacement constraint to be well posed, which means
\begin{Assumption}
Matrix $R$ is of full row rank, i.e.
\begin{equation}
{\rm rank}\,  R = q.
\label{eq:rowsR_indep}
\end{equation}
\label{ass:rowsR_indep}
\end{Assumption}
}
Equation \eqref{eq:bc} corresponds to the displacement boundary condition in classical continuum mechanics, see e.g. \cite[(4.1)]{PlasticityAndGeotechnics2006}.

{\color{black}
In the context of a lattice, defined by both \eqref{eq:gc} and \eqref{eq:bc} the following concept is a counterpart of rigidity. 
\begin{Definition}{\bf \cite[Sect. 2.2]{Lubensky2015}} A lattice endowed with constraint \eqref{eq:bc} is called {\it kinematically determinate} if and only if its feasible displacements of the nodes $\zeta\in \mathbb{R}^{nd}$ correspond one-to-one to elongations vectors $x\in\mathbb{R}^m$, i.e. when map $\zeta \mapsto x$, given by \eqref{eq:gc} and restricted to those $\zeta$ which satisfy \eqref{eq:bc}, {\color{black} can be inverted}.
\end{Definition}
Clearly, a lattice considered without the external displacement constraint would never be kinematically determinate in the sense of this definition, due to $\frac{d(d+1)}{2}$ dimensions of rigid motions we discussed above. However, in our modeling we would like to have a kinematically determinate system, thus we impose \eqref{eq:bc}. 

\begin{Definition} For a lattice endowed with additional constraint \eqref{eq:bc} we define the set of zero modes and the number of zero modes as, respectively, the kernel 
\[
{\rm Ker}\, (D_{\xi_0} \varphi) \cap {\rm Ker}\, R ={\rm Ker} \begin{pmatrix}D_{\xi_0} \varphi \\ R \end{pmatrix}  \subset \mathbb{R}^{nd},
\]
and its dimension, where {\color{black} we call} matrix $\begin{pmatrix}D_{\xi_0} \varphi \\ R \end{pmatrix}$ the {\it enhanced compatibility matrix}.
\label{def:constr_lattice_zero_modes}
\end{Definition}
In terms of the enhanced compatibility matrix we demand
\begin{Assumption} The lattice at the reference configuration is kinematically determinate, i.e.
\begin{equation}
{\rm Ker}\, (D_{\xi_0} \varphi) \cap {\rm Ker}\, R= {\rm Ker} \begin{pmatrix}D_{\xi_0} \varphi \\ R \end{pmatrix} = \{0\}
\label{eq:equivalent_condition}
\end{equation}
or, equivalently, 
\begin{equation}
{\rm rank}\,  \begin{pmatrix} (D_{\xi_0} \varphi)^T & R^T\end{pmatrix} = nd.
\label{eq:trivial_kernel_intersect}
\end{equation}
\label{ass:kinematic_determinacy}
\end{Assumption}

Constraint \eqref{eq:bc} restricts the motions of the nodes and, when \eqref{eq:equivalent_condition} holds, the number of zero modes in the constrained lattice is reduced to zero. Indeed, equivalent condition \eqref{eq:trivial_kernel_intersect} guarantees, that the enhanced compatibility matrix has a left inverse by Proposition \ref{prop:MP3}, therefore the lattice endowed with \eqref{eq:bc} is kinematically determinate (a similar procedure is mentioned at the end of Section 2.2 in \cite{Lubensky2015}). 
While Assumption \ref{ass:kinematic_determinacy} may seem restrictive, it yield {\color{black} several} benefits, both for the simplification of the {\color{black} computations} and {\color{black} for the} applicability of the model. We will summarize {\color{black} these benefits in Section \ref{ssect:msp-rigitity}} after we {\color{black} complete the main analytical derivations of the paper.}
}


\subsection{Additive decomposition and constitutive laws}
Each spring in the model is elasto-plastic, and the standard approach to modeling  of such behavior is the decomposition of total elongation $x_i$ of spring $i$ into elastic component $\varepsilon_i$ and plastic component $p_i$, so that 
\begin{equation}
x_i=\varepsilon_i+p_i,
\label{eq:x_individual}
\end{equation}
{\color{black} or}
\begin{equation}
x=\varepsilon+p
\tag{LSM3}\label{eq:x}
\end{equation}
{\color{black} as an equation in $\mathbb{R}^m$.}

\begin{figure}[h]\center
\includegraphics{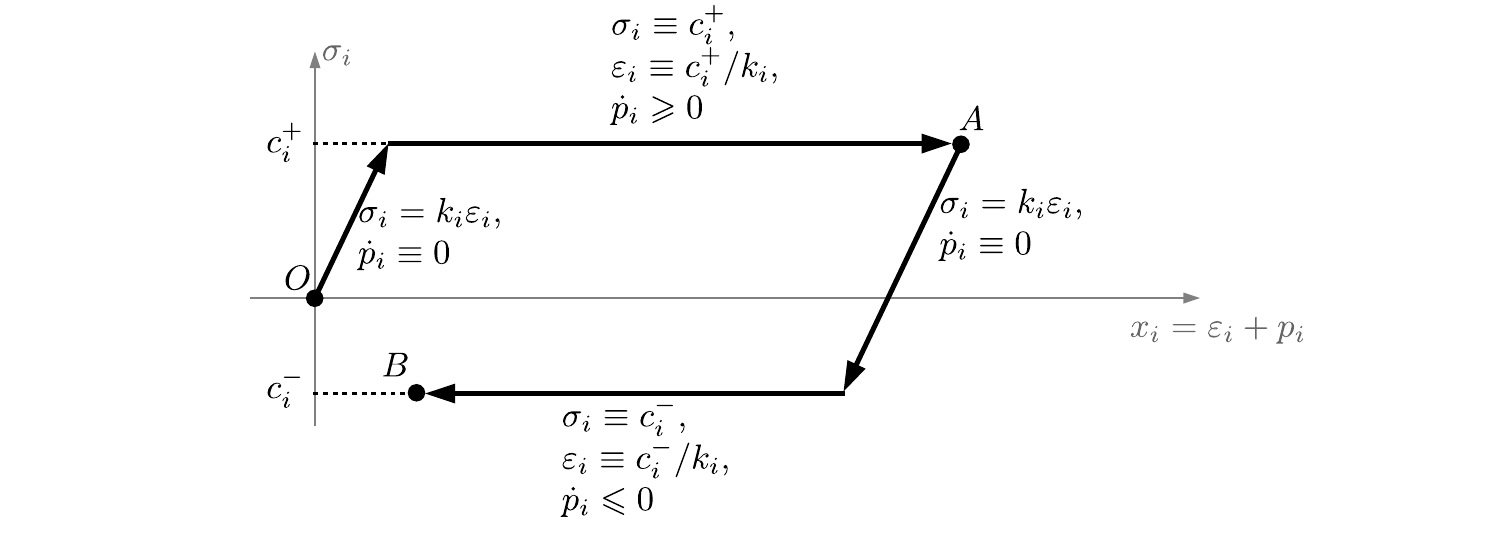}
\caption{\footnotesize Stress-strain behavior of an individual elasto-perfectly plastic spring with parameters $k_i, c_i^+, c_i^-$ under successive stretching ($OA$) and compressing ($AB$).} 
\label{fig:stress-strain}
\end{figure}

Each individual spring is characterized by its stiffness $k_i>0$ and its elasticity interval $(c_i^-, c_i^+)\subset \mathbb{R}$ with stretching and compressing yielding strengths $c_i^+$ and $c_i^-$, respectively (see Fig. \ref{fig:stress-strain}).

 Elastic energy of the individual spring is $\frac{1}{2}k_i \varepsilon_i^2$, and the elastic energy of the whole system is 
\[
E(\varepsilon)=\frac{1}{2}\sum_{i=1}^m k_i \varepsilon_i^2 = \frac{1}{2}\varepsilon^T K \varepsilon,
\]
where $K={\rm\, diag}(k_i),{i\in \overline{1,m}}$ is the $m\times m$ diagonal matrix of the stiffness coefficients. Elastic elongation $\varepsilon_i$ and stress $\sigma_i$ of an individual spring are connected via Hooke's law
\begin{equation}
\sigma_i=k_i\varepsilon_i
\label{eq:individual_hooke}
\end{equation}
the constitutive law of elasticity for the entire system is
\begin{equation}
\sigma = K \varepsilon,
\tag{LSM4}\label{eq:hooke}
\end{equation}
where $\sigma\in \mathbb{R}^m$ is the vector of stress values $\sigma_i$ of the springs $t\in \overline{1,m}$. In classical continuum mechanics, \eqref{eq:hooke} corresponds to the constitutive law of a linearly elastic solid \cite[(5.2.3)]{MechanicsTextbook}, \cite[(2.56)]{han}.

\begin{figure}[h]\center
\includegraphics{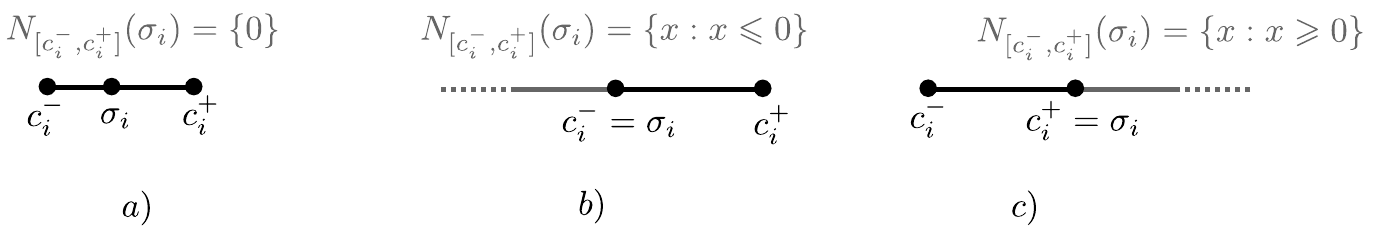}
\caption{\footnotesize The normal cone to the interval $[c_i^-,c_i^+]\subset \mathbb{R}$ at $\sigma_i$ for different locations of $\sigma_i$.} 
\label{fig:1d_normal_cone}
\end{figure}

The constitutive law for the plastic part (also called {\it plastic flow rule} in the literature, {\color{black} see} e.g. \cite{PlasticityAndGeotechnics2006}) of an individual spring is 
\begin{equation}
\begin{array}{rl}
\dot p_i = 0 &\text{ if } \sigma_i \in (c_i^-, c_i^+),\\
\dot p_i \geqslant 0 &\text{ if } \sigma_i = c_i^+,\\
\dot p_i \leqslant 0  &\text{ if } \sigma_i = c_i^-,
\end{array}
\label{eq:individual_plastic1}
\end{equation}
which is a common description of a plastic element, see e.g. \cite[(7)]{Martins2007}. Using the notation of the normal cone \eqref{eq:normal_cone} in $\mathbb{R}^1$ we rewrite \eqref{eq:individual_plastic1} as 
\begin{equation}
\dot p_i \in N_{[c_i^-, c_i^+]}(\sigma_i),
\label{eq:1d_normal_cone}
\end{equation}
Figure \ref{fig:1d_normal_cone} illustrates geometrically how \eqref{eq:individual_plastic1} and \eqref{eq:1d_normal_cone} are equivalent. 

{\color{black}As a side note, an equivalent description of the nonlinear behavior of an individual elasto-plastic spring \eqref{eq:x_individual}, \eqref{eq:individual_hooke}, \eqref{eq:1d_normal_cone} illustrated by Fig. \ref{fig:stress-strain} is given by the Stop operator in the theory of hysteresis \cite{KrasnoselskyPokrovsky, Krejci1996}. We will not use it in the current paper, but some new results from the theory of hysteresis \cite{BrokateKrejci2015} may be useful for optimization of elastoplastic media in future research.}

We can combine such constitutive relations for all $i\in\overline{1,m}$ into a single expression
\begin{equation}
\dot p\in N_C(\sigma),
\tag {LSM5} \label{eq:plasticity}
\end{equation}
where $C$ is the Cartesian product of all the intervals $[c_i^-,c_i^+]$:
\begin{equation}
C= [c_1^-,c_1^+]\times [c_2^-,c_2^+]\times \dots \times [c_m^-,c_m^+] =  \prod_{i=1}^m [c_i^-,c_i^+] \subset \mathbb{R}^m.
\label{eq:admissible stresses}
\end{equation}
The geometric meaning of the normal cone in \eqref{eq:plasticity} can be observed from Figure \ref{fig:normal-cone-def} a-d.
The constitutive law \eqref{eq:plasticity} can also be shown to follow from the {\it principle of maximal plastic work} (see e.g. \cite[p. 57]{han}) similarly to its counterpart in the continuum plasticity theory, {\it the plastic flow rule in the normality form}, \cite[4.35]{han}, \cite[cf5$'$]{DalMaso2006}.
{\color{black}
\subsection{Overview of the static properties in the setting of \eqref{eq:gc}}
\label{ssect:statics_gc_only}
Our model of the lattice is {\it quasi-static},  meaning that the lattice stays at an  equilibrium at all times. Specifically, we refer to the following general definition
\begin{Definition}\cite[p. 17]{Goldstein2002}
\label{def:equilibrium}
A system of $m$ particles is said to be at an {\it equilibrium} when the total force on each particle vanishes. 
\end{Definition}
However, so far the only ``force'' term introduced in the model was stress variable $\sigma\in \mathbb{R}^m$ of \eqref{eq:hooke}, which is related to the springs, not particles (the nodes in our case). In this section we construct the realizations of stresses at the nodes and discuss the main concepts on the {\it statics} of the lattice. The statics is tightly connected to rigidity {\color{black}that} we discussed in Section \ref{ssect:rigidity_properties1}. In a similar manner, we begin with a lattice defined by \eqref{eq:gc} only (without \eqref{eq:bc}) and discuss {\color{black} its} static properties, so in this sense the current section is a counterpart of Section \ref{ssect:rigidity_properties1}. In turn, similarly to Section \ref{ssect:additional_constraint}, we will amend the construction by taking into account the external constraint \eqref{eq:bc} in Section \ref{ssect:static_properties_with_bc} below.

Let us start by establishing a formula connecting scalar variables $\sigma_i$ of stresses in springs $i\in \overline{1,m}$ to the corresponding vector forces at the nodes.
\begin{Proposition}
\label{prop:stresses_to_forces}
For a vector $\sigma\in \mathbb{R}^m$ of stresses produced by \eqref{eq:hooke} the corresponding forces, exerted by the springs at the nodes are given by 
\begin{equation}
-(D_{\xi_0}\varphi)^T\sigma \in \mathbb{R}^{nd}.
\label{eq:stress_realizations}
\end{equation}
in which the $(d(j-1)+k)$-th component is applied to node $j\in \overline{1,n}$ along axis $k\in \overline{1,d}$.
\end{Proposition}
\noindent{\bf Proof.}
Note, that due to Hooke's law in the form \eqref{eq:hooke} with positive diagonal matrix $K$, the situation of $\sigma_i>0$ corresponds to a positive elongation $\varepsilon_i$, i.e. it means a {\it contraction force} in a particular spring $i\in \overline{1,m}$.  Observe that at an individual node $j\in\overline{1,n}$, the stresses of the incident springs add up to vector 
\[
\left(\sum_{i=1}^m -\sigma_i \mathcal{D}_{ik}Q_{ji}\right )_{k\in\overline{1,d}}\in \mathbb{R}^d.
\]
Indeed, if node $j$ is accounted as a terminus of spring $i$, then a contraction force in spring $i$ would act with the magnitude $\sigma_i$ in the direction $(\mathcal{D}_{ik})_{k\in\overline{1,d}}$ (a unit vector, given by \eqref{eq:unit_vectors_along_springs}), but $Q_{ji}=-1$ at a terminus, thus we have an extra minus sign. The similar argument can be done when $Q_{ij}=1$ (node $i$ is an origin) and in case of $\sigma_i<0$ (the stress is an expansion force). Due to \eqref{eq:D_xi_phi}, the stress realizations for all $n$ nodes can be written as \eqref{eq:stress_realizations}. $\blacksquare$

Matrix $(D_{\xi_0} \varphi)^T$ is known in structural mechanics as the {\it equilibrium matrix} \cite{Lubensky2015}. 

Let $f(t)\in \mathbb{R}^{nd}$ be the external force (stress load), in which $f_{d(j-1)+k}$ is the $k$-th component of the force vector, applied to node $j$.  According to Definition \ref{def:equilibrium} and Proposition \ref{prop:stresses_to_forces}, the equation of equilibrium in a lattice defined by \eqref{eq:gc} is
\[
-(D_{\xi_0} \varphi)^T\sigma +f(t)=0,
\]
and we can see from {\color{black} here} that not every stress load $f(t)\in \mathbb{R}^{nd}$ {\color{black} can be balanced by} a corresponding stress vector $\sigma$, but {\color{black} only} those from ${\rm Im}\, (D_{\xi_0} \varphi)^T$. Such stress loads are said to be {\it resolvable} by the stresses in the lattice \cite[p. 10]{Whiteley1992}, \cite[p. 424]{RothWhiteley1981}. An important special class of stress loads, related to rigid motions can be described as the following:
\begin{Definition} {\bf\cite[Sect. 9.3]{Alfakih2018_rigidity}}
Let $\mathcal{L}$ be the linear subspace of all infinitesimal rigid motions of the lattice as in \eqref{eq:zero_modes_space}. Then stress loads 
$
f\in \mathcal{L}^\perp
$
are called {\it equilibrium (stress) loads}.
\label{def:equilibrium_loads}
\end{Definition}

\noindent In particular, when $d=3$, Definition \eqref{def:equilibrium_loads} yields the following criterion:  $f\in \mathbb{R}^{nd}$ is an equilibruim load if and only if
\begin{equation}
\sum_{j=1}^n \left(f_{3(j-1)+k}\right)_{k\in\overline{1,3}} = 0, \qquad \sum_{j=1}^n \left(\xi^0_{3(j-1)+k}\right)_{k\in\overline{1,3}} \times \left(f_{3(j-1)+k}\right)_{k\in\overline{1,3}} = 0.
\label{eq:sum_of_forces_and_torques}
\end{equation}
\noindent Equations \eqref{eq:sum_of_forces_and_torques} have the physical interpretation which agree with natural intuition: applied stress load should have zero net force and zero net torque. For concrete formulas in the case of a general $d$ we refer to \cite[Sect. 9.3]{Alfakih2018_rigidity}.
 

The concepts of resolvable loads and equilibrium loads lead to an equivalent definition of a rigid lattice in terms of stress loads instead of displacements:
\begin{Definition}{\bf \cite[Def. 4.2]{RothWhiteley1981}, \cite[Sect. 9.3]{Alfakih2018_rigidity}} A lattice is called {\it statically rigid} if and only if every equilibrium stress load
is resolvable.
\end{Definition}

\begin{Theorem}{\bf(Whiteley and Roth \cite[Th. 4.3]{RothWhiteley1981})} A lattice is infinitesimally rigid if and only if it is statically rigid.
\label{th:rigidity_equivalence}
\end{Theorem} 
\noindent However, in the situation where the lattice is not infinitesimally rigid, the set of resolvable loads depends on the graph structure of the lattice, specifically on the mechanisms it has.

Another fundamental theorem which connects kinematic and static properties of the lattices is called {\it index theorem}. Instead of resolvable loads (i.e. the image of the equilibrium matrix) it focuses on the kernel of the equilibrium matrix.  The kernel of $(D_{\xi_0}\varphi)^T$ consist of {\it states of self-stresses}, i.e. the stress vectors $\sigma$ which produce zero resultant forces at the nodes, and the dimension of the kernel (the nullity of the equilibrium matrix) is called {\it the number of states of self-stresses}, see Fig. \ref{fig:states-of-self-stresses}.

\begin{figure}[h]\center
\includegraphics{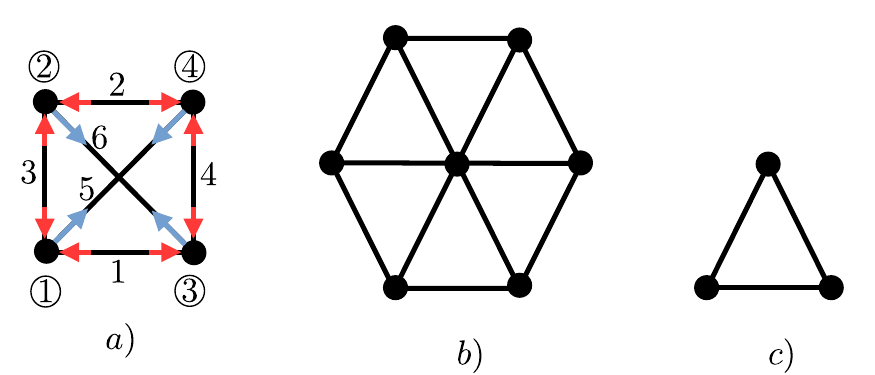}
\caption{
\footnotesize a) An example of a lattice in $\mathbb{R}^2$ which can sustain a single state of self-stress: outer edges 1,2,3 and 4 are pushing the vertices apart, but edges 5 and 6 compensate for that and bind the lattice together at the equilibrium with $\sigma\neq 0$. b) A part of a triangular lattice which in $\mathbb{R}^2$ sustains a single state of self-stress by the same principle as a. c) A single triangle has no states of self-stress (it is statically determinate).} 
\label{fig:states-of-self-stresses}
\end{figure} 

Due to the rank-nullity theorem, the numbers of zero modes and states of self-stresses are tightly connected, which fact is known as the follows:
\begin{Theorem}{\bf(Index theorem, \cite{Lubensky2015})} Let a lattice be defined in the sense of relation \eqref{eq:gc}. Then
\begin{equation*}
(\text{\# zero modes})-(\text{\# of states of self-stresses})= nd-m,
\label{eq:index_th}
\end{equation*}
i.e.
\begin{equation*}
{\rm dim\, Ker}\, (D_{\xi_0}\varphi) - {\rm dim\, Ker}\, (D_{\xi_0}\varphi)^T=nd-m.
\end{equation*}
\label{th:index_th}
\end{Theorem}

\subsection{Derivation of the equation of equilibrium in the final form}
\label{ssect:equation_of_equilibrium}

\noindent In the current section we derive the equation of equilibrium of the lattice with external displacement constraint \eqref{eq:bc} imposed along with \eqref{eq:gc}. In the context of kinematical constraints such as \eqref{eq:bc} the equation of equilibrium can be handled via
\begin{Proposition}{\bf  The Principle of Virtual Work \cite[Ch. III.1]{LanczosMechanics}.}
We assume that the given external forces $F_1, F_2,\dots, F_n$ act at the points $P_1, P_2,\dots , P_n$ of
the system. The virtual displacements of these points will be denoted by $\delta R_1, \delta R_2,\dots , \delta R_n$. These virtual displacements must be in harmony with the given kinematical constraints, and we shall assume that they are reversible, i.e. the given constraints do not prevent us from changing an arbitrary $\delta R_i$ into $-\delta R_i$. Now the principle of virtual work asserts that the given mechanical
system will be in equilibrium if, and only if, the total virtual work of all the impressed forces vanishes:
\[
F_1\cdot \delta R_1 + F_2 \cdot \delta R_2 + \dots + F_n \cdot \delta R_n = 0.
\]
\label{prop:principle_of_virtual work}
\end{Proposition}
\noindent We also refer to  \cite[\S1.4]{Goldstein2002} and \cite[\S 2.2]{Haslach2011} regarding the principle of virtual work.

In addition, we will require the following  technical construction, which will allow to write the equation of equilibrium in a desired form. Consider the $(nd)\times(m+q)$ matrix $\begin{pmatrix} (D_{\xi_0} \varphi)^T & R^T\end{pmatrix}$ from condition \eqref{eq:trivial_kernel_intersect}.
Its Moore-Penrose pseudoinverse (see Proposition \ref{prop:MP1}) is an $(m+q)\times(nd)$-matrix, so we can write the pseudoinverse as $\begin{pmatrix} H\\ H'\end{pmatrix}$, where $H$ and $H'$ are $m\times(nd)$ and $q\times(nd)$ matrices, respectively:
\begin{equation}
\begin{pmatrix} (D_{\xi_0} \varphi)^T & R^T\end{pmatrix}^+=\begin{pmatrix} H\\ H'\end{pmatrix}=\begin{pmatrix} D_{\xi_0} \varphi \\ R\end{pmatrix}\left(\begin{pmatrix} (D_{\xi_0} \varphi)^T & R^T\end{pmatrix}\begin{pmatrix} D_{\xi_0} \varphi \\ R\end{pmatrix}\right)^{-1},
\label{eq:new_equilibrium_matrix_pseudoinverse}
\end{equation}
where the last equality holds true due to condition \eqref{eq:trivial_kernel_intersect} and Proposition \eqref{prop:MP3}.

\begin{Proposition}
Consider the lattice, described by an incidence matrix $Q$ at a reference configuration $\xi_0$ with external displacement constraint \eqref{eq:bc}. Assume that the stresses are  connected to elastic elongations by Hooke's law \eqref{eq:hooke} and the stress load is given as $f(t)\in \mathbb{R}^{nd}$. If the lattice is at the equilibrium, then 
\begin{equation}
(D_{\xi_0} \varphi)^T (-\sigma+H f(t))\in {\rm Im}\,R^T.
\tag{LSM6}
\label{eq:balanceLaw} 
\end{equation}
Provided with the condition \eqref{eq:trivial_kernel_intersect}, the converse is also true for any $f(t)$: if \eqref{eq:balanceLaw} holds, then the lattice is at the equilibrium.  
\label{prop:balance_law_derivation}
\end{Proposition}
\noindent{\bf Proof.}
First we show that equilibrium state of the lattice implies \eqref{eq:balanceLaw}. 
From the above considerations in Section \ref{ssect:statics_gc_only} and Proposition \ref{prop:stresses_to_forces} in particular we note that forces $F_1, \dots F_n$ in the formulation of Proposition \ref{prop:principle_of_virtual work} in our case are the groups of components of size $d$ of 
\[
-(D_{\xi_0} \varphi)^T\sigma +f(t) \in \mathbb{R}^{nd}
\]
(the term ``external forces'' in the formulation of Proposition \ref{prop:balance_law_derivation} has the meaning ``external to the point $P_i$ and the kinematical constraint'' in the abstract setting of \cite[Ch. III.1]{LanczosMechanics}, and it should not be confused with $f(t)$, the external forces to the lattice in this paper). 
The kinematical constraint in our case is \eqref{eq:bc}, which is {\it affine} (hence, it is {\it reversible}), and for our system the principle of virtual work yields the following {\it criterion of equilibrium}:
\[
\left (-(D_{\xi_0} \varphi)^T\sigma +f(t)\right)^T \delta\zeta = 0 \qquad \text{for any } \delta\zeta\in {\rm Ker}\, R.
\]
Notice, that time-dependence of \eqref{eq:bc} does not play a role here, since displacement $\delta\zeta$ is virtual, see \cite[p. 17]{Goldstein2002}.
Equivalently written, the equilibrium condition is
\[
-(D_{\xi_0} \varphi)^T\sigma +f(t)\in ({\rm Ker}\, R)^\perp,
\]
which is, in turn, equivalent to (see e.g. \cite[Th. 4.45]{OlverLinearAlgebra2018})
\begin{equation}
-(D_{\xi_0} \varphi)^T\sigma +f(t)\in {\rm Im}\, R^T.
\label{eq:balance-law-derivation1}
\end{equation}
Finally, we want to express $f(t)$ so that we can factor out $(D_{\xi_0} \varphi)^T$. For the $(nd)\times(m+q)$ matrix $
\begin{pmatrix} (D_{\xi_0} \varphi)^T & R^T\end{pmatrix}$
we compute its Moore-Penrose pseudoinverse matrix of dimensions  $(m+q)\times(nd)$ in the form 
\[
\begin{pmatrix} H\\ H'\end{pmatrix}=\begin{pmatrix} (D_{\xi_0} \varphi)^T & R^T\end{pmatrix}^+,
\] where $H$ and $H'$ are defined as, respectively, $m\times(nd)$ and $q\times(nd)$ parts of the pseudoinverse matrix. It follows from \eqref{eq:balance-law-derivation1} and \eqref{eq:MP1} that there is $\rho\in \mathbb{R}^q$, such that
\begin{multline*}
f(t)=\begin{pmatrix} (D_{\xi_0} \varphi)^T & R^T\end{pmatrix}\begin{pmatrix}\sigma\\ \rho\end{pmatrix}=\begin{pmatrix} (D_{\xi_0} \varphi)^T & R^T\end{pmatrix}\begin{pmatrix} H\\ H'\end{pmatrix}\begin{pmatrix} (D_{\xi_0} \varphi)^T & R^T\end{pmatrix}\begin{pmatrix}\sigma\\ \rho\end{pmatrix}=\\=\begin{pmatrix} (D_{\xi_0} \varphi)^T & R^T\end{pmatrix}\begin{pmatrix} H\\ H'\end{pmatrix} f(t)= (D_{\xi_0} \varphi)^T H f(t) + R^T H' f(t).
\end{multline*}
Plug this expression for $f(t)$ back to \eqref{eq:balance-law-derivation1} to obtain
\[
-(D_{\xi_0} \varphi)^T\sigma + (D_{\xi_0} \varphi)^T H f(t) +  R^T H' f(t) \in {\rm Im}\, R^T.
\] 
Because $R^T H' f(t) \in {\rm Im}\, R^T$ we obtain the equation of equilibrium in the final form \eqref{eq:balanceLaw}.

Conversely, let \eqref{eq:trivial_kernel_intersect} and \eqref{eq:balanceLaw} hold. Recall that $ \begin{pmatrix} (D_{\xi_0} \varphi)^T & R^T\end{pmatrix}\begin{pmatrix} H\\ H'\end{pmatrix}$ is a projection matrix on ${\rm Im} \begin{pmatrix} (D_{\xi_0} \varphi)^T & R^T\end{pmatrix}$  (see Proposition \ref{prop:MP2}), but  \eqref{eq:trivial_kernel_intersect} means that ${\rm Im} \begin{pmatrix} (D_{\xi_0} \varphi)^T & R^T\end{pmatrix}$ spans the entire $\mathbb{R}^{nd}$, therefore the projection matrix is the identity matrix and we have
\[
f(t)=\begin{pmatrix} (D_{\xi_0} \varphi)^T & R^T\end{pmatrix}\begin{pmatrix} H\\ H'\end{pmatrix} f(t)= (D_{\xi_0} \varphi)^T H f(t) + R^T H' f(t).
\]
Plug this in \eqref{eq:balanceLaw} to obtain
\[
-(D_{\xi_0} \varphi)^T\sigma + (D_{\xi_0} \varphi)^T H f(t) +  f(t)\in {\rm Im}\,R^T+  (D_{\xi_0} \varphi)^T H f(t) + R^T H' f(t),
\]
equivalent to \eqref{eq:balance-law-derivation1}, which is, in turn, equivalent to the equation of equilibrium.  $\blacksquare$

Equation \eqref{eq:balanceLaw}  is a slightly modified version of equations from the literature \cite[(2.4)]{Lubensky2015}, \cite[(3.23)]{Moreau1973}. The counterpart of \eqref{eq:balanceLaw} in classical continuum mechanics is the equation of equilibrium \cite[2.61]{han}, \cite[cf3]{DalMaso2006}. {\color{black}Additionally, in the particular case of absent plastic deformation ($x = \varepsilon$) formulation of elasticity \eqref{eq:hooke} (with linearization \eqref{eq:gc}, \eqref{eq:D_xi_phi} and realizations of stresses \eqref{eq:stress_realizations}) coincides with the law of pairwise interactions between particles in microscopic elasticity theory \cite[(15)-(16)]{Goldhirsch2002}.}

\subsection{Static properties of the full model}
\label{ssect:static_properties_with_bc}
At last we discuss the static properties of a lattice defined by \eqref{eq:gc},\eqref{eq:bc} and \eqref{eq:balanceLaw}.
\begin{Definition}
For a lattice endowed with external constraint \eqref{eq:bc}, we call {\it resolvable (stress) loads}  the vectors from  
\[
{\rm Im}\, \begin{pmatrix} (D_{\xi_0} \varphi)^T & R^T\end{pmatrix}\subset \mathbb{R}^{nd},
\]
and we call matrix $\begin{pmatrix} (D_{\xi_0} \varphi)^T & R^T\end{pmatrix}$ the {\it enhanced equilibrium matrix}.
\end{Definition}

\noindent Observe from \eqref{eq:balanceLaw} that, while removing zero modes, the additional constraint  leads to a larger set of resolvable stress loads compared to the resolvable stress loads in the sense of Section \ref{ssect:statics_gc_only}. Moreover, along with kinematic determinacy, condition \eqref{eq:trivial_kernel_intersect} also means that \eqref{eq:bc} is tight enough so that {\it all} stress loads from $\mathbb{R}^{nd}$ are resolvable (and we, in fact, used this in the second part of the proof of Proposition \ref{prop:balance_law_derivation}):
\begin{Theorem}
A lattice endowed with external constraint \eqref{eq:bc} is kinematically determinate if and only if any stress load from $\mathbb{R}^{nd}$ is resolvable.
\label{th:all_loads_are_resolvable}
\end{Theorem}
\noindent {\bf Proof.} Indeed, all loads are resolvable if and only if
\begin{equation}
{\rm Im}\, \begin{pmatrix} (D_{\xi_0} \varphi)^T & R^T\end{pmatrix}= \mathbb{R}^{nd}
\label{eq:all_loads_are_resolvable}
\end{equation}
which is equivalent to \eqref{eq:equivalent_condition}, which, in turn, means kinematic determinacy. $\blacksquare$

Our goal here is not only to give the physical meaning to the terms of \eqref{eq:balanceLaw} and \eqref{eq:trivial_kernel_intersect}, but also to stress  that Theorem \ref{th:all_loads_are_resolvable} is   a counterpart of Theorem \ref{th:rigidity_equivalence}. Indeed, Theorem \ref{th:rigidity_equivalence} claims that {\color{black} the elongations-invariant motions in a lattice are limited to the set $\mathcal{L}$ of rigid motions if and only if the lattice can balance any stress load with zero} $\mathcal{L}$-component; loosely speaking, Theorem \ref{th:all_loads_are_resolvable} is the similar claim about set $\{0\}$ instead of $\mathcal{L}$. {\color{black} This shows that kinematic determinacy, which we require in Assumption \ref{ass:kinematic_determinacy}, is the appropriate analogue of the concept of infinitesimal rigidity for lattices with external constraints, as both concept link kinematics and statics in the similar way. In turn, this universal link gives a proper explanation behind the useful fact that it is enough to know that a) the lattice is infinitesimally rigid b) the external constraint \eqref{eq:bc} determines $\frac{d(d+1)}{2}$ components of $\zeta$ to guarantee the validity of   \eqref{eq:balanceLaw} as an equilibrium equation for {\it any} external force $f(t)\in \mathbb{R}^{nd}$.}

The following concept of self-stresses is needed for physical interpretation of the analytic constructions of Sections \ref{ssect:fundamental_spaces} and \ref{sect:sweeping_solution}.
\begin{Definition}
For a lattice endowed with an external constraint \eqref{eq:bc} we define the {\it states of self-stresses} as vectors $\sigma\in \mathbb{R}^m$ such that for some vector $\rho\in \mathbb{R}^q$ (called a {\it reaction} of \eqref{eq:bc}) 
\begin{equation}
\begin{pmatrix} (D_{\xi_0} \varphi)^T & R^T\end{pmatrix} \begin{pmatrix}\sigma \\ \rho\end{pmatrix}=0.
\label{eq:self-stresses_def_const_lattice}
\end{equation}
We call the dimension of the set of such $\sigma$'s the {\it number of states of self-stresses}.
\label{def:self-stresses_def_const_lattice}
\end{Definition}
Observe from the definition that the states of self-stresses constitute a linear hyperplane of $\sigma$'s satisfying \eqref{eq:balanceLaw} when $f(t)=0$. 

\begin{Remark}
We should clarify, that it would be in better agreement with the general approach to call states of self-stresses the vectors from the kernel of the enhanced equilibrium matrix $\begin{pmatrix} (D_{\xi_0} \varphi)^T & R^T\end{pmatrix}$, i.e. the whole vectors $\begin{pmatrix}\sigma \\ \rho\end{pmatrix}\in \mathbb{R}^{m+q}$ satisfying \eqref{eq:self-stresses_def_const_lattice}, but we are not interested in the component $\rho$. If needed, $\rho$ can be easily computed from $\sigma$ because there is a one-to-one correspondence between the states of self-stresses as in Definition \ref{def:self-stresses_def_const_lattice} and ${\rm Ker}\,\begin{pmatrix} (D_{\xi_0} \varphi)^T & R^T\end{pmatrix}$ due to condition \eqref{eq:rowsR_indep}.
\end{Remark}
}

\subsection{Combined equations of quasi-static evolution of an elastic - perfectly plastic Lattice Spring Model}
To summarize, the governing equations of the Lattice Spring Model are
\begin{align}
\text{Geometric constraint:}&&\left(D_{\xi_0} \varphi\right) \zeta&=x& \tag{\ref{eq:gc}}\\
\text{External displacement constraint:}&&R(\zeta+\xi_0) + r(t)&=0& \tag{\ref{eq:bc}}\\
\text{Additive decomposition:}&&x&=\varepsilon+p&\tag{\ref{eq:x}}\\
\text{Hooke's law:}&&\sigma&=K \varepsilon& \tag{\ref{eq:hooke}}\\
\text{Flow rule of perfect plasticity:}&&\dot p &\in N_{C}(\sigma)& \tag{\ref{eq:plasticity}}\\
\text{Equation of equilibrium:}&&(D_{\xi_0} \varphi)^T (-\sigma+H f(t))&\in {\rm Im}\,R^T &\tag{\ref{eq:balanceLaw}}
\end{align}
These equations can be viewed as discrete networks analogues of the corresponding equations for an elasto-plastic continuous medium, see e.g. \cite{han},\cite{Meyer} or \cite{DalMaso2006} and a particular case of an abstract problem, described in \cite[6a]{Moreau1973}. {\color{black} Details on conversion to the abstract problem can be found in \cite[Appendix A]{GM-1}, where we analyzed a system similar to \eqref{eq:gc}-\eqref{eq:balanceLaw} with one spatial dimension.}
{\color{black}
\subsection{The fundamental spaces of the lattice}
\label{ssect:fundamental_spaces}
To lay the groundwork for solving the evolution problem \eqref{eq:gc}-\eqref{eq:balanceLaw} we {\color{black} will consider} linear spaces {\color{black} that are} fundamentally connected with equations \eqref{eq:gc},\eqref{eq:bc}, \eqref{eq:hooke} and \eqref{eq:balanceLaw}, and, therefore, with the structure of the lattice. 

{\color{black} Recall that $K$ is a diagonal matrix of stiffness values $k_i$, hence it is symmetric, positive definite and invertible.} Define the following subspaces of $\mathbb{R}^m$.
\begin{gather}
\mathcal{U}=(D_{\xi_0} \varphi) {\rm Ker}\, R=\{(D_{\xi_0} \varphi) \zeta: R\zeta=0\} 
\label{eq:space_u_def} \\
\mathcal{V}= \{K^{-1}\sigma: (D_{\xi_0}\varphi)^T\sigma\in {\rm Im}\, R^T\}
\label{eq:space_v_def}
\end{gather}
Interpreted mechanically,  {\it $\mathcal{U}$ consists of vectors of total elongations which correspond to feasible displacements with $r(t)=0$ in \eqref{eq:bc}.} Mechanical interpretations of $\mathcal{U}^\perp$ and $\mathcal{V}$ are given by the following proposition. 
\begin{Proposition} Let Assumptions \ref{ass:rowsR_indep}, \ref{ass:kinematic_determinacy} hold true. Then
\begin{enumerate}[\it i)]
\item 
The orthogonal complement $\mathcal{U}^\perp$ consists of states of self-stresses (Definition \ref{def:self-stresses_def_const_lattice}):
\begin{equation}
\mathcal{U}^\perp  = \{\sigma\in \mathbb{R}^m:(D_{\xi_0}\varphi)^T \sigma\in {\rm Im}\, R^T\}.
\label{eq:U_perp}
\end{equation}
Furthermore,
\begin{equation}
\mathcal{V} = K^{-1}\mathcal{U}^{\perp},
\label{eq:V_perp_standard}
\end{equation}
i.e. members of $\mathcal{V}$ are vectors of elastic elongations, corresponding to the states of self-stresses by the Hooke's law. Thus ${\rm dim} \, \mathcal{U}^\perp = {\rm dim} \, \mathcal{V}$ is the number of states of self-stresses.
Moreover, $\mathcal{U}$ and $\mathcal{V}$ are orthogonal complements in sense of weighted inner product \eqref{eq:S-inner-prod} with weights $K$:
\begin{equation}
\mathcal{V} = \mathcal{U}^{\perp_{K}}.
\label{eq:orthogonality_K}
\end{equation}
\item The set of $\sigma$'s satisfying equilibrium equation \eqref{eq:balanceLaw} is an affine translation of $\mathcal{U}^\perp$. Specifically, 
for any $\sigma\in \mathbb{R}^m$ and any $t$
\begin{equation}
\sigma \text{ satisfies } \eqref{eq:balanceLaw} \quad \Longleftrightarrow \quad -\sigma+Hf(t)\in \mathcal{U}^\perp.
\label{eq:equiv_balance_law}
\end{equation}

\item The dimensions of $\mathcal{U}$ and $\mathcal{V}$ are
\begin{equation}
{\rm dim}\,\mathcal{U} = {\rm dim \, Ker} \, R = nd-q,\qquad {\rm dim}\,\mathcal{V}=m-nd+q.
\label{eq:dimensions_U_V}
\end{equation}
\end{enumerate}
\label{prop:U_and_V_orth}
\end{Proposition}
\noindent{\bf Proof.}
\begin{enumerate}[\it i)]
\item We show \eqref{eq:U_perp} directly
\begin{multline}
U^\perp=((D_{\xi_0} \varphi)\,{\rm Ker}\, R)^\perp=\{x\in \mathbb{R}^m: x^T(D_{\xi_0} \varphi)y=0 \text{ for all } y\in {\rm Ker}\, R\}=\\
=\{x\in \mathbb{R}^m: \left((D_{\xi_0} \varphi)^T x\right)^Ty=0 \text{ for all } y\in {\rm Ker}\,R\}=\{x\in \mathbb{R}^m: (D_{\xi_0} \varphi)^T x\in ({\rm Ker}\, R)^\perp\}=\\
=\{x\in \mathbb{R}^m:  (D_{\xi_0} \varphi)^T x\in {\rm Im}\, R^T\},
\end{multline}
where the second equality follows from \eqref{eq:space_u_def} and \cite[Th. 4.45]{OlverLinearAlgebra2018}.
 
\item This can be verified by substituting $-\sigma + H(t)$ as $\sigma$ into \eqref{eq:U_perp}.
\item 
By rank-nullity theorem (see e.g. \cite[Th. 2.49]{OlverLinearAlgebra2018}), Assumption \ref{ass:rowsR_indep} implies that ${\rm dim\, Ker} \,R = nd - q$. From \eqref{eq:space_u_def} we immediately see that ${\rm dim}\, \mathcal{U}\leqslant {\rm dim\, Ker}\, R$. In turn, with Assumption~\ref{ass:kinematic_determinacy} we can guarantee that ${\rm dim}\, \mathcal{U} = {\rm dim\, Ker}\, R$. 

Indeed, consider a basis of ${\rm Ker} \,R$ and arrange its vectors as columns in a matrix $R_0$.  Assume that ${\rm dim}\,\mathcal{U} < {\rm dim\, Ker}\, R$, which means that there exists a nontrivial linear combination of basis vectors of ${\rm Ker}\, R$ (i.e. $R_0 z$ with some nonzero vector of coefficients $z\in \mathbb{R}^{nd-q}$) such that  $\left((D_{\xi_0} \varphi) R_0\right) z = 0$. Hence $R_0 z$ is from ${\rm Ker}\, (D_{\xi_0} \varphi)$. Since $R_0 z$ is a nonzero vector from ${\rm Ker}\, R$ due to the choice of $R_0$ and $z$, we have a contradiction with \eqref{eq:equivalent_condition}.
\end{enumerate}
The proof of the proposition is complete. $\blacksquare$

Observe from \eqref{eq:U_perp}, that the underlying graph structure of the lattice and the external displacement constraint are the two factors which influence the number of states of self-stresses ${\rm dim}\, \mathcal{V}$. In particular, if \eqref{eq:bc} is chosen to be tighter, (i.e. with larger $q$), then the set ${\rm Im}\, R^T$ is also larger. This means that with tighter external displacement constraint more combinations of internal stresses can be balanced out by reactions of the constraint and under the condition of kinematic determinacy \eqref{eq:trivial_kernel_intersect} we have a precise formula for ${\rm dim }\, \mathcal{V}$ in Proposition \ref{prop:U_and_V_orth} (compare with the general Theorem \ref{th:index_th}). 

If $\mathcal{V}=\{0\}$ then the evolution of stresses is uniquely determined by \eqref{eq:balanceLaw} alone, the situation known as {\it static determinacy} \cite[Sect. 2.2]{Lubensky2015}. We consider this case degenerate, because the evolution in plastic regime (i. e. yielding) becomes impossible in such a lattice. {\color{black} In the following assumption we use formula \eqref{eq:dimensions_U_V} to exclude this degenerate situation.}

\begin{Assumption} There are non-trivial states of self-stresses, i.e
\begin{equation}
{\rm dim}\, \mathcal{V} = m-nd+q >0.
\label{eq:at_least_one_self-stress}
\end{equation}
\label{ass:at_least_one_self-stress}
\end{Assumption}
}

\subsection{Example 1: a simple toy network}
\label{ssect:ex1}
\begin{figure}[h]\center
\includegraphics{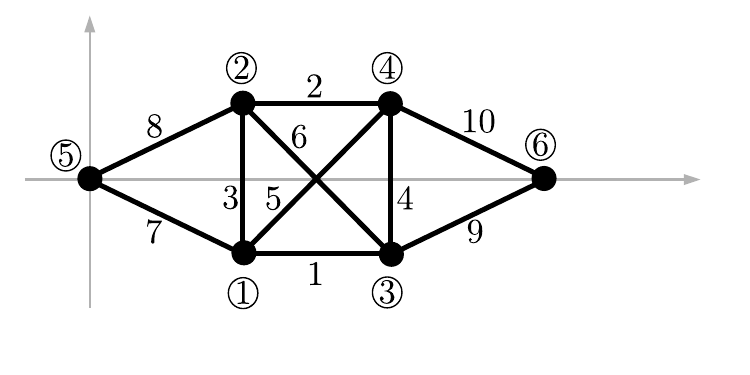}
\caption{
\footnotesize Toy example with 6 nodes (circled numbers) and 10 springs. Displacement loading is applied to nodes 5 and 6.
} \label{fig:ex1}
\end{figure} 
We illustrate the construction with the toy example of the truss shown in Fig. \ref{fig:ex1}. We have $m=10, n=6$ and $d=2$, the incidence matrix and the reference configuration are
\begin{align*}
Q=\begin{pmatrix}
1 & 0 & 1 & 0 & 1 & 0 & 1 & 0 & 0 & 0\\
0 & 1 &-1 & 0 & 0 & 1 & 0 & 1 & 0 & 0\\
-1& 0 & 0 & 1 & 0 &-1 & 0 & 0 & 1 & 0\\
0 &-1 & 0 &-1 &-1 & 0 & 0 & 0 & 0 & 1\\
0 & 0 & 0 & 0 & 0 & 0 & -1&-1 & 0 & 0\\
0 & 0 & 0 & 0 & 0 & 0 & 0 & 0 &-1 &-1      
\end{pmatrix},\\
\xi_0=\left(\begin{array}{cccccccccccc}2 & -1 & 2 & 1 & 4 & -1 & 4 & 1 & 0 & 0 & 6 & 0\end{array}\right)^T.
\end{align*}
{\color{black} As elasto-plastic parameters of springs we take
\begin{gather*}
K={\rm diag}(1,1,1,1,1,1,1,1,1,1)=I_{10}, \qquad C=\{c\in \mathbb{R}^{10}: c_i^-\leqslant c_i\leqslant c_i^+\}, \text{ where }\\
c^-=-c_0\, \widetilde{c},\qquad c^+= c_0\, \widetilde{c},\\
c_0=0.001, \qquad \widetilde{c}=\begin{pmatrix}1& 1& 1& 1& 1/\sqrt{2}& 1/\sqrt{2}& 10& 10& 10& 10\end{pmatrix}^T,
\end{gather*}
i.e. all springs have stiffness $1$ and, up to the scale factor $c_0$, the axis-aligned, diagonal and side springs have stress limits $1$,  $1/\sqrt{2}$ and $10$ respectively.}
In this example we fix the position of node 5 and move node 6 along the x-axis, which can be written as 
\[
\begin{cases}
\zeta_9 +\xi_9^0   =0,\\
\zeta_{10}+\xi_{10}^0 =0,\\
\zeta_{11}+\xi_{11}^0 - (r_0 c_0) t =0,\\
\zeta_{12}+\xi_{12}^0 =0,
\end{cases} \qquad r_0=10.
\]
i.e. external displacement constraint \eqref{eq:bc} contains $q=4$ equations with
\begin{equation}
R=\left(\begin{array}{cccccccccccc}0&0&0&0&0&0&0&0& 1&0&0&0\\
0&0&0&0&0&0&0&0& 0&1&0&0\\
0&0&0&0&0&0&0&0& 0&0&1&0\\
0&0&0&0&0&0&0&0& 0&0&0&1
 \end{array}\right), \quad r(t)=\begin{pmatrix}0\\0\\-(r_0 c_0) t \\ 0\end{pmatrix}.
 \label{eq:example1_constr}
\end{equation}

Now we will verify conditions \eqref{eq:rowsR_indep}, \eqref{eq:trivial_kernel_intersect},\eqref{eq:at_least_one_self-stress} for the example constructed.
Condition \eqref{eq:rowsR_indep} is satisfied trivially.

Meanwhile, condition \eqref{eq:trivial_kernel_intersect} can be verified numerically after one calculates $D_{\xi_0} \varphi$ by formula \eqref{eq:D_xi_phi}. Alternatively, one can use the insight into the rigidity theory and structural properties of the lattice to observe the following. The lattice in Fig. \ref{fig:ex1} is similar to Fig. \ref{fig:states-of-self-stresses}a and, if considered without the external constraint, it has only one state of self-stress (the triangles on the sides do not increase the number  of states of self-stress as long as nodes 5 and 6 are unconstrained). The external constraint \eqref{eq:example1_constr} with $q=4$ makes the lattice kinematically determinate and adds 1 more state of self-stresses.

Finally, condition \eqref{eq:at_least_one_self-stress} holds as
\[
{\rm dim}\, \mathcal{V} = m-nd+q=10-6\cdot 2+4=2>0.
\]

In the next Section \ref{sect:sweeping_intro} we introduce the sweeping process and then in Sections \ref{sect:sweeping_solution} and \ref{sect:leapfrog} we will show how it helps to solve the problem \eqref{eq:gc}-\eqref{eq:balanceLaw}

\section{A Brief Introduction to Sweeping Process Theory}
\label{sect:sweeping_intro}
In this section we explain the concept of the sweeping process as an abstract problem in $\mathbb{R}^n$.  An interested reader is referred to \cite{Kunze2000} for a more comprehensive introduction to sweeping process with detailed theorems and proofs. 

\begin{figure}[h]\center
\includegraphics[scale=0.5]{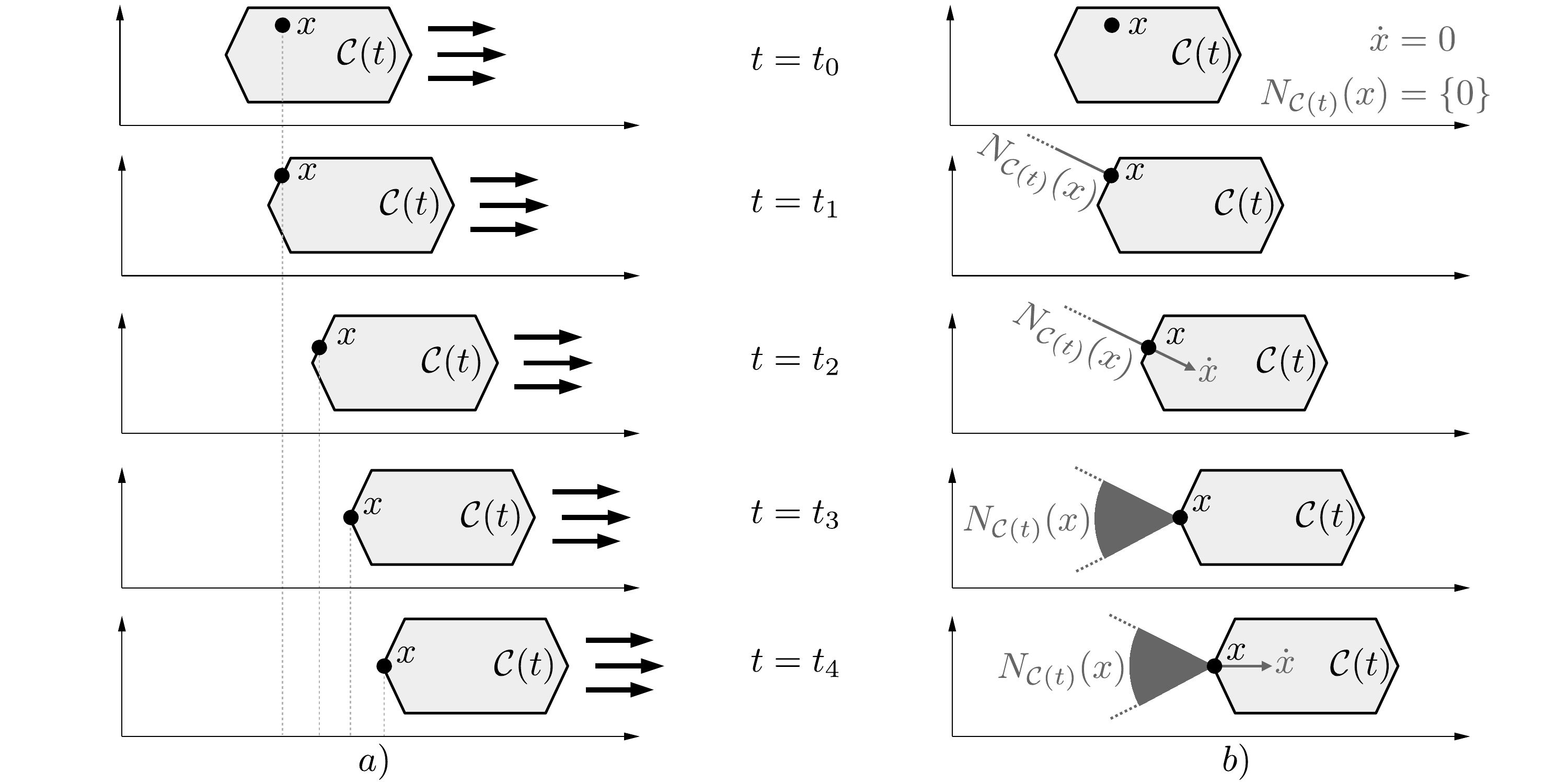}
\caption{
\footnotesize An example of a sweeping process, in which $\mathcal{C}(t)$ is the hexagonal set moving to the right. a) different stages of the solution: when $t_0\leqslant t < t_1$, point $x$ stays immovable as it belongs to the interior of $\mathcal{C}(t)$. At $t=t_1$, point $x$ is hit by the edge of $\mathcal{C}(t)$. When $t_1< t< t_3$ point $x$ is swept by the edge and it also slides towards the left vertex along the edge. $x$ meets the vertex at $t=t_3$ and moves together with it for $t>t_3$. b) We illustrate how \eqref{eq:sp} holds for all $t$ except at isolated moments $t_1$ and $t_3$. When $x$ is in the interior of $\mathcal{C}(t)$, we have $\dot x =0$ and $N_{\mathcal{C}(t)}(x)$ being a singleton set of zero vector. When $t_1< t< t_3$ the normal cone $N_{\mathcal{C}(t)}(x)$ is a one-dimensional ray, orthogonal to the edge of $\mathcal{C}(t)$, and $\dot x$ is directed opposite to that ray. When $t>t_3$ the normal cone $N_{\mathcal{C}(t)}(x)$ is a 2-dimensional cone and it contains the opposite direction to $\dot x$. 
} \label{fig:sweeping}
\end{figure}

 Let us be given with a set-valued function of time:
\begin{equation}
\mathcal{C}: [0,T] \to \mathcal{P}(\mathbb{R}^n)
\label{eq:moving_set_sweeping}
\end{equation}
where $\mathcal{P}(\mathbb{R}^n)$ is the power set (the collection of all subsets) of $\mathbb{R}^n$. In addition, the value $\mathcal{C}(t)$ is assumed to be a nonempty closed convex set for every $t\in [0,T]$. We call $\mathcal{C}(t)$ the {\it moving set}. For an initial condition $x_0\in \mathcal{C}(0)$ the {\it sweeping process} describes the trajectory of a point $x$, originally placed at $x_0$ for $t=0$ and constrained within the moving set $\mathcal{C}(t)$ for all $t\in [0, T]$. The motion of the point $x$ can be characterized as follows: $x$ remains at rest, unless it is ``swept'' by the boundary of $\mathcal{C}(t)$ in order to remain within the moving set (see Fig. \ref{fig:sweeping} a).

The equation of the sweeping process is 
\begin{equation}
-\dot x \in N_{\mathcal{C}(t)}(x)
\label{eq:sp}
\end{equation}
with the initial condition
\begin{equation}
x(0)=x_0\in \mathcal{C}(0),
\label{eq:ic}
\end{equation}
where $N_{\mathcal{C}(t)}(x)$ is the outward normal cone from Definition \ref{def:normal_cone}. Equation \eqref{eq:sp} is understood as being held for {\it almost all} $t\in[0,T]$ (e.g. $\dot x$ of Fig. \ref{fig:sweeping} b is undefined at $t=t_3$ where $x$ just reached the corner). 

\begin{figure}[h]\center
\includegraphics[scale=0.5]{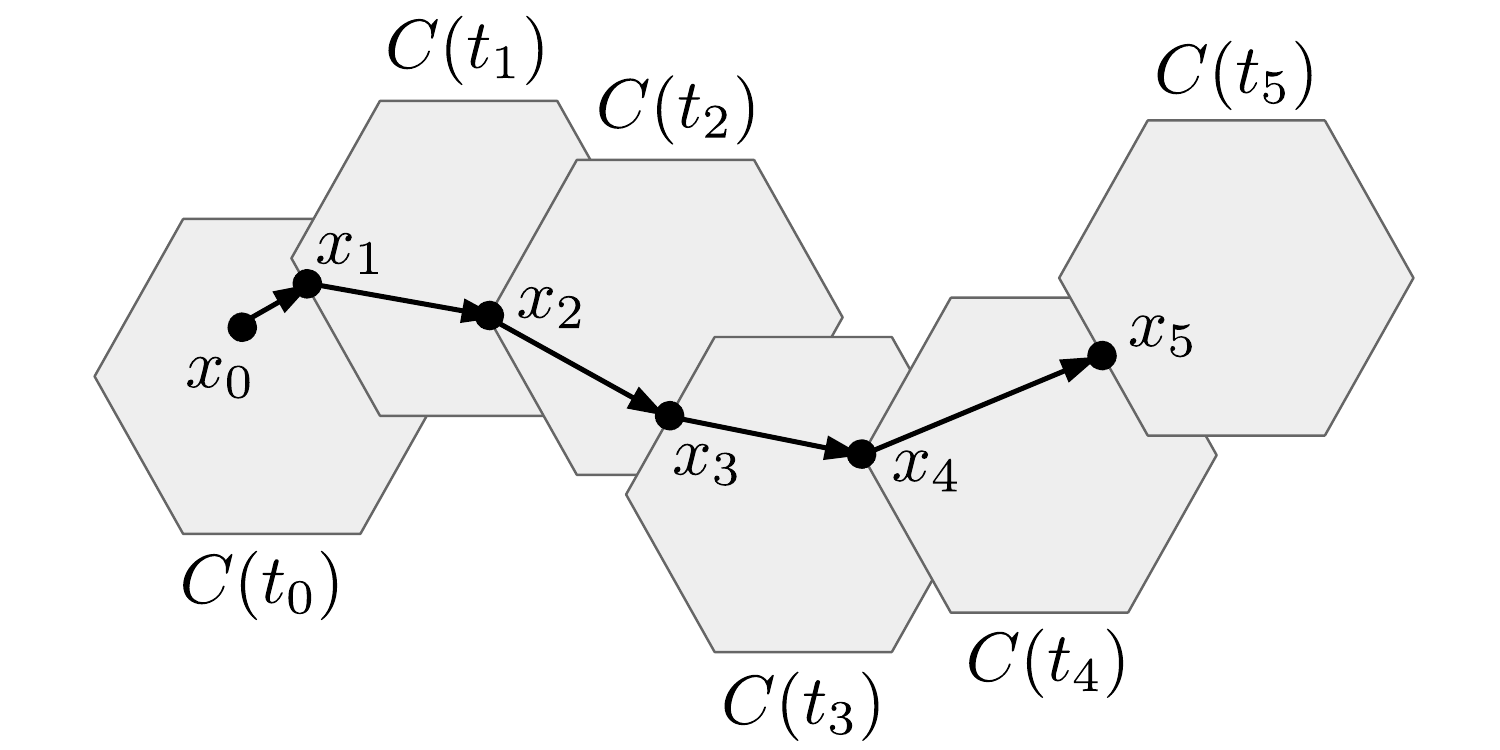}
\caption{
\footnotesize A geometric interpretation of the catch-up algorithm.
} \label{fig:catch-up}
\end{figure}

{\color{black}
While problem \eqref{eq:sp}-\eqref{eq:ic} may appear discouraging when looked at as a differential equation $\dot x = f(x,t)$ with discontinuous set-valued and unbounded right-hand side $f$ (see the evolution of $N_{\mathcal{C}(t)}(x)$ in Fig. \ref{fig:sweeping} b), the properties of the normal cone yield the well-posedness of the problem, under a reasonable regularity assumption on $\mathcal{C}(t)$:
\begin{Theorem}
\label{th:sp_abstract_existence_uniqueness}
Let $\mathcal{C}$ be a set-valued function of the form \eqref{eq:moving_set_sweeping} with a closed, convex, nonempty value for each $t\in[0,T]$. Moreover, let $\mathcal{C}$ be Lipschitz-continuous with respect to the Hausdorff distance, i.e. there is $L>0$ such that for any $t_1, t_2\in[0,T]$
\[
\max\left(\sup_{x\in \mathcal{C}(t_2)}{\rm dist}(x, \mathcal{C}(t_1)),\,\sup_{x\in \mathcal{C}(t_1)}{\rm dist}(x, \mathcal{C}(t_2))\right)\leqslant L|t_1-t_2|,
\]
where $\rm dist$ is defined by \eqref{eq:dist_def}.
Then there exists a unique solution to \eqref{eq:sp}-\eqref{eq:ic}.
\end{Theorem}

The proof of Theorem \ref{th:sp_abstract_existence_uniqueness} is already considered classical (see e.g. \cite{Kunze2000},\cite[Sect. 5]{Moreau1973}), and we will only explain the general approach. The uniqueness follows from the {\it monotonicity} of the normal cone, \cite[Th. 3]{Kunze2000}. To show the existence one considers a sequence of approximating problems, extracts a limiting function of solutions to the approximating problems and demonstrates that the limit satisfies the original problem  \eqref{eq:sp}-\eqref{eq:ic}. One way to construct such approximating problems is } a so-called {\it catch-up algorithm}, in which the approximating trajectory is found via consecutive projections (see Fig. \ref{fig:catch-up}).
 
Specifically, for a partition $0=t_0<t_1<\dots< t_{k-1}<t_k=T$ of the time interval $[0,T]$ into $k$ segments,  the catch-up algorithm computes
\begin{equation}
x_{i+1}={\rm proj}(x_{i}, \mathcal{C}(t_{i+1})). \qquad i\in\overline{0,k-1}
\label{eq:catch-up}
\end{equation}
As $\max\{t_i-t_{i-1}:i\in \overline{1,k}\}\to 0$, the values $x_i$ uniformly approach the corresponding values $x(t_i)$ of the exact solution of \eqref{eq:sp} (see e.g. \cite[\S 5.h]{Moreau1973}). 

{\color{black}Informally, it is intuitive to view \eqref{eq:catch-up} as an ``Euler step'' of \eqref{eq:sp} with $\dot x \approx \frac{x_{i+1}-x_{i}}{t_{i+1}-t_i}$ :
\[
-\frac{x_{i+1}-x_{i}}{t_{i+1}-t_i}\in N_{\mathcal{C}(t_{i+1})}(x_{i+1})
\]
\[
x_{i}- x_{i+1}\in (t_{i+1}-t_i) N_{\mathcal{C}(t_{i+1})}(x_{i+1}).
\]
Because of property \eqref{eq:normal_cone_is_a_cone} and $t_{i+1}-t_i>0$ we have
\[
x_{i}- x_{i+1}\in N_{\mathcal{C}(t_{i+1})}(x_{i+1}),
\]
which is equivalent to \eqref{eq:catch-up} due to relation  \eqref{eq:nc_proj_relation}.
}

For a given symmetric positive definite $n\times n$ matrix $S$, the sweeping process can be defined in the sense of weighted inner product \eqref{eq:S-inner-prod} using normal cone \eqref{eq:S-normal_cone}. We write the corresponding sweeping process as
\begin{equation}
-\dot x \in N^S_{\mathcal{C}(t)}(x).
\label{eq:S-sp}
\end{equation} 
The catch-up algorithm for such sweeping process is \eqref{eq:catch-up} with the projection replaced by \eqref{eq:S-proj}. The algorithm does not only help to prove existence of the solution, but it can also be used as a practical numerical scheme, to which we refer as Algorithm \ref{alg:catch-up1}.


\begin{algorithm}[H]
\SetAlgoLined
\SetKwComment{Comment}{//}{}
\Comment{Given $\mathcal{C}(t)$  as in Theorem \ref{th:sp_abstract_existence_uniqueness}}
\Comment{and a partition $0=t_0<t_1<\dots< t_{k-1}<t_k=T$}
 \For{$i:=0$ \KwTo $k-1$}{
  \Comment{the projection of the type \eqref{eq:S-proj}:}
  $x_{i+1} := {\rm proj}^S(x_i, \mathcal{C}(t_{i+1}))$\;
   }  
 \caption{Catch-up algorithm for general sweeping process \eqref{eq:ic},\eqref{eq:S-sp}}
 \label{alg:catch-up1}
\end{algorithm}

\section{Stresses in the Lattice Spring Model via a sweeping process}
\label{sect:sweeping_solution}
It turns out that the evolution of vector of elastic elongations $\varepsilon$ (and, respectively, vector of stresses $\sigma$) governed by equations \eqref{eq:gc}-\eqref{eq:balanceLaw} boils down to a sweeping process of the type \eqref{eq:S-sp}. The derivation of sweeping process from governing equations is based on the ideas of J.-J. Moreau \cite{Moreau1973} and it is similar to \cite[Th. 3.1]{GM-1}, where one-dimensional lattices are considered. In the current section we explicitly derive the sweeping process from \eqref{eq:gc}-\eqref{eq:balanceLaw}, discuss the numerical schemes to solve it and present examples of lattices that are solved for stresses via this approach.

\subsection{Derivation of the sweeping process from the governing equations of a LSM}
\label{ssect:sp-derivation}

{\color{black}The sweeping process itself is defined in the space $\mathbb{R}^m$ equipped with weighted inner product \eqref{eq:S-inner-prod} with $S=K$. We write the sweeping process and the initial condition as
\begin{equation}
\begin{cases}
-\dot y \in  N^K_{\mathcal{C}(t)}(y),\\
y(0)=y_0,
\end{cases}
\label{eq:elastoplastic_sp}
\end{equation}
where the right-hand side of the inclusion is the normal cone \eqref{eq:S-normal_cone} defined in accordance with the inner product of the space. Loosely speaking, moving set $\mathcal{C}(t)$ represents the interplay between the affine constraint of  \eqref{eq:balanceLaw} and the unilateral constraint
\[
\sigma \in C
\]
implied by \eqref{eq:plasticity}, see Remark \ref{rem:hidden_inclusion}. In turn, the {\it sweeping variable } $y$ is directly tied to the {\it yielding variables} $\sigma$ and $\varepsilon$ via a change of variables. Specifically, we set
\begin{equation}
\mathcal{C}(t)=\left(K^{-1}C+Gr(t)-Ff(t)\right)\cap \mathcal{V},
\label{eq:elastoplastic_moving_set}
\end{equation}
\begin{equation}
y=\varepsilon+G r(t)- F f(t),
\label{eq:y_from_e}
\end{equation}
where $G$ and $F$ are known matrices constructed as described in Table \ref{table:matrices}. The proof of the following theorem shows how these quantities emerge step by step.
\begin{table}[h]
\begin{center}
\begin{tabular}{| r | l |}
\hline & \\[-3mm]
$R^+= R^T(RR^T)^{-1}$ & the $nd \times q$ Moore-Penrose pseudoinverse of $R$, \\ & where the equality is due to \eqref{eq:rowsR_indep} and Proposition \ref{prop:MP3}, \\[1mm]
\hline  & \\[-3mm]
$R_0$ & the matrix composed of columns, \\ & which  form a basis in the nullspace of $R$,\\[1mm]
\hline  & \\[-3mm]
$U = (D_{\xi_0} \varphi) R_0$ & the columns of the matrix form a basis in $\mathcal{U}$ \\ & due to Proposition \ref{prop:U_and_V_orth} {\it iii} \\[1mm]
\hline & \\[-3mm]
$V$ & the matrix composed of columns which form a basis in $\mathcal{V}$\\[1mm]
\hline & \\[-3mm]
$P_U = \left(U^T K U\right)^{-1}U^T K$ &  matrices of orthogonal in the sense of \eqref{eq:orthogonality_K} projections  \\
$P_V = \left(V^T K V\right)^{-1}V^T K$ &  on $\mathcal{U}$ and $\mathcal{V}$ respectively, expressed in terms of coordinates, \\
&in bases $U$ and $V$, see \cite[Sect. 5.3 and 5.4]{OlverLinearAlgebra2018}\\[1mm]
\hline & \\[-3mm]
$UP_U$ & $m\times m$ matrices of orthogonal in the sense of \eqref{eq:orthogonality_K} projections \\
$VP_V$ & on $\mathcal{U}$ and $\mathcal{V}$ respectively \\[1mm]
\hline & \\[-3mm]
$G=V\, P_V\,(D_{\xi_0}\varphi)\, R^+$ & \\
$F=U\,P_U\,K^{-1}\,H$ &\\[1mm]
\hline
\end{tabular}
\caption{Technical quantities to construct the sweeping process \eqref{eq:elastoplastic_sp}-\eqref{eq:elastoplastic_moving_set} with the unknown variable \eqref{eq:y_from_e}.}
\end{center}
\label{table:matrices}
\end{table}
}
\begin{Theorem} Consider a Lattice Spring Model given by an $n\times m$ incidence matrix $Q$, $m\times m$ positive diagonal matrix $K$ of stiffness coefficients, $q\times nd$ matrix $R$ (which describes the external displacement constraint), elasticity limits $c^-, c^+ \in \mathbb{R}^m$, a reference configuration $\xi_0\in \mathbb{R}^{nd}$ and Lipschitz-continuous functions $r:[0,T]\to \mathbb{R}^q$(displacement load) and $f:[0,T]\to \mathbb{R}^{nd}$(stress load), such that Assumptions \ref{ass:rowsR_indep}-\ref{ass:at_least_one_self-stress} hold.

If there are Lipschitz-continuous functions $\zeta:[0,T]\to\mathbb{R}^{nd},  x,\varepsilon, p,\sigma :[0,T]\to \mathbb{R}^m$ satisfying the system  \eqref{eq:gc}-\eqref{eq:balanceLaw}, then function $y$ given by \eqref{eq:y_from_e} solves sweeping process \eqref{eq:elastoplastic_sp}-\eqref{eq:elastoplastic_moving_set} with initial condition 
\begin{equation}
y_0=K^{-1}\sigma_0+G r(0) - F f(0).
\label{eq_y_0}
\end{equation}
\label{th:lsm_to_sp}
\end{Theorem}
\noindent{\bf Proof.} Matrices $U$ and $V$ are constructed so that their columns form bases in, respectively, $\mathcal{U}$ and $\mathcal{V}$. In turn, as shown in \cite[Sect. 5.3 and 5.4]{OlverLinearAlgebra2018}, matrices $P_U$ and $P_V$  applied to a vector from $\mathbb{R}^m$ give coordinates of the vector's orthogonal projections onto $\mathcal{U}$ and $\mathcal{V}$ in bases $U$ and $V$ respectively, where orthogonality is meant in the sense of the weighted inner product \eqref{eq:S-inner-prod} with $S=K$. Thus $U P_U$ and $V P_V$ is a pair of corresponding orthogonal projection matrices and this means, in particular, that 
\begin{equation}
a=UP_U a+ VP_V a\in \mathcal{U} +V P_V a \qquad \text{ for any } a\in \mathbb{R}^m.
\label{eq:projections_decomposition}
\end{equation}

Now we derive the sweeping process \eqref{eq:elastoplastic_sp}-\eqref{eq:elastoplastic_moving_set} from \eqref{eq:gc}-\eqref{eq:balanceLaw}.
Fix a.e. $t\in[0,T]$, take time-derivative of \eqref{eq:bc} and apply the pseudoinverse matrix $R^+$ of $R$:
\[
R^+ R \dot \zeta + R^+\dot r(t)=0.
\]
By Proposition \ref{prop:MP2} matrix $R^+R$ is the orthogonal projection matrix onto ${\rm Im}\, R^T$, therefore there is $z\in ({\rm Im}\, R^T)^\perp={\rm Ker}\, R$, such that $\dot \zeta = R^+ R \dot \zeta + z$ and 
\[
\dot \zeta -z +R^+\dot r(t)=0.
\]
We then apply $D_{\xi_0} \varphi$ to get,
\[
(D_{\xi_0} \varphi)\dot \zeta - (D_{\xi_0} \varphi) z + (D_{\xi_0} \varphi)R^+\dot r(t)=0,
\]
and from \eqref{eq:gc}, \eqref{eq:space_u_def} we deduce that
\[
\dot x \in \mathcal{U} -  (D_{\xi_0} \varphi)R^+\dot r(t).
\]
Due to \eqref{eq:projections_decomposition}
\[
\dot x \in \mathcal{U} -  VP_V(D_{\xi_0} \varphi)R^+\dot r(t),
\]
i.e.
\begin{equation}
\dot x \in \mathcal{U} - G\dot r(t) 
\label{eq:x-in-U-G}
\end{equation}
with $G=VP_V(D_{\xi_0} \varphi)R^+$ as defined earlier.

Now consider equation of equilibrium \eqref{eq:balanceLaw} and its equivalent form in \eqref{eq:equiv_balance_law}. By applying $K^{-1}$ to its both sides and using \eqref{eq:hooke}, \eqref{eq:space_v_def} we get
\[
\varepsilon-K^{-1}H f(t) \in-\mathcal{V}=\mathcal{V}.
\]
Using \eqref{eq:projections_decomposition} we rewrite the latter equality as
\[
\varepsilon-UP_UK^{-1}Hf(t)\in \mathcal{V},
\]
or 
\[
\varepsilon - Ff(t) \in \mathcal{V},
\]
with $F=UP_UK^{-1}H$ as defined earlier. Here we use the change of variables \eqref{eq:y_from_e} to get 
\[
y - Gr(t) \in \mathcal{V},
\]
and, since $Gr(t)\in \mathcal{V}$ by construction of $G$,
\begin{equation}
y \in \mathcal{V}.
\label{eq:y_in_V}
\end{equation}

Now we put together all the parts to obtain the sweeping process. From \eqref{eq:x} we have
\[
\dot x =\dot \varepsilon + \dot p,
\]
i.e.
\[
-\dot \varepsilon =\dot p -\dot x,
\]
and by \eqref{eq:plasticity} and \eqref{eq:x-in-U-G}
\[
-\dot \varepsilon \in N_{C}(\sigma) -\mathcal{U}+ G\dot r(t).
\]
It follows from the definition of the normal cone \eqref{eq:normal_cone} and \eqref{eq:hooke} that 
\[
N_C(\sigma)=N_C(K\varepsilon)=N^K_{K^{-1}C}(\varepsilon).
\] 
Therefore
\[
-\frac{\rm d}{\rm dt}(\varepsilon + Gr(t))\in N^K_{K^{-1}C}(\varepsilon) -\mathcal{U}.
\]
For $y$ and $\varepsilon$, which are connected via the change of variables \eqref{eq:y_from_e}, we have 
\begin{align*}
-\frac{\rm d}{\rm dt}(y+Ff(t))&\in N_{K^{-1}C}^K(y-Gr(t)+Ff(t))-\mathcal{U},\\
-\frac{\rm d}{\rm dt}(y) &\in N_{K^{-1}C}^K(y-Gr(t)+Ff(t))-\mathcal{U}+F\dot f(t),
\end{align*}
Observe that $F\dot f(t)\in \mathcal{U}$ (by construction of $F$) and $\mathcal{U}=\mathcal{-U}$ (as a linear space). Hence, due to \eqref{eq:y_in_V} and orthogonality of $\mathcal{U}$ and $\mathcal{V}$ in the  sense of weighted inner product \eqref{eq:S-inner-prod} with $S=K$, we have 
\[
-\mathcal{U}+F\dot f(t) = \mathcal{U}=N^K_\mathcal{V}(y).
\]
Therefore,
\[
-\dot y\in N_{K^{-1}C}^K(y-Gr(t)+Ff(t))+N^K_\mathcal{V}(y),
\]
Again, from the definition of the normal cone \eqref{eq:S-normal_cone} it can be shown that $N^K_C(a+b)=N^K_{C-b}(a)$ for any $a,b\in \mathbb{R}^m$, so we have
\[
-\dot y\in N_{K^{-1}C+Gr(t)-Ff(t)}^K(y)+N^K_{\mathcal{V}}(y).
\]
Finally, the additive property of the normal cone to polyhedral sets \cite[Corollary 23.8.1]{Rockafellar} yields 
\[
-\dot y\in N^K_{\left(K^{-1}C+Gr(t)-Ff(t)\right)\cap \mathcal{V}}\,(y),
\]
which is the same as \eqref{eq:elastoplastic_sp}-\eqref{eq:elastoplastic_moving_set}
$\blacksquare$

\begin{figure}[H]\center
\includegraphics{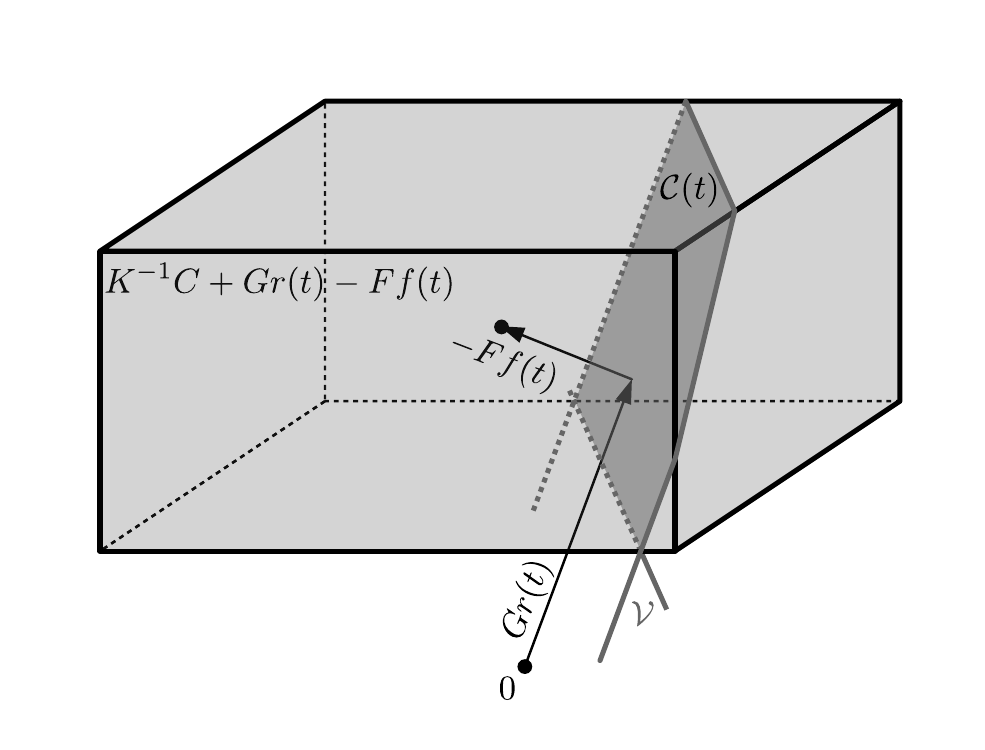}
\caption{
\footnotesize A schematic representation of the moving set $\mathcal{C}(t)$ in the sweeping process \eqref{eq:elastoplastic_sp}-\eqref{eq:elastoplastic_moving_set} and its dependence on the displacement load $r(t)$ and external forces $f(t)$. The positions of the origin $0$ and hyperplane $\mathcal{V}$ are fixed, while the hyperrectangle $K^{-1}C+Gr(t)-Ff(t)$ moves by translation.
} \label{fig:moving-set}
\end{figure}

\begin{Remark}
Notice that by construction of $G$, we have $Gr(t)\in \mathcal{V}$ (see Fig. \ref{fig:moving-set}), so \eqref{eq:elastoplastic_moving_set} can be rewritten as 
\[
\mathcal{C}(t) = \left(K^{-1}C -Ff(t)\right)\cap\mathcal{V}+ Gr(t).
\]
Therefore, a change {\color{black} of displacement load} $r(t)$ results only in a {\it translation} of the whole moving set $\mathcal{C}(t)$ in $\mathcal{V}$ without any change of its {\it shape}. In contrast, $Ff(t)\notin\mathcal{V}$ (except when $Ff(t)=0$), so that a change of external force $f(t)$ affects the shape of the intersection $\left(K^{-1}C -Ff(t)\right)\cap\mathcal{V}$, i.e. it affects the shape of $\mathcal{C}(t)$.
\end{Remark}

\begin{Remark}
\label{remark:slc1}
As mentioned in Section \ref{sect:sweeping_intro}, {\color{black} for a sweeping process to have a solution, its moving set $\mathcal{C}(t)$ must} be nonempty for all $t\in[0,T]$. In particular, the requirement for the $\mathcal{C}(t)$ in \eqref{eq:elastoplastic_moving_set} to be nonempty, {\color{black} or, equivalently, the requirement
\begin{equation}
\left(K^{-1}C -Ff(t)\right)\cap\mathcal{V}\neq \varnothing
\label{eq:slc1}
\end{equation}
 is called the {\it safe load condition}, \cite[(3.3)-(3.4)]{DalMaso2006}, \cite[Sect. 6b, Assumpt. 3]{Moreau1973}. From the physical point of view, the nonempty intersection means that external forces $f(t)$ can be potentially balanced by the internal stresses in the lattice. Vice versa,} the safe load condition is violated when the magnitude of $f(t)$ is too big and the equilibrium \eqref{eq:balanceLaw} cannot be  achieved with all individual stresses staying within their elasticity bounds. This corresponds to the geometric fact, observable from Fig. \ref{fig:moving-set}: if the magnitude of $f(t)$ is large enough, the hyperrectange $K^{-1}C+Gr(t)-Ff(t)$ will no longer intersect with the hyperplane $\mathcal V$.
\end{Remark}

\begin{Remark}In Theorem \ref{th:lsm_to_sp} we established a one-sided implication saying that for every solution of \eqref{eq:gc}-\eqref{eq:balanceLaw} there is a corresponding solution of sweeping process \eqref{eq:elastoplastic_sp}-\eqref{eq:elastoplastic_moving_set}. It is also possible to prove the reversed statement, where for every solution $y$ of the sweeping process one can not only find $\varepsilon$ and $\sigma$ (trivially obtained from \eqref{eq:y_from_e}, \eqref{eq:hooke}), but also obtain trajectories for $p,x$ and $\zeta$ satisfying \eqref{eq:gc}-\eqref{eq:balanceLaw} altogether. Specifically, one can construct the differential inclusion similar to \cite[(3.12)]{GM-1} to get $p$, use \eqref{eq:x} to get $x$ and then take advantage of kinematic determinacy, which allows to solve \eqref{eq:gc}-\eqref{eq:bc} for $\zeta$. However, in the situation of perfect plasticity the trajectory of plastic elongation $p$ is, generally, not defined uniquely. Instead, a continuum of possible trajectories of $p$ exists (see e.g. \cite{Meyer2}), and some of those trajectories develop shear bands, where shear deformation concentrates \cite[p. 239]{DalMaso2006}. Reliable numerical computation of a possible trajectory of plastic elongation $p$ in the context of the sweeping process approach is a nontrivial task which requires special attention and is beyond the scope of the current paper. Instead, we focus on obtaining stress trajectories by using the current formulation of Theorem \ref{th:lsm_to_sp}. Such treatment of stress trajectory alone is called ``reduced solution'' \cite{Meyer} and the corresponding problem is called ``the dual problem'' or  ``the reduced form of the problem'' \cite[Sect. 8.2]{han}.
\end{Remark}

\subsection{The relation between rigidity properties of the lattice and its corresponding sweeping process}
\label{ssect:msp-rigitity}
One of the questions coming from the construction of the sweeping process as described in this paper, is how the sweeping process depends on the underlying graph structure of the lattice. In particular, it is critical to understand what determines the dimension of subspace $\mathcal{V}$. 
Let us {\color{black} summarize} some implications of the rigidity theory, {\color{black} which we discussed in Sections \ref{ssect:rigidity_properties1}, \ref{ssect:additional_constraint}, \ref{ssect:statics_gc_only}, \ref{ssect:static_properties_with_bc} and \ref{ssect:fundamental_spaces}}, for our construction of the sweeping process:
\noindent\begin{itemize}
\item The space $\mathcal{V}$, which contains the sweeping process we constructed, has a clear physical interpretation {\color{black} (see Proposition \ref{prop:U_and_V_orth} {\it i}) and its dimension coincides with the number of states of self-stresses in the lattice}.
\item It is convenient to require the external displacement constraint to be tight enough for the lattice to be kinematically determinate (Assumption \ref{ass:kinematic_determinacy}), as such a requirement ensures the following properties:
\begin{itemize}
\item every stress load is resolvable, so we don't have to restrict stress load for the well-posedness of the problem (apart from the safe load condition{\color{black}, which comes from the limitations of perfect plasticity, see Remark \ref{remark:slc1}}),
\item there is an explicit formula \eqref{eq:dimensions_U_V} for the dimension ${\rm dim}\, \mathcal{V}$ of the sweeping process and, since kinematic determinacy means that matrix $\begin{pmatrix} (D_{\xi_0} \varphi)^T & R^T\end{pmatrix}$ has right inverse, matrix $H$ in \eqref{eq:balanceLaw} can also be expressed by an explicit formula (the last equality in \eqref{eq:new_equilibrium_matrix_pseudoinverse}),
\item it corresponds to the physically correct situation in a real experiment when a specimen is properly secured with no freely moving parts,
\item while it is beyond the scope of the current paper, in a kinematically determinate lattice  the displaced positions of the nodes are uniquely defined by elastic and plastic elongations, i.e. it is possible to compute $\zeta$ from $\varepsilon+p$.
\end{itemize}
\item Condition \eqref{eq:at_least_one_self-stress} excludes the case of a statically determinate lattice which leads to a degenerate sweeping process with $\mathcal{C}(t)$ being a single point set.
\item Structural mechanics and rigidity theory can provide valuable results, which can be used to estimate the 
dimension of the sweeping process a priory and help design the appropriate additional constraint satisfying conditions \eqref{eq:rowsR_indep}, \eqref{eq:at_least_one_self-stress}. For example, if we know a priory that the lattice is infinitesimally rigid, then it is enough to have constraint \eqref{eq:bc} with $q=\frac{d(d+1)}{2}$. Such a priory results are classical for lattices with nodes and edges {\color{black} placed at vertices and edges of a convex polyhedron} \cite{Dehn1916, Roth1981, alexandrov_book,AsimowRothRigidityII} and for lattices in $\mathbb{R}^2$ (Laman's theorem, see e.g. \cite{CombinatorialRigidityBook}). For periodic graphs the analysis of zero modes is a topic of current interest in science \cite{Lubensky2015, Borcea2021}, which is stimulated, in particular, by  applications in crystallography, see e.g. \cite{Giddy1993}.
\end{itemize}

\begin{Remark}
When there is only one spatial dimension ($d=1$) the question of determination of zero modes is trivial compared to the general case. If all the springs are chosen to be oriented in the same direction, equilibrium matrix $ (D_{\xi_0} \varphi)^T$ becomes just the incidence matrix of the graph of springs (up to the sign) and the set of states of self-stresses becomes the cycle space of the graph, see \cite{bapat} for definitions. Then the number of zero modes equals the number of connected components of the graph (it follows from \cite[Th. 2.3]{bapat} and the rank-nullity theorem), the set of resolvable loads is described by just the first equation of \eqref{eq:sum_of_forces_and_torques} (which is scalar when $d$=1, see \cite[Lemma 2.4]{bapat}) and rigidity is equivalent to connectedness of the graph \cite[Prop. 1.1.2]{Whiteley1992}. Detailed studies of \eqref{eq:gc}-\eqref{eq:balanceLaw} in the case of a single spatial dimension can be found in our previous works \cite{GM-1, GM-2, GMR}.
\end{Remark}

\subsection{Catch-up algorithm for the sweeping process coming from the Lattice Springs Model}
\label{ssect:catch-up_elastoplastic}
To use the catch-up algorithm on the sweeping process \eqref{eq:elastoplastic_sp}-\eqref{eq:elastoplastic_moving_set} we explicitly rewrite the moving set \eqref{eq:elastoplastic_moving_set} in the form \eqref{eq:polyhedral_set}
as 
\begin{equation}
\mathcal{C}(t)=\left\{x\in \mathbb{R}^m: \begin{array}{c} K^{-1}c^- + Gr(t)-F f(t)\leqslant x \leqslant K^{-1}c^+  +Gr(t)-F f(t),\\ U^T K x=0 \end{array}\right\},
\label{eq:elastoplastic_moving_set2}
\end{equation}
where
\[
c^-=(c_i^-)_{i\in \overline{1,m}}, \qquad c^+=(c_i^+)_{i\in \overline{1,m}},
\]
i.e. in terms of \eqref{eq:polyhedral_set} we have
\begin{equation}
\begin{array}{ll}
A=\begin{pmatrix} I_{m\times m}\\ -I_{m\times m}\end{pmatrix},&
b(t)=\begin{pmatrix}K^{-1}c^+ + Gr(t)-F f(t) \\ -\left(K^{-1}c^- + Gr(t)-F f(t)\right)\end{pmatrix}\\[0.5cm]
A_{eq}= U^T K,&
b_{eq}(t) = 0\in \mathbb{R}^{{\rm dim}\, \mathcal{U}}.
\end{array}
\label{eq:ABinRm}
\end{equation}
Algorithm \ref{alg:catch-up-elastoplastic1} below is an adaptation of the catch-up algorithm (Algorithm \ref{alg:catch-up1}) to solve for stresses in a Lattice Spring Model, in which the unknown variables $\varepsilon, \sigma$ of the lattice are obtained via
\begin{equation}
\varepsilon(t)=y(t)-G r(t)+ F f(t),\qquad \sigma(t) = K \varepsilon(t).
\label{eq:e_s_from_y}
\end{equation}

\begin{algorithm}[H]
\SetAlgoLined
\SetKwComment{Comment}{//}{}
\Comment{Given $\sigma_0, K, G, F, r, f$ and $\mathcal{C}(t)$ as \eqref{eq:elastoplastic_moving_set2}-\eqref{eq:ABinRm}}
\Comment{and a partition $0=t_0<t_1<\dots< t_{k-1}<t_k=T$}
 $y_0:=K^{-1}\sigma_0+Gr(0)-Ff(0)$\;
 \For{$i:=0$ \KwTo $k-1$}{
  \Comment{the projection of the type \eqref{eq:S-proj-polyhedral} with constraints \eqref{eq:ABinRm}:}
  $y_{i+1} := {\rm proj}^K(y_i, \mathcal{C}(t_i))$\;
  \Comment{recover the elastic elongations and the stresses from the solution of the sweeping process:}
  $\varepsilon_{i+1}:=y_{i+1}-Gr(t_{i+1})+Ff(t_{i+1})$\;
  $\sigma_{i+1}:=K\varepsilon_{i+1}$
   }  
 \caption{Catch-up algorithm to compute stresses via the sweeping process \eqref{eq:elastoplastic_sp}-\eqref{eq:ABinRm}}
 \label{alg:catch-up-elastoplastic1}
\end{algorithm}

In Appendix \ref{sect:appendix_reduced_dimension} we discuss an equivalent sweeping process of reduced dimension, which leads to a significant raise in the performance of numerical algorithms.
Also, in Section \ref{sect:leapfrog} we will discuss an event-based approach, which allows to skip a large number of non-essential steps of the catch-up algorithm in case when the displacement load $r$ is a piecewise linear function and the stress load $f$ is constant.
{\color{black} 
\begin{Remark}
\label{remark:slc2}
We can add to Remark \ref{remark:slc1} that, on top of physical and geometric interpretations, from the numerical point of view the violation of the safe load condition means that the corresponding projection step of the catch-up algorithm is an infeasible quadratic programming problem. Using \eqref{eq:elastoplastic_moving_set2}, the safe load condition \eqref{eq:slc1} can be written in a more explicit form
\begin{equation*}
\text{there is $\sigma\in \mathbb{R}^m$ such that} \quad c^--KFf(t)\leqslant \sigma\leqslant c^+ -KFf(t)\quad\text{ and }\quad U^T \sigma=0.
\label{eq:slc2}
\end{equation*}

\end{Remark}
}

\begin{Remark}
\label{remark:reduction_of_dimensions}
As one can see from the construction of \eqref{eq:elastoplastic_moving_set} and Fig. \ref{fig:moving-set}, the unknown variable $y$ never leaves linear subspace $\mathcal{V}\subset \mathbb{R}^m$, which is independent of time. In other words, the equality constraint in \eqref{eq:elastoplastic_moving_set2} is independent on time. Using basis $V$ in $\mathcal{V}$, we can derive a fully equivalent sweeping process in the space of coordinates $\mathbb{R}^{{\rm dim}\, \mathcal{V}}$. We give the derivation in Appendix \ref{sect:appendix_reduced_dimension}, together with the corresponding numerical schemes including the catch-up algorithm (Appendix \ref{ssect:catch-up method for lattices in V}).

Computation in the space of reduced dimension is much cheaper: e.g. for a lattice of $m=496$ springs with ${\rm dim}\, \mathcal{V} = m-nd+q = 496 - 198 \cdot 2+56=156$ the catch-up algorithm in the space of reduced dimension $\mathbb{R}^{{\rm dim}\, \mathcal{V}}$ took about seven times less time then the direct implementation of Algorithm \ref{alg:catch-up-elastoplastic1} in $\mathbb{R}^m$. From our observations, this is a typical performance increase for lattices of such size when passing to reduced dimensions. We will examine the above-mentioned lattice in Section \ref{sect:example_with_the_hole}. 

On top of the performance boost, for some lattices it is possible to visualize the sweeping process in $\mathbb{R}^{{\rm dim} \mathcal{V}}$, but not in $\mathbb{R}^m$. In particular, this is the case for Example 1 with ${\rm dim}\, \mathcal{V}=2$ and $m=10$.
\end{Remark}

\subsection{Example 1 continued.}
\begin{figure}[H]\center
\includegraphics{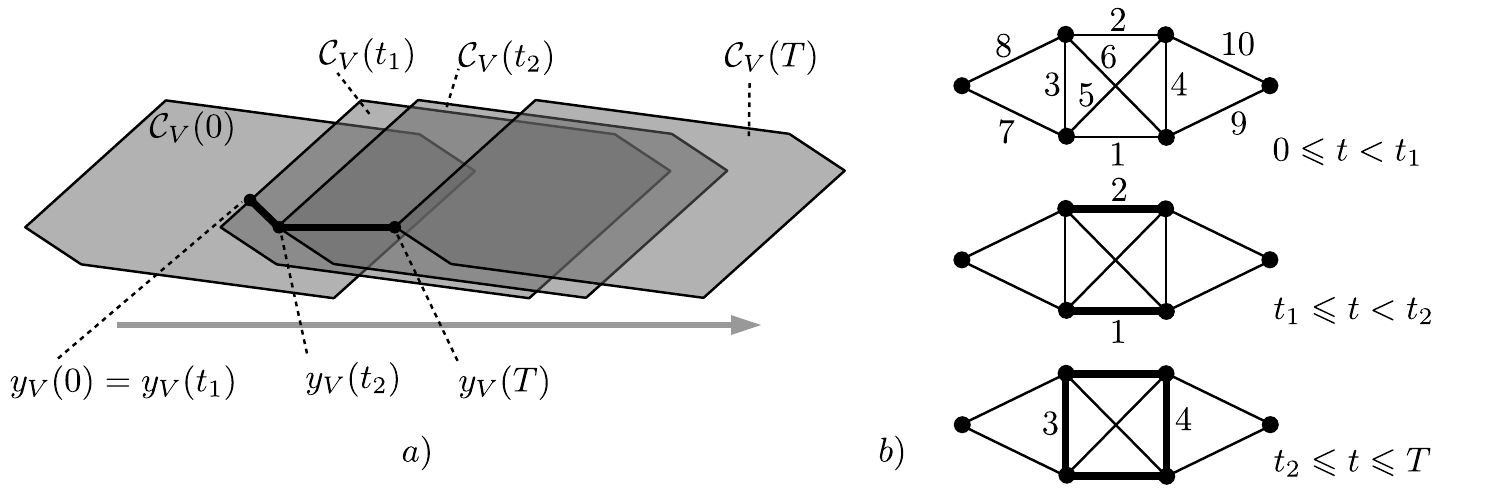}
\caption{
\footnotesize a) The representation of the sweeping process corresponding to Example 1 with zero initial condition and zero stress load on the interval $[0, T]=[0, 0.08]$. The moving set is shown at instructional moments: at $t=0$ we have $y(0)=0$ at the center of the moving set, at $t=t_1=0.042$ the solution meets the edge of the moving set and starts to slide along the edge, until it reaches the vertex at $t=t_2=0.055$. For the rest of  time it is swept by that vertex. b) We indicate yielding springs at different stages of the evolution.  The moments $t_1$ and $t_2$, when the solution of the sweeping process meets a new edge of the moving set, correspond to new springs (bold lines) reaching the plastic stage of deformation.
} \label{fig:moving-set-example}
\end{figure}

\begin{figure}[H]\center
\includegraphics{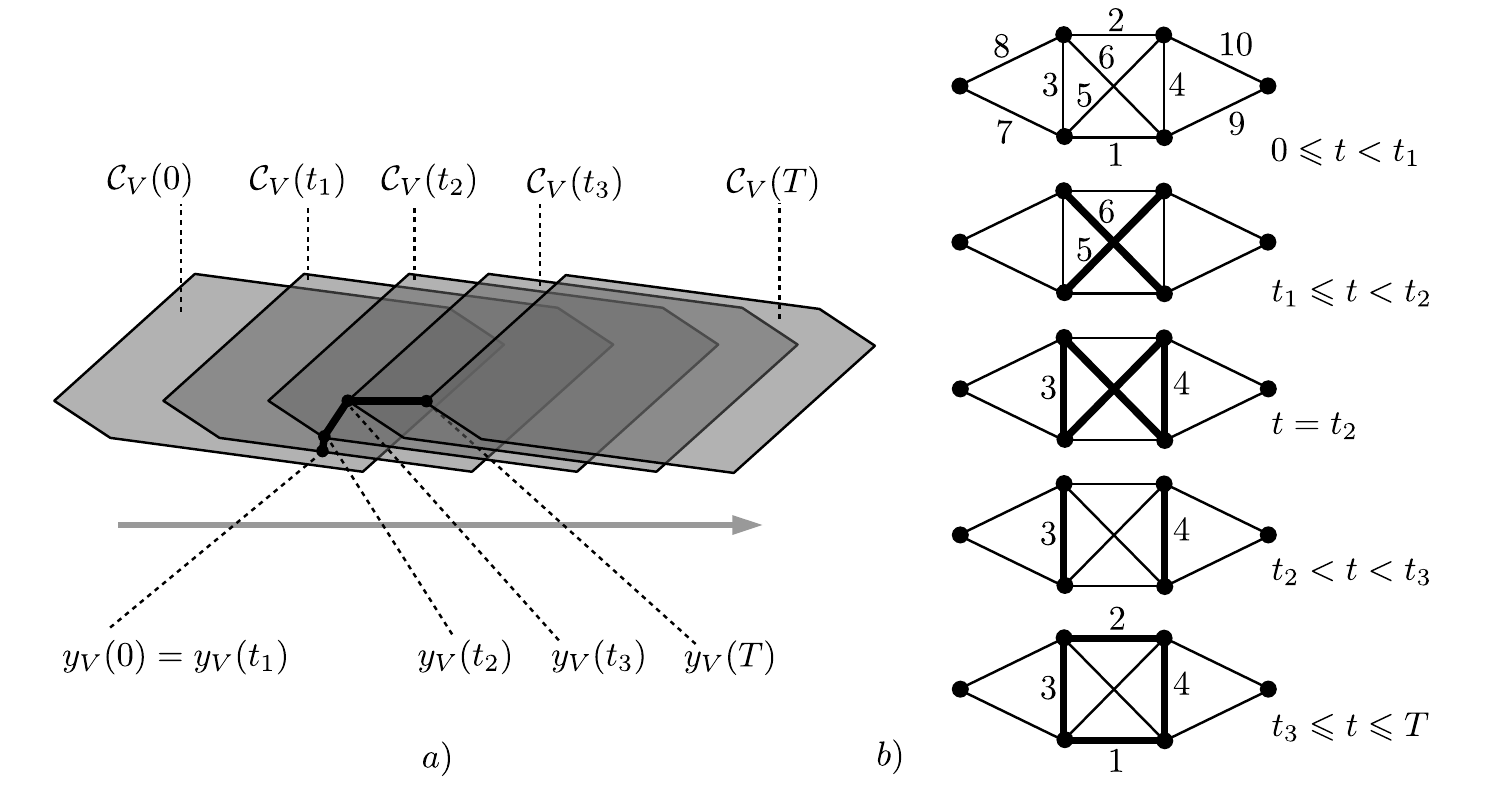}
\vspace{-3mm}
\caption{ \footnotesize  a) The representation of the sweeping process corresponding to Example 1 with the initial condition \eqref{eq:alt_initinal_cond} and with zero stress load on the interval $[0, T]=[0, 0.08]$. The moving set is shown at important moments:  at $t=0$ we have $y_V(0)$ in the interior of $\mathcal{C}_V(0)$, at $t=t_1=0.027$ the solution meets the edge of the moving set and starts to slide along the edge, until it reaches the vertex at $t=t_2=0.046$. However, the solution then leaves the vertex and continues to slide along the new edge, until at $t=t_3=0.064$ it finally reaches the same vertex as the sweeping process in Fig. \ref{fig:moving-set-example}. The fact that the solutions always arrive at the same final vertex (corresponding to the same distribution of stresses) independently of the initial condition (provided that the magnitude of the displacement loading is large enough) is a simple illustration of the {\it shakedown} phenomenon. b) The moments $t_1, t_2$ and $t_3$ when the solution of the sweeping process meets a new edge of the moving set correspond to new springs (bold lines) reaching plastic stage of deformation. Notice how diagonal springs 5 and 6 reach plastic mode at $t_1$, but then return back to elastic mode, when springs 3 and 4 begin to yield. 
} \label{fig:moving-set-example-pre-stressed}
\end{figure}

In Section \ref{ssect:ex1} we described the parameters of a simple toy network. It has $m=10$ and one computes ${\rm dim}\, \mathcal{U}=8$ and ${\rm dim}\, \mathcal{V}=2$, therefore the moving set $\mathcal{C}(t)$ in the sweeping process  \eqref{eq:elastoplastic_sp}-\eqref{eq:elastoplastic_moving_set} is a section of a 10-dimensional hyperrectangle by a 2-dimensional plane $\mathcal{V}$. {\color{black}Therefore,  solution $y(t)$ and moving set $\mathcal{C}(t)$ can be represented in $\mathbb{R}^2$ as, respectively,
\begin{gather*}
y_V(t)\in \mathbb{R}^2, \text{ such that } y(t)=Vy_{V}(t),\\
\mathcal{C}_{V}(t) \subset \mathbb{R}^2, \text{ such that } \mathcal{C}(t) = V \mathcal{C}_{V}(t),
\end{gather*}
and illustrated by} Figs. \ref{fig:moving-set-example}a and \ref{fig:moving-set-example-pre-stressed}a. The representations, however, depend on the choice of basis vectors in $\mathcal{V}$, which constitute matrix $V$. In our computations
\[
V=\begin{pmatrix}
-0.247641312342202 &  0.409252171336287\\
-0.247641312342202 &  0.409252171336286\\
-0.404312022261124 & -0.120407185624592\\
-0.404312022261124 & -0.120407185624592\\
 0.497928264380854 & -0.079402340433155\\
 0.497928264380854 & -0.079402340433155\\
 0.116775452487286 &  0.394784775694451\\
 0.116775452487287 &  0.394784775694451\\
 0.116775452487287 &  0.394784775694451\\
 0.116775452487286 &  0.394784775694451
 \end{pmatrix}.
\]
{\color{black} To verify that another basis representation of the subspace $\mathcal{V}$ agrees with the one above we suggest to compare} matrices $VP_V$ that are computed using different bases. {\color{black}The compare procedure} is based on the facts that $VP_V$ is uniquely defined for the space $\mathcal{V}$ and matrix $K$ and that $VP_V$ is independent of the choice of a basis.

\begin{figure}[H]\center
\includegraphics{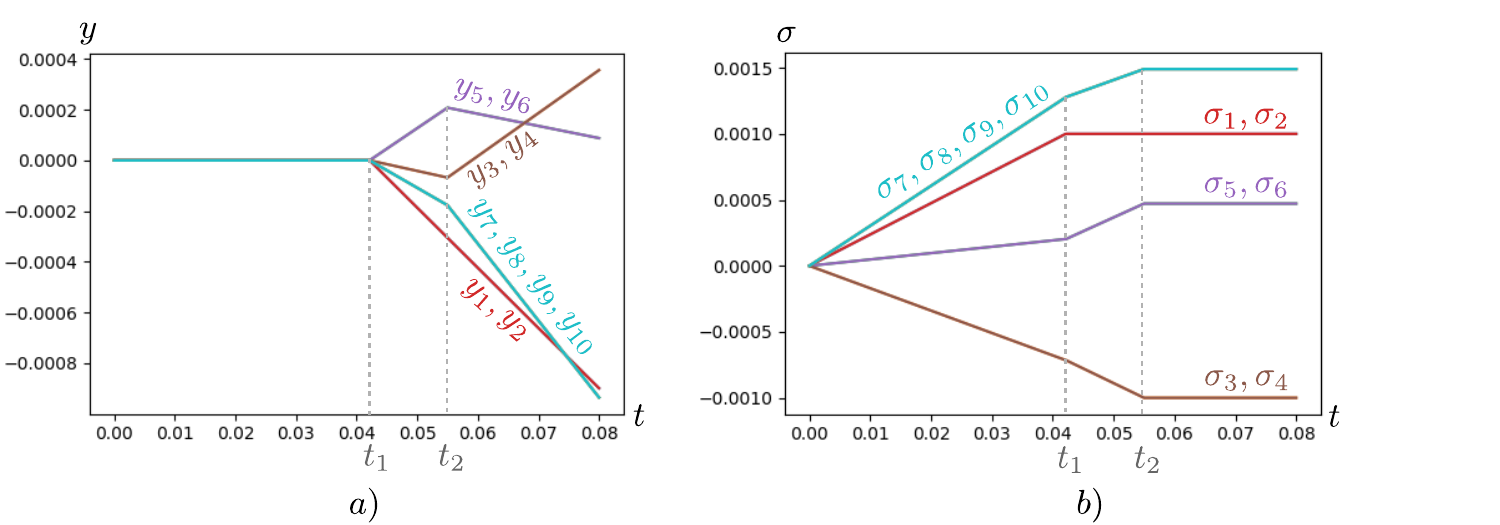}
\caption{
\footnotesize Graphs of the components of the solution $y$ of the sweeping process  \eqref{eq:elastoplastic_sp}-\eqref{eq:elastoplastic_moving_set} set up for Example 1 (a) and the corresponding stress trajectory $\sigma$ (b), connected by the expression \eqref{eq:e_s_from_y}. One can see from (a) that $y$, indeed, remained stationary until it was hit by the edge at time $t_1$, and observe from (b) that time-interval $(0, t_1)$ corresponds to purely elastic deformation of the system. Also from (b) one can see that further application of displacement-controlled loading after time $t_2$ no longer affects the stresses due to perfect plasticity. This corresponds to the solution $y$ being ``trapped'' at the vertex of $\mathcal{C}(t)$ for $t\geqslant t_2$, as shown in Fig. \ref{fig:moving-set-example}a. Also notice, that the signs of the stresses of the vertical springs are different from the signs of the rest of the springs, as the Poisson's ratio dictates.
} \label{fig:example1-graphs}
\end{figure}

\begin{figure}[H]\center
\includegraphics{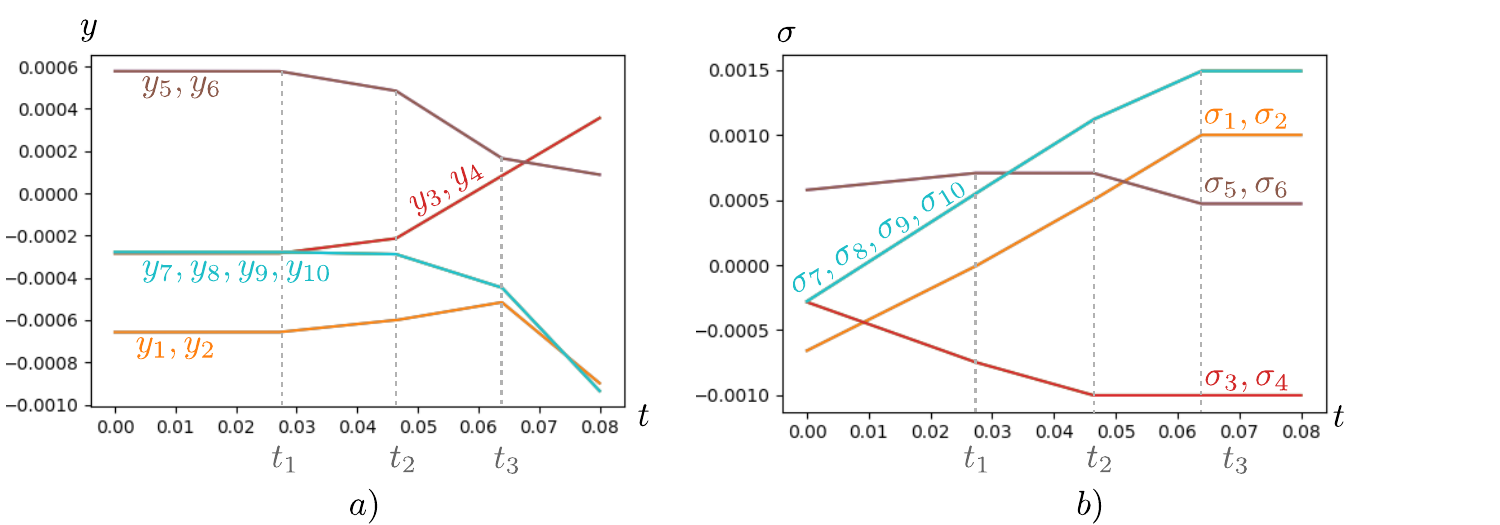}
\vspace{-3mm}
\caption{
\footnotesize
 Graphs of the components of the solution $y$ of the sweeping process \eqref{eq:elastoplastic_sp}-\eqref{eq:elastoplastic_moving_set} (a) and the corresponding stress trajectory $\sigma$ (b), connected by the expression \eqref{eq:e_s_from_y} for the initial condition condition \eqref{eq:alt_initinal_cond}. Notice, how the springs 5 and 6 become loaded and un-loaded depending on the state of the other springs, despite the monotonically increasing displacement loading of the overall system. Also notice, that the final stress values are the same as in Fig. \ref{fig:example1-graphs}, where a different initial condition was taken.
} \label{fig:example1-graphs-pre-stressed}
\end{figure}

In this example we do not apply any stress load, i. e. we consider
\[
f(t)\equiv 0,
\]
and the displacement load is already given by \eqref{eq:example1_constr}. We illustrate the example with two scenarios computed for two different initial conditions:
\begin{itemize}
\item  All springs are relaxed initially:
\[
\sigma_0=0\in \mathbb{R}^{10},
\] 
The trajectory of the sweeping process is shown in Fig. \ref{fig:moving-set-example}a and the corresponding sequence of springs reaching their elasticity limits is shown in Fig. \ref{fig:moving-set-example}b.    The resulting graphs of solution $y$ of \eqref{eq:elastoplastic_sp},\eqref{eq:elastoplastic_moving_set2}  and stress $\sigma$ are shown in Fig. \ref{fig:example1-graphs}.
\item Non-zero initial stresses, given by
{\color{black}
\begin{equation}
\sigma_0=c_0 \,KV\widetilde{\sigma}_0, \qquad \widetilde{\sigma}_0=\begin{pmatrix}1 \\ -1\end{pmatrix}.
\label{eq:alt_initinal_cond}
\end{equation}
}
In this case the solution $y$ takes a different route so that the springs reach plastic deformation in a different order. The result is shown in Fig. \ref{fig:moving-set-example-pre-stressed} and Fig. \ref{fig:example1-graphs-pre-stressed}.
\end{itemize}
\newpage

\section{Event-based Method}
\label{sect:leapfrog}
\subsection{Event-based method for an abstract sweeping process}
A particularly simple case of an abstract sweeping process \eqref{eq:S-sp} is when the polyhedral set $\mathcal{C}(t)$ moves monotonically by translation, i.e. 
\begin{equation}
\mathcal{C}(t) = \mathcal{C}_c + \dot c t
\label{eq:fixedShapeandDirection}
\end{equation}
for some time-independent vector $\dot c$ and a fixed set $\mathcal{C}_c$ of the type \eqref{eq:polyhedral_set}, see e.g. Fig. \ref{fig:sweeping} a. 

In this case we can skip many intermediate steps of the catch-up algorithm and jump directly between the ``events'', i.e. between the instances where the solution meets a new facet of the polyhedron $\mathcal{C}(t)$, for example from $t_0$ to $t_1$ and then to $t_3$ (in terms of Fig. \ref{fig:sweeping}). We refer to this as the ``event-based'' or ``leapfrog'' method and the full algorithm for an abstract sweeping process is given as Algorithm \ref{alg:event_based_abstract}. 
On the $i$-th step of the algorithm we first compute the derivative $\dot z_i$ of the solution {\it relative to the moving set}, which is related to $\dot x_i$ by
\[
\dot z_i:=\dot x_i-\dot c.
\]
To do so, we use the projection \eqref{eq:S-proj-polyhedral}, but accounting only for the currently active constraints of $\mathcal{C}_c$ (they form a so-called {\it tangent cone}, see \cite[p. 67]{Hiriart-UrrutyANDLemarechal}). Then we find the time of the next event by looking for the first intersection of the direction $\dot z_i$ with a new (not currently active) facet of $\mathcal{C}_c$ and get the next position $z_i$ (relative to $\mathcal{C}_c$). Overall, the algorithm finds the values $x_i$ of the solution at the time-moments $t_i$ of the events. The values of $x$ between the events can be found via 
\[
x(t)=x_i+(\dot z_i+\dot c)(t-t_i), \qquad i=\max\{j: t_j\leqslant t\}.
\]

\begin{algorithm}[H]
\SetAlgoLined
\SetKwComment{Comment}{//}{}

 \Comment{Assume that $\mathcal{C}_c$ is given as in \eqref{eq:polyhedral_set}: via matrices $A, A_{eq}$ of dimensions $l\times n, l_{eq}\times n$ respectively and vectors $b,b_{eq}$.}
 $i:=0$\;
 $t_0:=0$\;
 $z_0:=x_0$\;
terminate $:=$ {\bf false}\;
 \Repeat{\rm terminate}{
 \Comment{a tangent cone to $\mathcal{C}_c$, which is also a set of the type \eqref{eq:polyhedral_set}:}
  $T_{\mathcal{C}_c}(z_i):= \left\{x\in \mathbb{R}^n:\begin{array}{l}
  A_{eq}x=b_{eq},\\
  \text{for }j\in\overline{1,l} \text{ such that } (b-Az_i)_j=0 : (A_{jk})_{k\in \overline{1,n}} x\leqslant 0;
  \end{array} \right\}$\;  
  \Comment{use the projection of the type \eqref{eq:S-proj-polyhedral} to find the velocity relative to $\mathcal{C}_c$:}
  $\dot z_i := {\rm proj}^S(-\dot c, T_{\mathcal{C}_c}(z_i))$\;
  \eIf{$\dot z_i \not \approx 0$}{
    \Comment{finding the time until a new event:}
   $\tau_i:=\min\left\{(b-Az_i)_j/(A \dot  z_i)_j: j\in \overline{1,l}:\begin{array}{l}
   (A\dot z_i)_j>0,\\(b-Az_i)_j>0;
   \end{array}
   \right\}$ \;
   $t_{i+1}:=t_i+\tau_i$\;
   \eIf{$t_{i+1}\leqslant T$}{ 
   \Comment{Update the data for the next step and the stresses:}
   $z_{i+1}:=z_i+\dot z_i \tau_i$\;
   $i:=i+1$\;
   $x_i:=z_i+\dot c t_i$\;
   }{\Comment{Further events happen after the interval $[0,T]$}
   terminate $:=$ {\bf true}\;}
   }{\Comment{The solution have stabilized}
   terminate $:=$ {\bf true}\;}
   }  
 \caption{Event-based method for an abstract sweeping process  \eqref{eq:ic},\eqref{eq:S-sp} with a moving set \eqref{eq:fixedShapeandDirection}}
 \label{alg:event_based_abstract}
\end{algorithm}

\subsection{Event-based method for the Lattice Spring Model}
Event-based method can be directly applied to the sweeping process \eqref{eq:elastoplastic_sp},\eqref{eq:elastoplastic_moving_set2} coming form the Lattice Springs Model in the special case when the external force (stress load) is constant and the the displacement load changes at a constant rate, i.e. for all $t\in[0,T]$
\begin{align}
f_c := f(t) \equiv \text{const}, \label{eq:fc}\\
\dot r_c: = \dot r(t) \equiv \text{const}.\label{eq:rc}
\end{align}
These conditions guarantee that the set $\mathcal{C}(t)$ does not change its shape and only moves by translation along the constant direction $G \dot r$:
\begin{equation*}
\mathcal{C}(t)=\mathcal{C}_c + G\dot r_c t, \\
\end{equation*}
where the fixed shape of the set is
\begin{equation}
\mathcal{C}_c :=\left\{x\in \mathbb{R}^m: \begin{array}{c} b^-_{c}\leqslant x \leqslant b^+_c,\\ U^T K x=0 \end{array}\right\}, \qquad \begin{array}{l}b^-_c:= K^{-1}c^-+Gr(0) -F f_c, \\ b^+_c:=K^{-1}c^+ +Gr(0)-F f_c. \end{array}
\label{eq:not-so-moving-set-for-event-based-method}
\end{equation}

Similarly to the event-based method for an abstract sweeping process, we can jump directly between the initial moments of yielding (e.g. for the trajectory in Fig. \ref{fig:moving-set-example-pre-stressed} we go directly from $t=0$  to $t_1$ then to $t_2$ and then to $t_3$). This is especially useful to step over the lengthy initial phase of purely elastic evolution in larger networks (see Sections \ref{sect:example_with_the_hole} and \ref{sect:delaunay} below).
 
Algorithm \ref{alg:event_based1} is an adaptation of Algorithm \ref{alg:event_based_abstract} to the sweeping process \eqref{eq:elastoplastic_sp},\eqref{eq:elastoplastic_moving_set2} in $\mathbb{R}^m$ and, under conditions \eqref{eq:fc}-\eqref{eq:rc}  it computes the stresses $\sigma_i$ at times $t_i$, where each $t_i$ is a time-moment when a new spring starts to yield.

Finally, following Remark \ref{remark:reduction_of_dimensions} we can construct a practical adaptation of Algorithm \ref{alg:event_based_abstract} for the sweeping process in $\mathbb{R}^{{\rm dim}\, \mathcal{V}}$, see Appendix \ref{ssect:event-based method for lattices in V} . Similarly to Algorithm \ref{alg:event_based1}, it requires the same assumptions  \eqref{eq:fc}-\eqref{eq:rc} to hold, and it is significantly faster than Algorithm  \ref{alg:event_based1} for the same lattice, as it deals with much fewer dimensions.

\begin{algorithm}[H]
\SetAlgoLined
\SetKwComment{Comment}{//}{}
\Comment{Given $\sigma_0,m,T, K,U, G, F,r(0), \dot r_c, f_c, b^-_c, b^+_c$}
 $i:=0$\;
 $t_0:=0$\;
 $z_0:= K^{-1}\sigma_0 +Gr(0) - Ff_c$\;
terminate $:=$ {\bf false}\;
 \Repeat{\rm terminate}{
 \Comment{a tangent cone to $\mathcal{C}_c$, which is also a set of the type \eqref{eq:polyhedral_set}:}
  $T_{\mathcal{C}_c}(z_i):= \left\{x\in \mathbb{R}^m:\begin{array}{l}
  U^TKx=0,\\
  \text{for }j\in\overline{1,m} \text{ such that } (b_c^--z_i)_j=0 : x_j\geqslant 0,\\
  \text{for }j\in\overline{1,m} \text{ such that } (b_c^+-z_i)_j=0 : x_j\leqslant 0;
  \end{array} \right\}$\;  
  \Comment{use the projection of the type \eqref{eq:S-proj-polyhedral} to find the rate of change of the stresses:}
  $\dot z_i := {\rm proj}^K(-G\dot r_c, T_{\mathcal{C}_c}(z_i))$\;
  \eIf{$\dot z_i \not \approx 0$}{
    \Comment{finding the time when a new spring starts yielding:}
   $m^-:=\min\left\{(b^-_c-z_i)_j/(\dot z_i)_j: j\in \overline{1,m}:\begin{array}{l}
   (\dot z_i)_j<0,\\(b_i^--z_i)_j<0;
   \end{array}
   \right\}$ \;
   $m^+:=\min\left\{(b^+_c-z_i)_j/(\dot z_i)_j: j\in \overline{1,m}:\begin{array}{l}
   (\dot z_i)_j>0,\\(b_i^+-z_i)_j>0;
   \end{array}
   \right\}$\;   
   $\tau_i:=\min\left(m^-, m^+\right)$\;
   $t_{i+1}:=t_i+\tau_i$\;
   \eIf{$t_{i+1}\leqslant T$}{ 
   \Comment{Update the data for the next step and the stresses:}
   $z_{i+1}:=z_i+\dot z_i \tau_i$\;
   $i:=i+1$\;
   $\sigma_i:=K(z_i-Gr(0)+Ff_c)$\;
   }{\Comment{Further events happen after the interval $[0,T]$}
   terminate $:=$ {\bf true}\;}
   }{\Comment{The stresses have stabilized}
   terminate $:=$ {\bf true}\;}
   }  
 \caption{Event-based method for sweeping process \eqref{eq:elastoplastic_sp}-\eqref{eq:elastoplastic_moving_set2} in $\mathbb{R}^m$}
 \label{alg:event_based1}
\end{algorithm}

\subsection{Example 2: triangular grid with a hole}
\label{sect:example_with_the_hole}

\begin{figure}[h]\center
\includegraphics{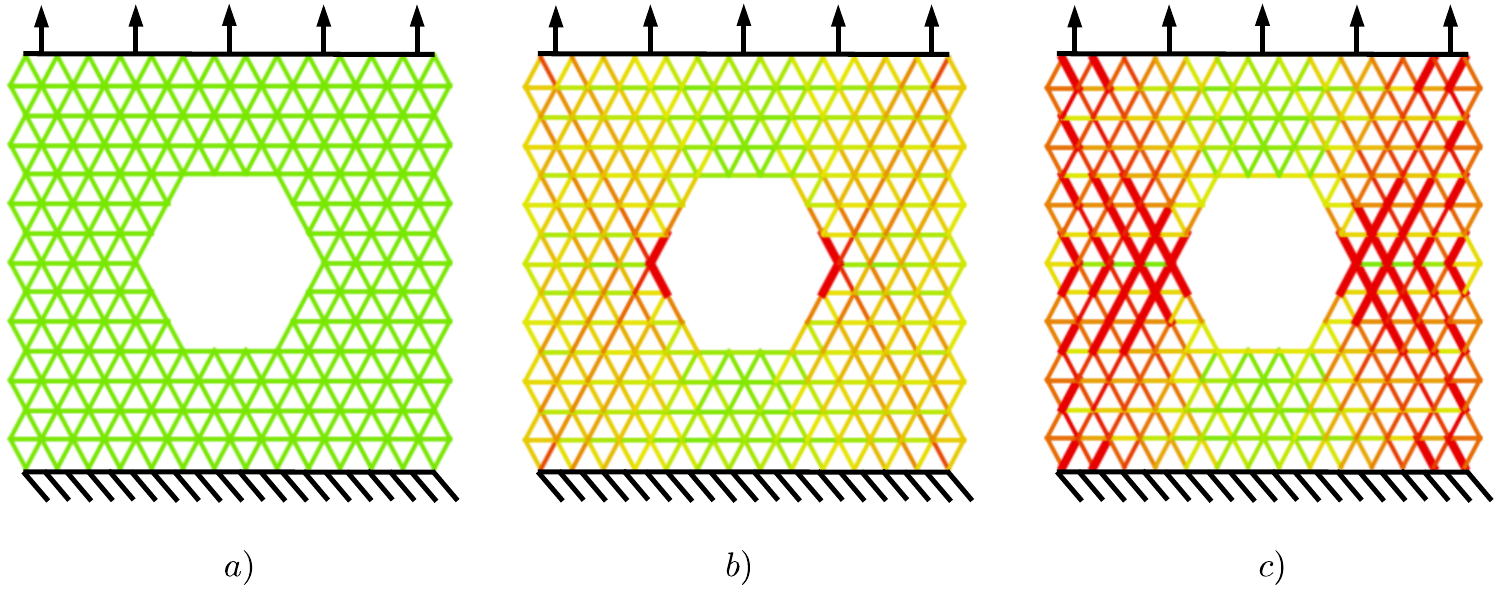}
\caption{
\footnotesize Example 2 at the initial relaxed state (a), at the first yielding event (b) and at the final state (c). A color indicates stress from green (relaxed) to red (maximal stress). {\color{black}The yielding springs are indicated in bold.}
} \label{fig:ex2}
\end{figure}

While the toy example (Example 1) serves as an illustration for the construction of the sweeping process, corresponding to the Lattice Spring Model, we would like to present Example 2, which is a $15 \times 15$ triangular grid  with a hole ($m=496$ springs, $n=198$ nodes, $d=2$), subject to vertical displacement load, see Fig. \ref{fig:ex2}. All the springs are set with stiffness $k_i = 1$ and elastic range $(c_i^-, c_i^+)=(-0.001, 0.001)$. In this example external displacement constraint \eqref{eq:bc} restricts the $x$- and $y$-coordinates of the nodes ($q = 56$ constraints in total) from the top and the bottom of the grid, and the $y$-coordinates of the nodes from the top monotonically increase. Fig. \ref{fig:ex2} shows the key moments of the evolution of stresses in Example 2. {\color{black} The full videos of the simulations via the catch-up and event-based algorithms can be found in the supplemental material \cite{SupplMat}.}

\section{Stress Analysis of Disordered Hyperuniform Networks}
\label{sect:delaunay}

\subsection{Construction of disordered hyperuniform networks}

In this section, we apply the even-based method to analyze the nonlinear mechanical behavior of a class of disordered ``hyperuniform'' networks in two-dimensional Euclidean space under uni-axial loading conditions. In particular, these networks are constructed as the Delaunay triangular network associated with a hyperuniform distribution of points, which by definition possesses vanishing infinite-wavelength density fluctuations. This condition is quantified as the vanishing number variance $\sigma^2_N$ associated with an infinite observation window with linear size $R$, i.e., $\lim_{R\rightarrow \infty} \sigma^2_N(R) = 0$ and equivalently the zero wavenumber limit in the structure factor $\lim_{k\rightarrow 0}S(k) = 0$ (the readers are referred to \cite{hyper00} for details). A unique feature of disordered hyperuniform systems is that they suppress large-scale fluctuations as in a perfect crystal, yet are statistically isotropic and do not possess Bragg peaks as in liquids and glasses. Such unique feature endow these systems with many exotic physical properties, such as large isotropic photonic band gaps \cite{hyper01}, nearly optimal transport properties \cite{hyper02}, and superior mechanical properties \cite{hyper03}.

The hyperuniform point configurations for the construction of the Delaunay networks are numerically generated via the ``collective coordinates'' method \cite{hyper04}, which is essentially a stochastic optimization by setting a target structure factor, i.e., $S(k) = 0$ for $k<K^*$. Starting from a random initial configuration of points, the positions of randomly selected points are continuously perturbed to generate new configurations that gradually converge to the target $S(k)$ (see \cite{hyper04} for details). By tuning the $K^*$ values, one can effectively control the degree of disorder in the generated configurations. In the literature, the parameter $\chi$ (which is $K^*$ normalized by the total degrees of freedom in the system) is typically used to quantify the degree of order, and higher $\chi$ values correspond to more ordered hyperuniform configurations.

The network is contained within a rectangle of width $w = 1$ and height $h=1$ (called the {\it domain}) and it is subject to periodic boundary conditions applied along the x and y edges of the rectangle. We apply uniaxial displacement load along the horizontal direction by increasing $w$, and obtain the stresses within the network using the sweeping process method.

{\color{black}
\subsection{Total stress of the lattice expressed via stresses of the springs}
To characterize the overall response of a lattice to applied load we compute the evolution of the {\it total stress} in the lattice.  As such we borrow the following formula from atomistic simulations (\cite[(6)]{stress_computation_paper}, see also the concept of the system-wide {\it virial stress} \cite[(1.2)]{virial_stress_1}, \cite[p. 6]{virial_stress_2} in which the velocity term ``vanishes for quasistatic deformations'' \cite[p. 246]{Goldhirsch2002}). For a system of pairwise interacting particles the total stress is a $d\times d$ matrix with components
\begin{equation}
\sigma_{total}^{kl}=\frac{1}{2V} \sum_{i,j\in \overline{1,n}}F^k_{ij}r^l_{ij}, \qquad k,l\in \overline{1,d},
\label{eq:total_stress_formula}
\end{equation}
where $F_{ij}$ is the force exerted on particle $i$ by particle $j$, $r_{ij}=r_j-r_i$ with $r_i, r_j$ being the positions of particles $i,j$ respectively, and $V$ is the area or the volume of the domain for 2D or 3D cases respectively.

Since in the current paper we only compute the evolution of the force variables (stresses) and not the spatial variables (displacements and positions), and the elongations of the springs are assumed to be small compared to their lengths by the linearization approach of \eqref{eq:gc}, we will use the reference configuration for the spatial terms in \eqref{eq:total_stress_formula}, namely $r$ and $V$. Furthermore,  the only interactions between particles (nodes) in lattices with periodic boundary conditions are the springs, hence \eqref{eq:total_stress_formula} can be rewritten in terms of springs as
\[
\sigma_{total}^{kl} = \frac{1}{V} \sum_{i\in \overline{1,m}} \left(\sigma_i \mathcal{D}_{ik}\right)\left((\varphi(\xi_0))_i\mathcal{D}_{il}\right),
\] 
where $\mathcal{D}$ is given by \eqref{eq:unit_vectors_along_springs},   the first parenthesis represent the stress of spring $i$ acting on its terminus and the second parenthesis is a vector from the terminus to the origin of the spring. The same product corresponding to the origin of spring $i$ has the same value and cancel out factor $1/2$. Therefore, we can rewrite
\begin{equation}
\sigma_{total} = \frac{1}{V}\, \mathcal{D}^T \, {\rm diag}(\sigma) \,{\rm diag}(\varphi(\xi_0))\, \mathcal{D},
\label{eq:total_stress_via_springs}
\end{equation}
from where we can see that $\sigma_{total}$ is symmetric as it should be.
}
\subsection{Results}

\begin{figure}[H]\center
\includegraphics{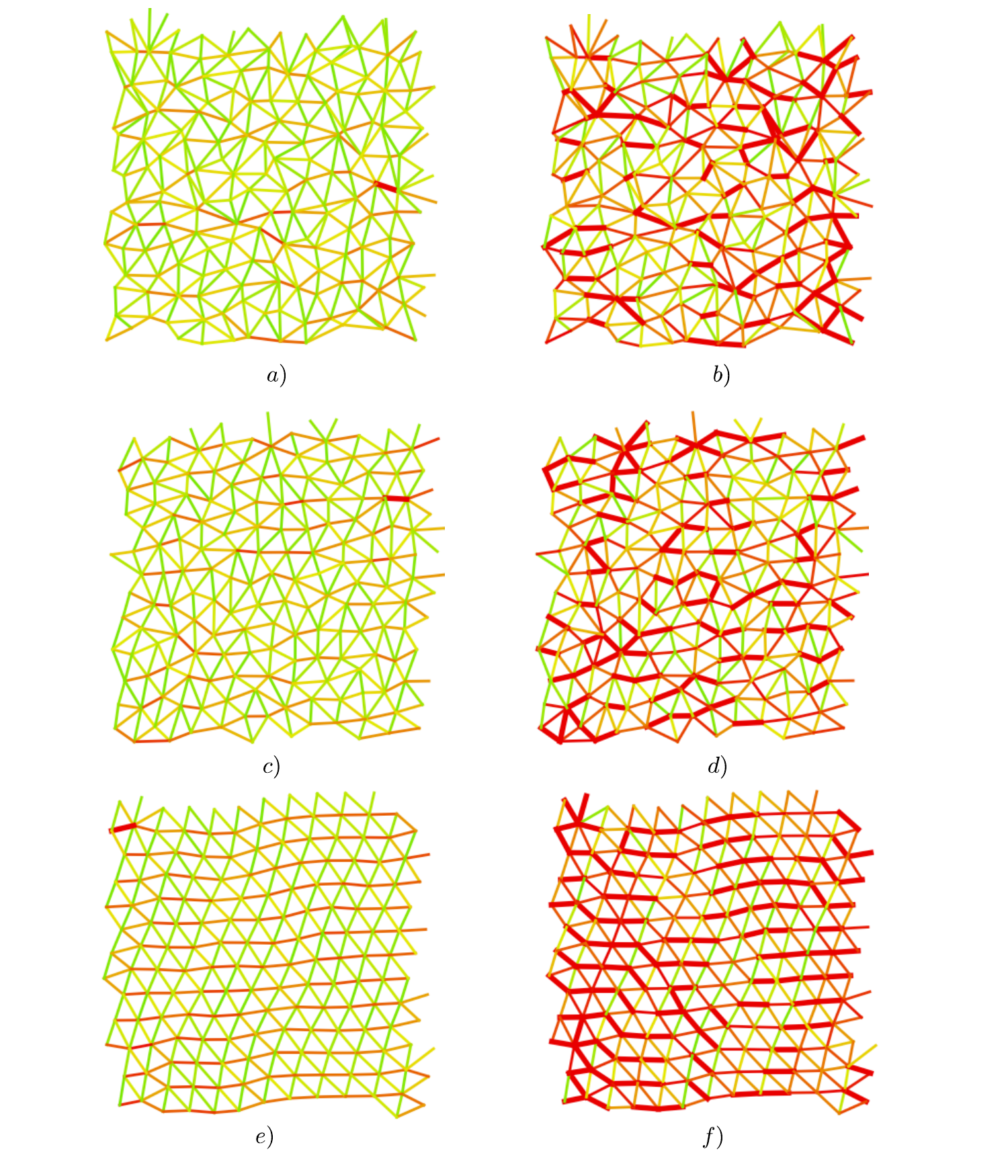}
\caption{
\footnotesize Disordered hyperuniform networks at the first yielding event (a for $\chi=0.3$, c for $\chi=0.4$, e for $\chi = 0.5$) and at the end of the simulation (b for $\chi=0.3$, d for $\chi=0.4$, f for $\chi = 0.5$). A color indicates stress from green (relaxed) to red (maximal stress). {\color{black}The yielding springs are indicated in bold.}
} \label{fig:Delaunay1}
\end{figure}

We explored the hyperuniform networks derived from configurations with $\chi = 0.3, 0.4, 0.5$ (three realizations of each type) and simulated the quasistatic evolution of stresses in each system from the relaxed state under the periodic boundary condition with length $w$ along the horizontal direction monotonically increasing with constant rate from $1$ to $1.04$ (which serves as a horizontal displacement load). The states of the systems (one of each type) at the first yielding event and at the end of the simulation are shown at Figure \ref{fig:Delaunay1} and the full videos of the simulations {\color{black} via the catch-up and the event-based ``leapfrog'' methods can be found in the supplemental material \cite{SupplMat}}. 

\begin{figure}[H]\center
\includegraphics[width=\textwidth,keepaspectratio]{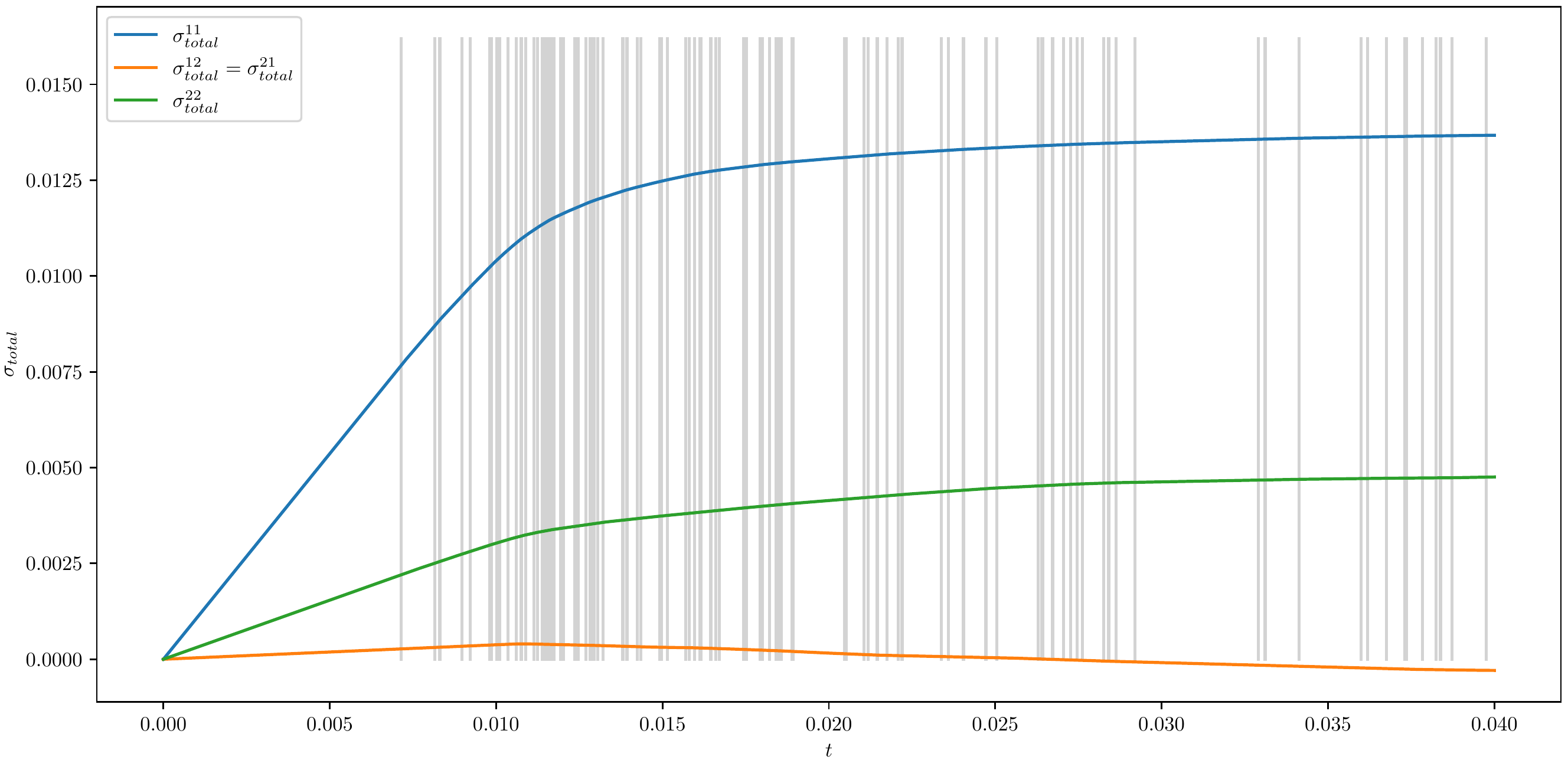}
\caption{
\footnotesize Components of total stress \eqref{eq:total_stress_via_springs} evolving under monotonically increasing horizontal displacement load. Observed increase of component $\sigma^{22}_{total}$ can be explained by Poisson's contraction being restricted by the periodic boundary condition. Gray vertical lines indicate individual yielding events. These graphs correspond to the configuration of Fig. \ref{fig:Delaunay1} a), b), the other configurations behave similarly.
} 
\label{fig:virial-stress-plot}
\end{figure}

The typical behavior is shown in Fig.~\ref{fig:virial-stress-plot}. Since the load is increasing with a constant rate, Fig.~\ref{fig:virial-stress-plot} essentially shows the stress-strain curves obtained in a typical tensile test. In the following analysis, we will focus on the stress-strain behavior in the loading direction.

\begin{figure}[H]\center
\includegraphics[width=\textwidth,keepaspectratio]{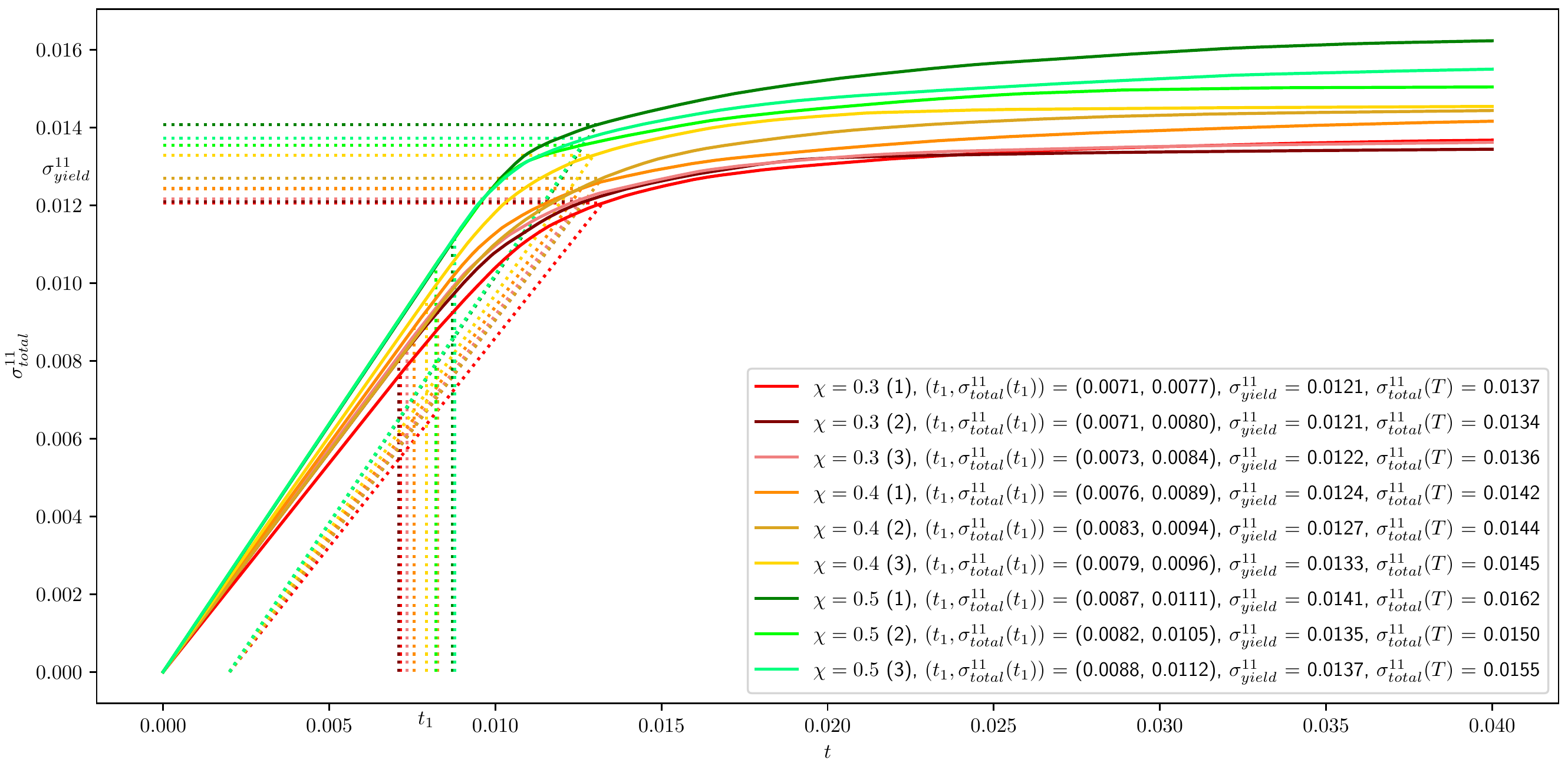}
\caption{
\footnotesize Evolution of component $\sigma^{11}_{total}$, compared for different systems (3 of each of 3 types). We also show the moment of the first yielding event $t_1$ and the corresponding value $\sigma^{11}_{total}(t_1)$, the ``0.2\% yield stresses'' and the stress component $\sigma_{total}^{11}(T)$ at the end of the simulation.
} 
\label{fig:all-virial-stresses}
\end{figure}

Fig. \ref{fig:all-virial-stresses} shows the mechanical behaviors of different network systems (each with three independent realizations). It can be clearly seen that as $\chi$ increases (i.e., the degree of order and hyperuniformity increase), the overall stiffness (i.e., the slope of the linear part of the curve before the first yielding event, indicated with vertical dashed lines), yield strength $\sigma^{11}_{yield}$ and the tensile strength $\sigma^{11}_{total}(T)$ also increase. We note the $\sigma^{11}_{yield}$ is computed following the conventional engineering approach, i.e., we select 0.2\% on the strain axis and construct a straight line with the slope determined by the stiffness, and then the intersection of the constructed line with the stress-strain curve provides the estimated $\sigma^{11}_{yield}$. In turn, $\sigma^{11}_{total}(T)$ is defined as the component of the total stress at the maximal elongation of the system.

To easily compare the properties of the lattices with different values $\chi$ we also refer to Table~\ref{table:dealunay_averaged}, which contains the observed macroscopic values, averaged per each type of the network.

\begin{figure}[H]\center
\includegraphics{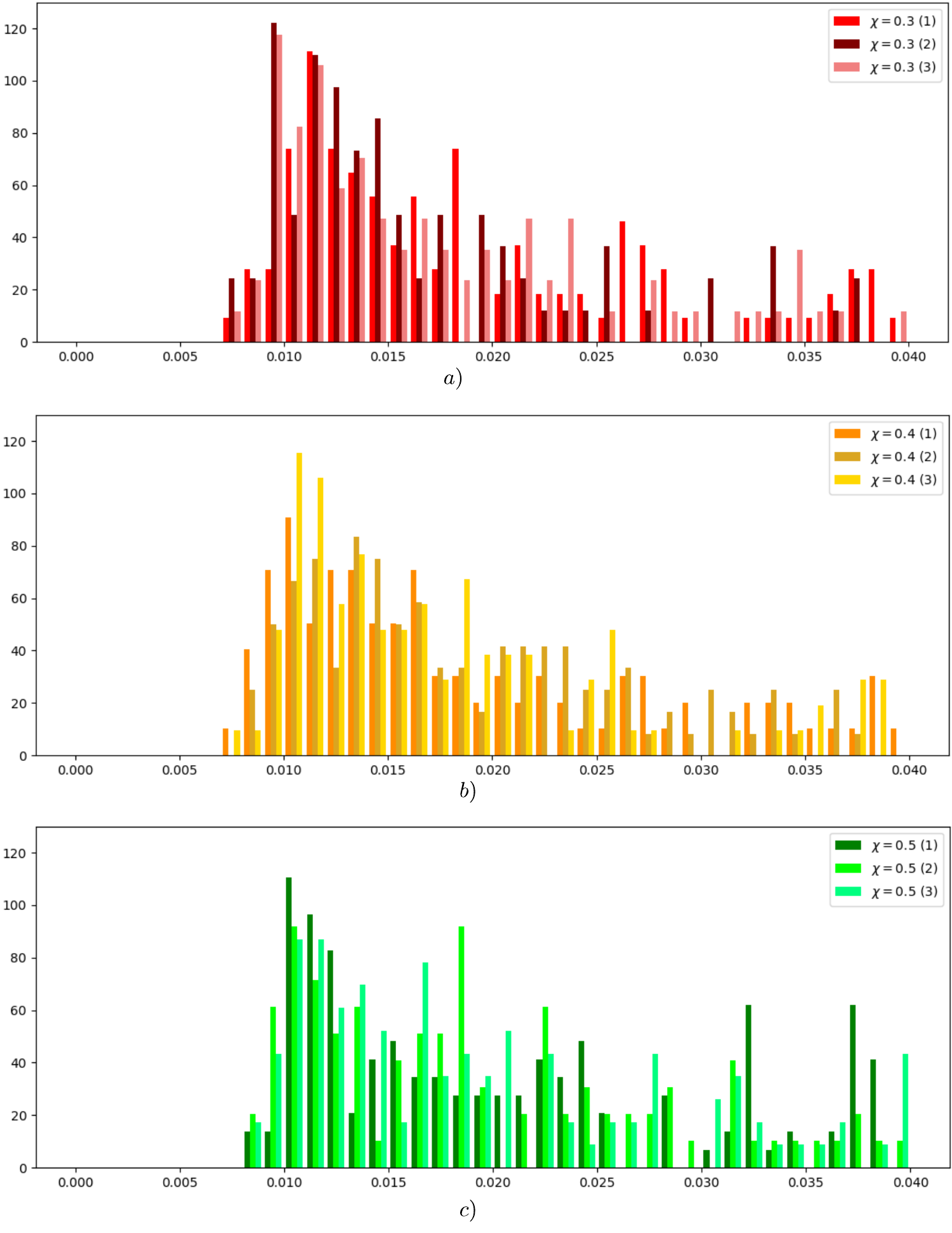}
\caption{
\footnotesize Distributions of yielding events over time during the evolution of the systems. a) corresponds to the systems with $\chi=0.3$, b) with $\chi=0.4$, and c) with $\chi=0.5$.
} 
\label{fig:distributions}
\end{figure}

These macroscopic behaviors can be well explained by the evolution of stress distribution in the systems. In the lattices with small $\chi$ (e.g., 0.3, see Fig. \ref{fig:Delaunay1} a and b) the less uniform distribution of spring (bond) lengths lead to a high degree of stress concentrations, leading to yielding of the springs (i.e., occurrence of the first yielding event) at relatively small overall tensile strains (indicated by the dashed vertical lines in Fig. \ref{fig:all-virial-stresses}). On the other hand, in lattices with high degree of hyperuniformity (with $\chi = 0.5$), the stress distribution is much more uniform, resulting in delayed plasticity and overall increase of stiffness. This analysis is also consistent with the distributions of yielding events in the systems shown in Fig. \ref{fig:distributions}. It can be seen that the yielding events in systems with $\chi = 0.3$ are mainly clustered in early loading stages, while those for $\chi = 0.5$ are spread over the entire loading history. 

\begin{table}[h]
\centering
\begin{tabular}{| c || c | c | c | c |  }
\hline & & & &\\[-3mm]
  $\chi$ & $E$ & $t_1$ & $\sigma^{11}_{yield}$ & $\sigma^{11}_{total}(T)$ \\[1mm] \hline \hline & & & &\\[-3mm]
  0.3 & 1.1174 & 0.0072 & 0.0121 & 0.0136 \\[0.5mm] 
\hline & & & &\\[-3mm]
  0.4 & 1.1731 & 0.0079 & 0.0128 & 0.0144 \\[0.5mm]
\hline & & & &\\[-3mm]
  0.5 & 1.2748 & 0.0086 & 0.0138 & 0.0156 \\[0.5mm]
  \hline
\end{tabular}
\caption{Averaged values of stiffness $E=\frac{\sigma^{11}_{total}(t_1)}{t_1}$, time of the first yielding event $t_1$, yield strength $\sigma^{11}_{yield}$ and tensile strength $\sigma^{11}_{total}(T)$ for each type of networks.}
\label{table:dealunay_averaged}
\end{table}






\section{Conclusions}
\label{sect:conclusions}

In this paper, we connected the mathematical framework of sweeping process to a generic class of Lattice Spring Models with plasticity for nonlinear stress analysis in complex network materials. We started with the governing equations of the quasi-static evolution of the Lattice Spring Model made of elasto-perfectly plastic springs with infinitesimal elongations. We then explicitly constructed a sweeping process to find the evolution of stresses in such Lattice Spring Models using J.-J. Moreau's approach. The sweeping process constructed is of a ``classical'' (unperturbed and convex) type, for which there is a plethora of mathematical research available, ready to be used to uncover properties of elastoplastic Lattice Spring Models. We also provided illustrative examples and established a time-stepping (``catch-up'') and a highly efficient event-based (``leapfrog'') computational frameworks that allow to rigorously track the progression of yielding events in a particular Lattice Spring Model. The utility of our framework has been demonstrated by analyzing the elastoplastic stresses in a novel class of disordered network materials exhibiting the property of hyperuniformity, in which the infinite wave-length density fluctuations associated with the distribution of network nodes are completely suppressed. We find enhanced mechanical properties such as increasing stiffness, yield strength and tensile strength as the degree of hyperuniformity of the material system increases. These results have implications for optimal network material design.


We note that our framework and the leapfrog method can be readily generalized for nonlinear stress analysis in other heterogeneous material systems, such as composites, alloys, porous materials to name a few. The key component in the generalization is the representation of the heterogeneous microstructure of these materials as (ordered) networks (e.g., the triangular network discussed in Sec. 6.3 or face-centered cubic networks in 3D). The nodes of the networks will be grouped according to different material phases they represent, and the constitutive equation governing the springs connecting different phase nodes will be calibrated so that the network system can accurately produce the overall mechanical behavior of the original material. We will explore these generalizations in our future work. {\color{black} Moreover, further extensions of the approach could cover more challenging types of nonlinearities, such as plasticity with softening and large deformations with hypoelasticity.

Apart from a purely applied purpose of computing the evolution of stresses and a theoretical goal of converting the problem into a well-defined sweeping process, in the current paper we presented a working and explicit model which carefully combines the concepts from hysteresis and sweeping process theory with the classical and contemporary studies of framework structures and rigidity.}

\appendix
\section*{Appendix}
\section{The sweeping process of reduced dimension}
\label{sect:appendix_reduced_dimension}
\subsection{Derivation of the sweeping process in $\mathbb{R}^{{\rm dim}\, \mathcal{V}}$}
\label{ssect:catch-up method for lattices in V}
\begin{figure}[h]\center
\includegraphics{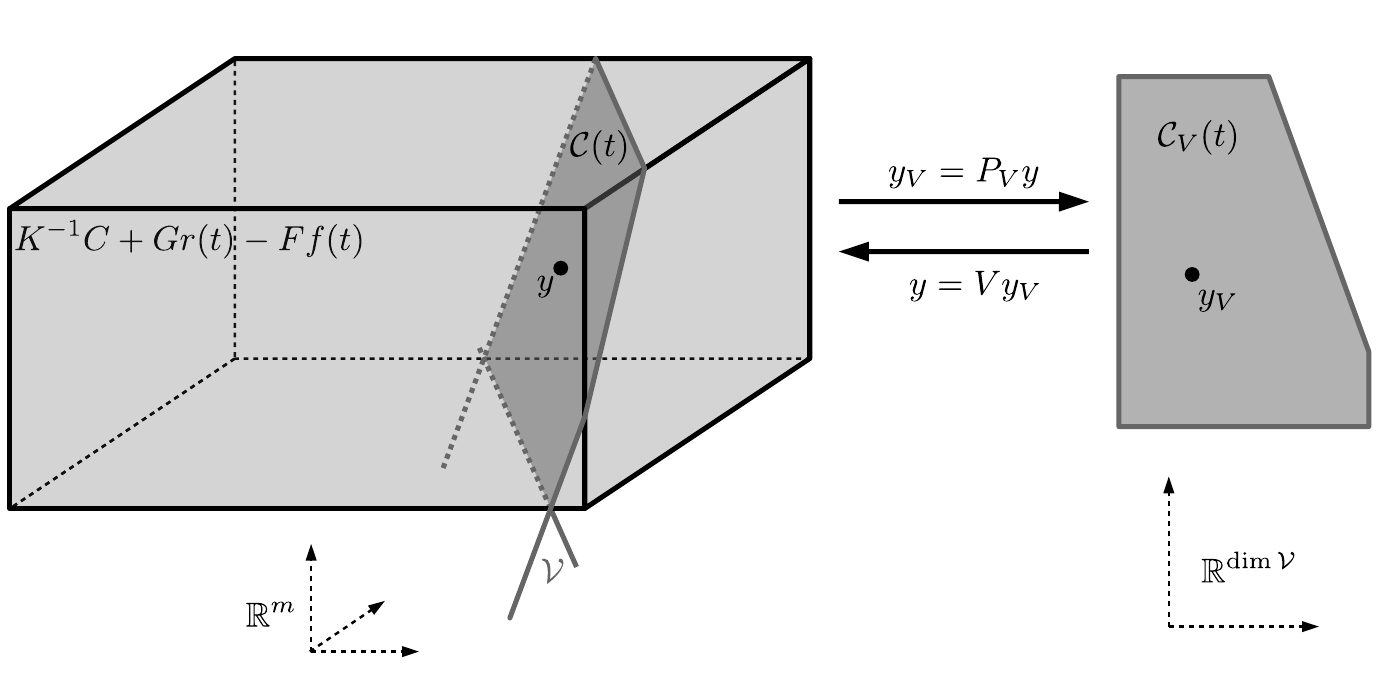}
\caption{
\footnotesize The moving set $\mathcal{C}(t)\subset\mathcal{V}\subset\mathbb{R}^m$ of the sweeping process \eqref{eq:elastoplastic_sp}, the moving set $\mathcal{C}_V(t) \subset \mathbb{R}^{{\rm dim}\, \mathcal{V}}$ of the sweeping process \eqref{eq:elastoplastic_sp_in_V}, and the transformations between them.
} \label{fig:moving-set-with-V-2}
\end{figure}

Here we will provide a sweeping process which is fully equivalent to \eqref{eq:elastoplastic_sp}-\eqref{eq:elastoplastic_moving_set}, but which is significantly cheaper computationally. This happens due to the fact that hyperplane $\mathcal{V}$ in \eqref{eq:elastoplastic_moving_set} is independent of time (hence the equality constraint in \eqref{eq:elastoplastic_moving_set2} is independent of time as well). So, instead of numerically solving the sweeping process in $\mathbb{R}^m$ we are going to formulate an equivalent sweeping process in space $\mathbb{R}^{{\rm dim}\, \mathcal{V}}$, which represents the coordinates of elements of $\mathcal{V}$ in basis $V$, see Fig. \ref{fig:moving-set-with-V-2}. Typically, this significantly reduces the amount of variables which the optimization algorithm has to deal with to compute the projection, and relieves the algorithm from handling the same equality constraint at each time-step of the catch-up algorithm.

In $\mathbb{R}^{{\rm dim}\, \mathcal{V}}$ we define the sweeping process with the unknown $y_V$:
\begin{equation}
\begin{cases}
-\dot y_V \in  N^{S_V}_{\mathcal{C}_V(t)}(y_V),\\
y_V(0)=y_{V0},
\end{cases}
\label{eq:elastoplastic_sp_in_V}
\end{equation}
where 
\begin{equation}
\begin{array}{rll}
S_V &=& V^T K V,\\[0.2cm]
\mathcal{C}_V(t) &=& \left \{x\in \mathbb{R}^{{\rm dim}\, \mathcal{V}}: c^- + KGr(t) -K Ff(t)\leqslant P_V ^T S_V x \leqslant c^+ + KGr(t) - KFf(t)\right\}=\\[2mm]
&=&\left \{x\in \mathbb{R}^{{\rm dim}\, \mathcal{V}}: \begin{array}{c}K^{-1}c^- + Gr(t) -Ff(t)\leqslant (P_V K^{-1})^T S_V x\\[1mm] K^{-1}c^+ + Gr(t) - Ff(t)\geqslant (P_V K^{-1})^T S_V x\end{array}\right\}
,\\[4mm]
y_{V0} &=& P_V K^{-1} \sigma_0 + G_V r(0),\\[1mm]
G_V&=& P_VG\, =\, P_V\,(D_{\xi_0}\varphi)\, R^+.
\end{array}
\label{eq:process_in_V_papameters}
\end{equation}

\begin{Proposition}
\label{prop:equivalent_process}
Let $y:[0,T]\to \mathbb{R}^m$ be a solution to \eqref{eq:elastoplastic_sp},\eqref{eq:elastoplastic_moving_set2}, then $y_V=P_V y$ is a solution to \eqref{eq:elastoplastic_sp_in_V}-\eqref{eq:process_in_V_papameters}. Conversely, for a solution $y_V:[0,T]\to \mathbb{R}^{{\rm dim}\, \mathcal{V}}$ we always have $y=Vy_V$ as a solution to \eqref{eq:elastoplastic_sp},\eqref{eq:elastoplastic_moving_set2}.
\end{Proposition}
\noindent Before proving the equivalence of the sweeping processes, we must show a technical fact on how to represent the projection of the normal cone in the coordinates of basis $V$ in $\mathcal{V}$.
\begin{Lemma}Let $\mathcal{C}\subset \mathcal{V}$ be a nonempty closed convex set and let $x\in\mathcal{C}$. Then 
\[
P_V N^K_\mathcal{C}(x) = N_{P_V \mathcal{C}}^{S_V}(P_V x).
\]
\label{lemma:proj_normal_cone}
\end{Lemma}
\noindent{\bf Proof.} Indeed, 
\begin{multline*}
P_V N^K_\mathcal{C}(x)  =\{P_V y: y \in \mathbb{R}^m, \forall c\in \mathcal{C}: y^T K(c-x)\leqslant 0 \}=\\
=\{P_V y: y \in \mathbb{R}^m, \forall c\in \mathcal{C}: (UP_Uy+ VP_Vy)^T K(c-x)\leqslant 0 \}=\\
=\{P_V y: y \in \mathbb{R}^m, \forall c\in \mathcal{C}: (VP_Vy)^T K(c-x)\leqslant 0 \}=\\
=\{y_V \in \mathbb{R}^{{\rm dim}\, \mathcal{V}}: \forall c\in P_V\mathcal{C}: (Vy_V)^T K(Vc-VP_Vx)\leqslant 0 \}=\\
=\{y_V \in \mathbb{R}^{{\rm dim}\, \mathcal{V}}: \forall c\in P_V\mathcal{C}: y_V^T (V^TKV)(c-P_Vx)\leqslant 0 \}= N_{P_V \mathcal{C}}^{S_V}(P_V x),
\end{multline*}
where the third equality is due to $(c-x)\in C\subset \mathcal{V}$ and the fourth equality is due to $P_V$ being a surjective linear map and $VP_Vx=x$ for any $x\in \mathcal{V}$. $\blacksquare$

\noindent{\bf Proof of Proposition \ref{prop:equivalent_process}.}
Let $y:[0,T]\to \mathbb{R}^m$ be a solution to \eqref{eq:elastoplastic_sp},\eqref{eq:elastoplastic_moving_set2} and put $y_V=P_V y$. By applying $P_V$ to both sides of the inclusion in \eqref{eq:elastoplastic_sp} we get
\[
	- \dot y_V \in P_V N^K_{\mathcal{C}(t)}(y) = N_{P_V\mathcal{C}(t)}^{S_V}(y_V),
\]
where the equality is proven above as Lemma \ref{lemma:proj_normal_cone}. Also notice, that $y(t)\in \mathcal{C}(t)$ for some $t$ if and only if $y_V(t)\in P_V\mathcal{C}(t)$ (since $\mathcal{C}(t)\subset \mathcal{V}$), therefore the normal cone in the above inclusion is well-defined.  Now we show that $P_V\mathcal{C}(t) = \mathcal{C}_V(t)$ with the latter  defined by \eqref{eq:process_in_V_papameters}. Recall that for diagonal matrix with positive coefficients $K$ we have $x\geqslant 0$ if and only if $Kx\geqslant 0$. Then from \eqref{eq:elastoplastic_moving_set2} we have:
\begin{multline*}
P_V\mathcal{C}(t)=P_V\left\{x\in \mathbb{R}^m: \begin{array}{c} K^{-1}c^- + Gr(t)-F f(t)\leqslant x \leqslant K^{-1}c^+  +Gr(t)-F f(t),\\[1mm] U^T K x=0 \end{array}\right\}=\\
=P_V\left\{Vx_V: x_V\in \mathbb{R}^{{\rm dim}\, \mathcal{V}},\forall i\in\overline{1,m}:  \begin{array}{c}\ c^- + KGr(t)- KF f(t)\leqslant e_i^T KVx_V,\\[1mm] c^+  +KGr(t)-KF f(t)\geqslant e_i^T KVx_V\end{array}\right\}=\\
=\left\{x_V\in \mathbb{R}^{{\rm dim}\, \mathcal{V}}: \forall i\in\overline{1,m}: \begin{array}{c} c^- + KGr(t)- KF f(t)\leqslant e_i^T KVx_V\\[1mm] c^+  +KGr(t)-KF f(t)\geqslant e_i^T KVx_V \end{array}\right\},
\end{multline*}
where $e_i$ are the standard basis vectors from $\mathbb{R}^m$.
For each $i\in \overline{1,m}$ recall, that by an equivalent definition of projection (see e.g. \cite[Corollary 5.4]{Brezis2011}),
there is $n_i\in \mathcal{V}$ such that for any $x\in \mathcal{V}$ we have $n_i^TKx=e_i^TKx$, namely, the orthogonal projection $n_i:=VP_Ve_i$ (in sense of the weighted inner product \eqref{eq:S-inner-prod} with $S=K$). Therefore we can continue:
\begin{multline*}
P_V\mathcal{C}(t) = \\
= \left\{x_V\in \mathbb{R}^{{\rm dim}\, \mathcal{V}}: \forall i\in\overline{1,m}: \begin{array}{c} c^- + KGr(t)- KF f(t)\leqslant (VP_Ve_i)^T KVx_V \\[1mm] c^+  +KGr(t)-KF f(t)\geqslant (VP_Ve_i)^T KVx_V \end{array}\right\}=\\
=\left\{x_V\in \mathbb{R}^{{\rm dim}\, \mathcal{V}}:  c^- + KGr(t)- KF f(t)\leqslant (VP_V)^T KVx_V \leqslant c^+  +KGr(t)-KF f(t)\right\}=\\
=\left\{x_V\in \mathbb{R}^{{\rm dim}\, \mathcal{V}}:  c^- + KGr(t)- KF f(t)\leqslant P_V^T S_V x_V \leqslant c^+  +KGr(t)-KF f(t)\right\}=\mathcal{C}_V(t).
\end{multline*}
We have proven the inclusion in \eqref{eq:elastoplastic_sp_in_V}. To prove the expression for the initial condition $y_{V0}$ apply the projection matrix $VP_V$ to both sides of \eqref{eq_y_0}:
\begin{multline*}
VP_Vy_0 = VP_V(K^{-1}\sigma_0+G r(0) - F f(0))=VP_V K^{-1}\sigma_0 +VP_V VG_V r(0)-VP_V Ff(0)=\\=VP_V K^{-1}\sigma_0 +VG_V r(0),
\end{multline*}
where the last equality is due to the facts that $P_V V=I_{{{\rm dim}\, \mathcal{V}}\times {{\rm dim}\, \mathcal{V}}}$ and that $(VP_V)x=0$ for any $x\in\mathcal{U}$, including $x\in {\rm Im}\, F$. 
Apply $P_V$ to both sides and observe that
\[
y_V(0) = P_V y(0)=P_V y_0= P_VVP_V y_0=P_VVP_V K^{-1}\sigma_0 +P_VVG_V r(0)= P_V K^{-1}\sigma_0 + G_Vr(0).
\]
Therefore, $y_V$ is, indeed, the solution of \eqref{eq:elastoplastic_sp_in_V}-\eqref{eq:process_in_V_papameters}.

Conversely, let $y_V:[0,T]\to \mathbb{R}^{{\rm dim}\,\mathcal{V}}$ be a solution to a well-defined process \eqref{eq:elastoplastic_sp_in_V}-\eqref{eq:process_in_V_papameters}. Then $\mathcal{C}(t)=V\mathcal{C}_V(t)$ and $y_0=Vy_{V0}$ define the process \eqref{eq:elastoplastic_sp}-\eqref{eq:elastoplastic_moving_set2} with its own solution $y$. By uniqueness of solution to \eqref{eq:elastoplastic_sp_in_V}-\eqref{eq:process_in_V_papameters} and the previous part of the proof, we must have $P_Vy=y_V$, therefore $Vy_V =VP_Vy=y$ (since $y(t)\in \mathcal{V}$ for all $t$). $\blacksquare$

\subsection{Practical version of the catch-up algorithm in $\mathbb{R}^{{\rm dim}\,\mathcal{V}}$}
{\color{black}Along the lines of Section \ref{ssect:catch-up_elastoplastic}, we rewrite moving set $C_V(t)$ given by \eqref{eq:process_in_V_papameters} in the form \eqref{eq:polyhedral_set}, where}
\begin{equation}
\begin{array}{ll}
A=\begin{pmatrix} P_V^T S_V\\ -P_V ^T S_V\end{pmatrix},&
b(t)=\begin{pmatrix}c^+ + KGr(t)- KF f(t) \\ -\left(c^- + KGr(t)-KF f(t)\right)\end{pmatrix};
\end{array} 
\label{eq:ABinRm_V1}
\end{equation}
or, equivalently, with the same $b$ as in \eqref{eq:ABinRm}
\begin{equation}
\begin{array}{ll}
A=\begin{pmatrix} (P_V K^{-1})^TS_V\\ -(P_V K^{-1})^TS_V\end{pmatrix},&
b(t)=\begin{pmatrix}K^{-1}c^+ + Gr(t)-F f(t) \\ -\left(K^{-1}c^- + Gr(t)-F f(t)\right)\end{pmatrix}
\end{array} 
\label{eq:ABinRm_V2}
\end{equation}
with no equality {\color{black}constraints of \eqref{eq:polyhedral_set}, i.e. without $A_{eq}, b_{eq}$}. The evolution of stress $\sigma$ and elastic elongation $\varepsilon$ can be obtained by 
\[
\varepsilon(t)=Vy_V(t)-G r(t)+ F f(t),\qquad \sigma(t) = K \varepsilon(t).
\]
In turn, Algorithm \ref{alg:catch-up-elastoplastic2} is the adaptation of the catch-up algorithm for the problem \eqref{eq:elastoplastic_sp_in_V}-\eqref{eq:process_in_V_papameters}:

\begin{algorithm}[H]
\SetAlgoLined
\SetKwComment{Comment}{//}{}
\Comment{Given $\sigma_0, K,V, P_V,S_V,G_V,G, F,r,f$ and $\mathcal{C}_V(t)$ via  \eqref{eq:ABinRm_V1} or \eqref{eq:ABinRm_V2}}
\Comment{and a partition $0=t_0<t_1<\dots< t_{k-1}<t_k=T$}
$y_{V0} := P_V K^{-1} \sigma_0 + G_V r(0)$\;
 \For{$i:=0$ \KwTo $k-1$}{
  \Comment{the projection of the type \eqref{eq:S-proj-polyhedral} with constraints \eqref{eq:ABinRm_V1} or \eqref{eq:ABinRm_V2}:} 
  \Comment{no $A_{eq}, b_{eq}$}
  $y_{Vi+1} := {\rm proj}^{S_V}(y_{Vi}, \mathcal{C}_V(t_i))$\;
  \Comment{recover the elastic elongations and the stresses from the solution of the sweeping process:}
  $\varepsilon_{i+1}:=Vy_{Vi+1}-Gr(t_{i+1})+Ff(t_{i+1})$\;
  $\sigma_{i+1}:=K\varepsilon_{i+1}$
   }  
 \caption{Practical catch-up algorithm to compute stresses via the sweeping process \eqref{eq:elastoplastic_sp_in_V}-\eqref{eq:process_in_V_papameters}:}
 \label{alg:catch-up-elastoplastic2}
\end{algorithm}

\subsection{Practical version of the event-based method in $\mathbb{R}^{{\rm dim}\, \mathcal{V}}$}
\label{ssect:event-based method for lattices in V}
Under assumptions  \eqref{eq:fc}-\eqref{eq:rc} moving set of the sweeping process \eqref{eq:elastoplastic_sp_in_V}-\eqref{eq:process_in_V_papameters} takes the form
\begin{equation*}
\mathcal{C}_V(t) = \mathcal{C}_{Vc} + G_V\dot r_c t,
\end{equation*}
with
\[
\mathcal{C}_{Vc}=\left \{x\in \mathbb{R}^v: b^-_c\leqslant W x \leqslant b^+_c\right\}, \qquad W=(P_V K^{-1})^T S_V,
\]
where $b^-_c, b^+_c$ are as in \eqref{eq:not-so-moving-set-for-event-based-method}.
The corresponding adaptation of the event-based method of Section~\ref{sect:leapfrog} is Algorithm~\ref{alg:event_based2}.

\begin{algorithm}[H]
\SetAlgoLined
\SetKwComment{Comment}{//}{}
\Comment{Given $\sigma_0,m,v,T, K,V, P_V,S_V,G_V,G, F,r(0),\dot r_c,f_c,W, b^-_c, b^+_c$}
 $i:=0$\;
 $t_0:=0$\;
 $z_{0}:= P_V K^{-1} \sigma_0 + G_V r(0)$\;
terminate $:=$ {\bf false}\;
 \Repeat{\rm terminate}{
 \Comment{a tangent cone to $\mathcal{C}_{Vc}$, which is also a set of the type \eqref{eq:polyhedral_set}:}
  $T_{\mathcal{C}_{Vc}}(z_{i}):= \left\{x\in \mathbb{R}^v:\begin{array}{l}
  \text{for }j\in\overline{1,m}  \text{ s. t. } (b_c^--Wz_{i})_j=0 : (W_{jk})_{k\in \overline{1,v}} x\geqslant 0,\\
  \text{for }j\in\overline{1,m} \text{ s. t. } (b_c^+- Wz_{i})_j=0 : (W_{jk})_{k\in \overline{1,v}}x\leqslant 0;
  \end{array} \right\}$\;  
  \Comment{use the projection of the type \eqref{eq:S-proj-polyhedral} to find the rate of change of the stresses:}
  $\dot z_{i} := {\rm proj}^{S_V}(-G_V\dot r_c, T_{\mathcal{C}_{Vc}}(z_{i}))$\;
  \eIf{$\dot z_{i} \not \approx 0$}{
    \Comment{finding the time until a new spring starts yielding:}
   $m^-:=\min\left\{(b^-_c-Wz_{i})_j/(W \dot z_{i})_j: j\in \overline{1,m}:\begin{array}{l}
   (W\dot z_{i})_j<0,\\(b_i^--Wz_{i})_j<0;
   \end{array}
   \right\}$ \;
   $m^+:=\min\left\{(b^+_c-Wz_{i})_j/(W\dot z_{i})_j: j\in \overline{1,m}:\begin{array}{l}
   (W\dot z_{i})_j>0,\\(b_i^+-Wz_{i})_j>0;
   \end{array}
   \right\}$\;   
   $\tau_i:=\min\left(m^-, m^+\right)$\;
   $t_{i+1}:=t_i+\tau_i$\;
   \eIf{$t_{i+1}\leqslant T$}{ 
   \Comment{Update the data for the next step and the stresses:}
   $z_{i+1}:=z_{i}+\dot z_{i} \tau_i$\;
   $i:=i+1$\;
   $\sigma_i:=K(Vz_{i}-Gr(0)+Ff_c)$\;
   }{\Comment{Further events happen after the interval $[0,T]$}
   terminate $:=$ {\bf true}\;}
   }{\Comment{The stresses have stabilized}
   terminate $:=$ {\bf true}\;}
   }  
 \caption{Practical event-based method for sweeping process \eqref{eq:elastoplastic_sp_in_V}-\eqref{eq:process_in_V_papameters} in $\mathbb{R}^{{\rm dim}\, \mathcal{V}}$}
 \label{alg:event_based2}
\end{algorithm}

\section*{Acknowledgements}
\hspace{\parindent}
The authors thank Josean Albelo-Cortes for related useful scientific discussions and for the the suggestion to use Moore-Penrose pseudoinverse in particular. {\color{black} Ivan Gudoshnikov thanks Pavel Krej\v{c}\'{i}, Giselle Antunes Monteiro and \v{S}\'{a}rka Ne\v{c}asov\'{a} from IM CAS for helpful scientific discussions.
The authors also thank anonymous referees for providing insightful comments which helped to significantly improve the quality of the paper. } 

Ivan Gudoshnikov (the first author) was successively supported by the NSF Grant CMMI-1916878, the GA\v{C}R project 20-14736S and the project L100192151 funded by the "Programme to support prospective human resources – post Ph.D. candidates" of the Czech Academy of Sciences, and also supported by RVO: 67985840.

Yang Jiao (the second author) is supported by grant NSF CMMI-1916878.

Oleg Makarenkov (the third author) is supported by grant NSF CMMI-1916876. 

\printbibliography

\end{document}